\documentclass[12pt]{amsart}
\usepackage{amsmath}
\usepackage[margin=1in]{geometry}
\usepackage{mathtools}
\usepackage{amssymb} 
\usepackage{fancyhdr}
\usepackage{color}
\usepackage{graphicx}
\usepackage{tabularx}
\usepackage{dcolumn} 
\usepackage{bm}      
\usepackage{array}   
\usepackage{amsfonts}
\usepackage[bookmarks=false,colorlinks,citecolor=red]{hyperref}
\usepackage[numbers, sort]{natbib}
\usepackage{subcaption}
\usepackage{placeins}
\mathtoolsset{showonlyrefs}

\newtheorem{theorem}{Theorem}
\newtheorem{lemma}{Lemma}

\theoremstyle{definition}\newtheorem{remark}{Remark}

\title{Stochastic Solution of Elliptic and Parabolic Boundary Value Problems for the Spectral Fractional Laplacian}
\author{Mamikon Gulian}
\address{Division of Appled Mathematics, Brown University}
\email{\texttt{mamikon\_gulian@brown.edu}}

\author{Guofei Pang}
\address{Division of Appled Mathematics, Brown University}
\email{\texttt{guofei\_pang@brown.edu}}
\begin{document}

\maketitle

\begin{abstract}
We prove and implement stochastic solution (or Feynman-Kac) formulas for boundary value problems involving the spectral fractional Laplacian $(-\Delta_{\Omega,g})^{\alpha/2}$ with nonzero Dirichlet boundary condition $g$. 
The main tools used in the proofs are the abstract Cauchy problem for Feller semigroups together with Balakrishnan's theory of fractional powers. We show the operator
$(-\Delta_{\Omega,g})^{\alpha/2}$  is the generator of an appropriate Feller semigroup for subordinate stopped Brownian motion 
\begin{equation}
X^{\Omega,\alpha}_t = \sqrt{2} B_{T_{\alpha/2}(t) \wedge \tau_\Omega}, 
\quad
\tau_\Omega = \text{inf}\left\{ t \ \big| \ \sqrt{2} B_t \not\in\Omega \right\}
\end{equation}
and obtain a stochastic solution formula
\begin{align}
u(t,x) = \mathbb{E}_{X^{\Omega,\alpha}_0 = x} 
\Bigg[ f\left(X^{\Omega,\alpha}_t\right) \chi_{\tau_\Omega > T_{\alpha/2}(t)} 
& + g\left(X^{\Omega,\alpha}_t\right) \chi_{\tau_\Omega \le T_{\alpha/2}(t)}\Bigg] \\
& + \mathbb{E}_{X^{\Omega,\alpha}_0 = x}  \left[ \int_0^t  r\left(t-s, X^{\Omega,\alpha}_{s}\right)  \chi_{\tau_\Omega > T_{\alpha/2}(s)}  ds \right]
\end{align}
for the fractional heat equation in a bounded domain $\Omega$, 
\begin{align}
\begin{cases}
&\partial_t u(t,x) + (-\Delta_{\Omega,g})^{\alpha/2}u (t,x) 
= r(x,t)\text{ for $t > 0, x \in \Omega$} \\
&u(0,x) = f(x) \text{ for $x \in \Omega$, }\quad u(t,x) = g(x) 
 \text{ for $x \in \partial\Omega$}. 
\end{cases}
\end{align}
Here, $\tau_\Omega$ denotes (Brownian) exit time from $\Omega$, and $T_{\alpha/2}$ is the standard $\alpha/2$-stable subordinator starting at zero. 
We then obtain precise regularity and steady-state convergence properties of the parabolic problem using the eigenfunction expansion of the classical solution, which leads to estimates for the survival probability of subordinate stopped Brownian motion. These results allow us to take $t \rightarrow \infty$ in the parabolic formula to establish a stochastic solution formula 
\begin{equation}
u(x) = \mathbb{E}_{X_0^{\Omega,\alpha} = x} \left[ g\left(X^{\Omega,\alpha}_{T^{-1}_{\alpha/2}(\tau_\Omega)}\right) \right]
+ \mathbb{E}_{X_0^{\Omega,\alpha} = x}  \left[ \int_0^{T^{-1}_{\alpha/2}(\tau_\Omega)}  r\left(X^{\Omega,\alpha}_{s}\right) ds \right],
\end{equation}
for the Dirichlet boundary value problem 
\begin{align}
\begin{cases}
&(-\Delta_{\Omega,g})^{\alpha/2}u (t,x) = r(x) \text{ for $t > 0, x \in \Omega$} \\
&u(x) = g(x) \quad \text{ for $x \in \partial\Omega$}. 
\end{cases}
\end{align}
These stochastic solution formulas for the operator $(-\Delta_{\Omega, g})^{\alpha/2}$ (i.e., in the setting of nonzero boundary conditions) are novel, and allow for efficient, embarrassingly parallel local solution of the above boundary value problems. We discuss the discretization of these formulas, and verify them in dimensions two and three with benchmark examples. We study the effect of the number of path samples and the time step size for path discretization on the accuracy of the solution. 
\end{abstract}

\newpage

\tableofcontents

\section{Introduction.}
The deep connection \cite{meerschaert_sikorskii, metzler2000random, kolokoltsov} between continuous-time random walks (CTRWs) and fractional-order partial differential equations (FPDEs) can be utilized to establish stochastic solution formulas, or Feynman-Kac formulas, for FPDEs. Such formulas forge direct connections between stochastic processes at a microscopic level and an FPDE at a macroscopic level, providing a physical basis for using fractional-order models. At the same time, they provide a simple, embarrassingly parallel method for computing the solution of the FPDE locally at a point without having to generate a grid or otherwise solve for the solution at other points. Monte Carlo method based on stochastic solution formulas scale favorably to high dimensions. In general, various methods can be considered to accelerate numerical implementation of such Feynman-Kac formulas, such as walk-on-spheres \cite{muller1956, cai2013, cai2017, Zhou2016} or quasi-Monte Carlo sampling \cite{moskowitz1996smoothness}. 

In \cite{meerschaert2002stochastic} and \cite{baeumer2001stochastic}, such stochastic solution formulas were studied for a general time-fractional equations involving generators of Feller semigroups. Stochastic solutions of equations involving both first-order and fractional-order time derivatives were studied in \cite{chen2017time}.
As regards fractional Laplacians in bounded domains, stochastic solution formulas for the regional fractional Laplacian were studied \cite{guan2005boundary} and \cite{Guan2006reflected}; formulas for the time-fractional Cauchy problem for this operator have recently been obtained in \cite{toniazzi2018stochastic}. In \cite{chen2012space}, a stochastic solution formula for the time fractional Cauchy problem for Riesz fractional Laplacian with zero (exterior) boundary condition was proven. In \cite{kyprianou2017unbiased} and \cite{shardlow2018walk}, stochastic solution formulas for the Dirichlet problem for the Riesz fractional Laplacian were obtained and walk-on-spheres algorithms were developed. 
For a discussion of these different fractional Laplacians and their connections to stochastic processes, see \cite{big_laplacian}.
Finally, we mention that this line of research is not limited to FPDEs; stochastic connections for nonlocal equations have been studied, e.g., in \cite{du2014nonlocal}. Recently, \cite{du2018stochastic} gave stochastic representations for a general, nonlocal-in-time evolution equation. 

The present article is motivated by recent advances \cite{AntilPfeffererRogovs, Cusimano2017, big_laplacian} in defining the spectral fractional Laplacian with nonzero boundary conditions, including Dirichlet, Neumann, and Robin. These advances have made clear the appropriate definition for the operator, and have established well-posedness of the associated fractional elliptic boundary value problems and fractional parabolic initial-boundary value problems. Therefore, in this article, we develop and implement stochastic solution formulas for such problems for the case of Dirichlet boundary conditions.  

\begin{figure}
\centering
\begin{subfigure}{.5\textwidth}
  \centering
  \includegraphics[width=0.9\linewidth]{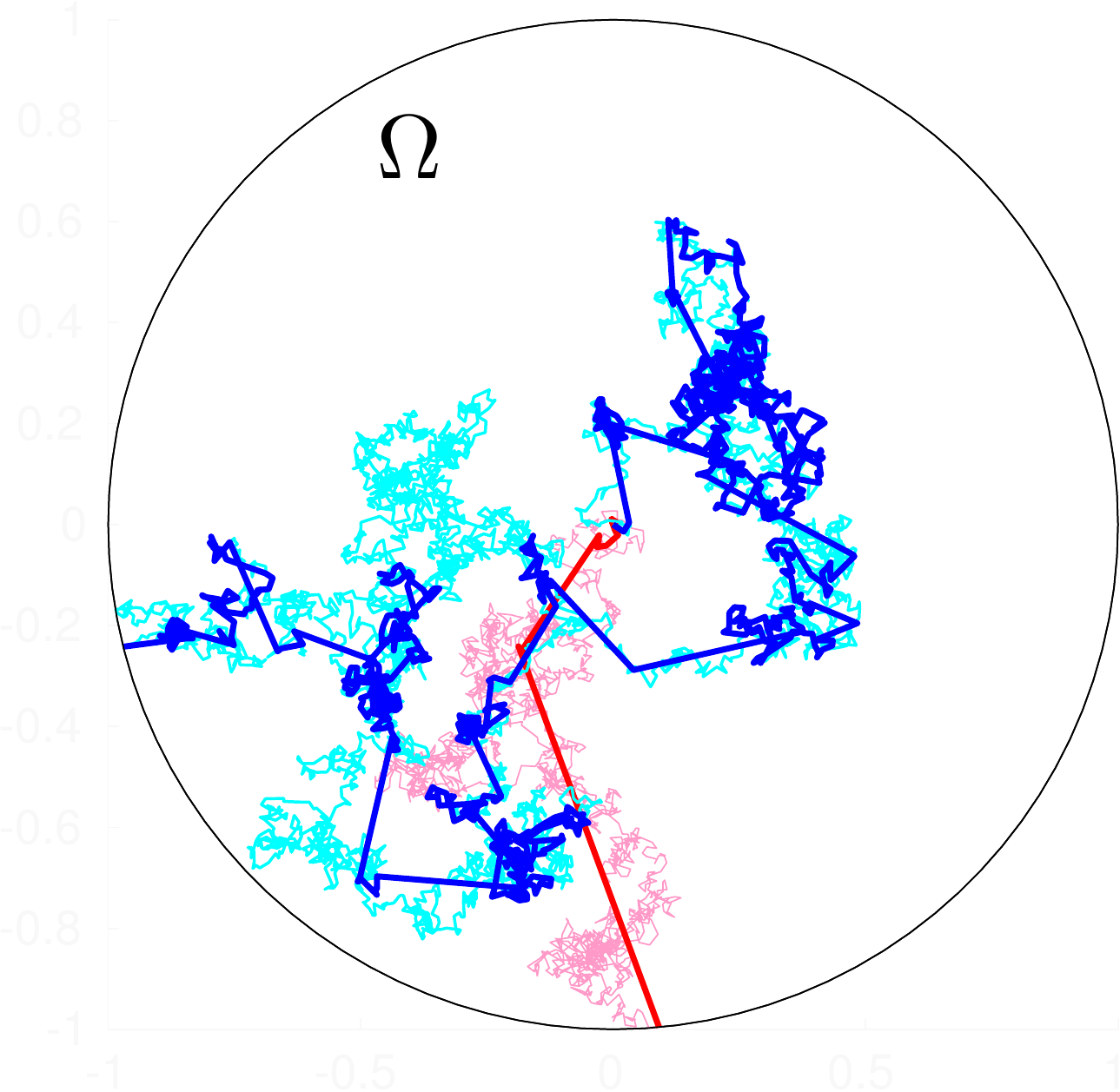}
\end{subfigure}%
\begin{subfigure}{.5\textwidth}
  \centering
  \includegraphics[width=0.9\linewidth]{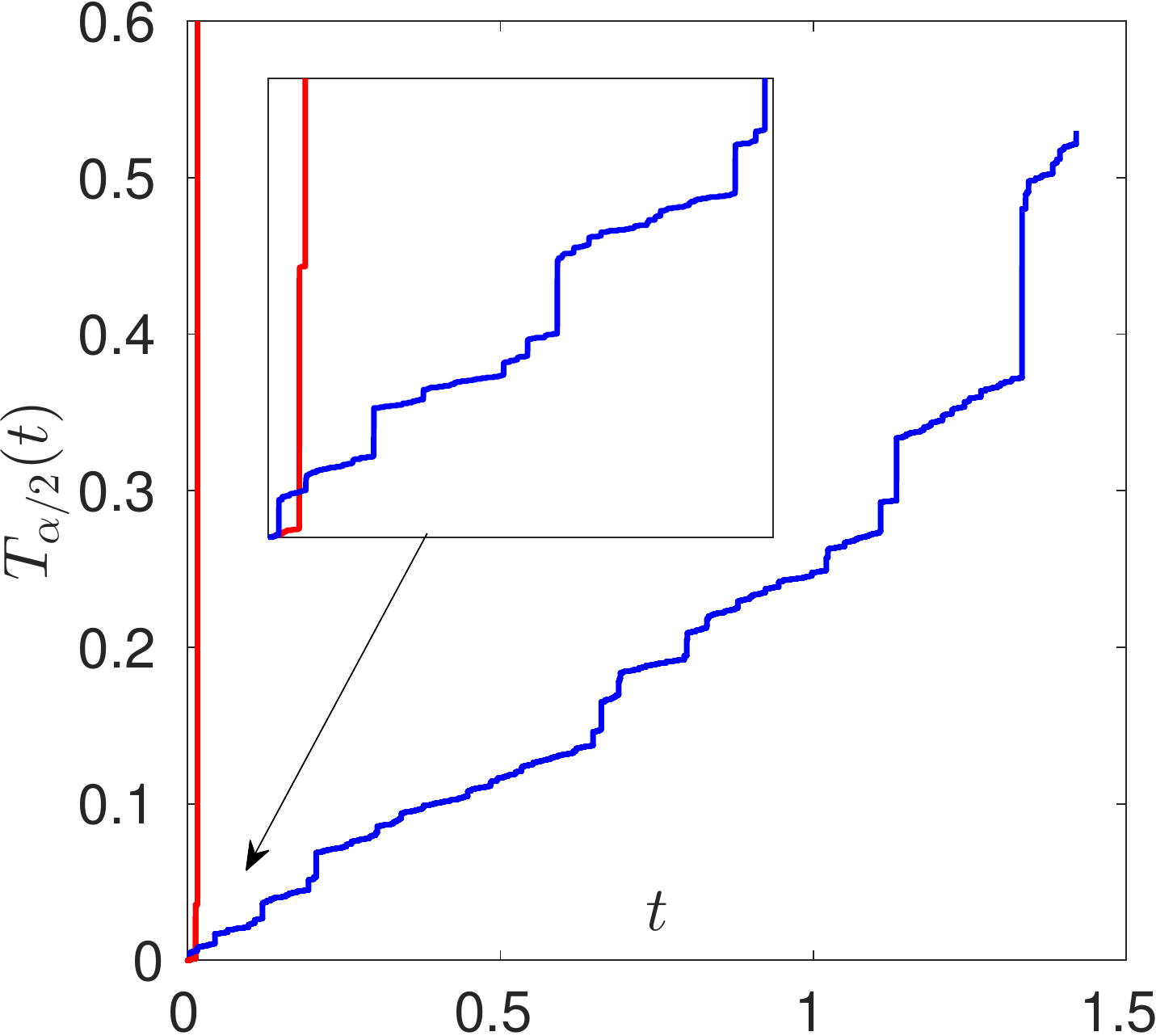}
\end{subfigure}
\caption{\small Illustration of two sample paths of $X^{\Omega,\alpha}_t$ on the unit disk $\Omega$. In the left subfigure, the cyan/pink curves are two samples of $X^{\Omega, 2}_t = X^{2}_{t \wedge \tau_\Omega} = \sqrt{2} B_{t \wedge \tau_\Omega}$.
The right subfigure shows two samples of $T_{\alpha/2}(t)$ that run until they exceed the respective exit times $\tau_\Omega$ of these two samples of $X^{\Omega, 2}_t$. The blue/red curves on the left subfigure are two resulting samples of $X^{\Omega,\alpha}_t = X^{\Omega, 2}_{T_{\alpha/2}(t)}$.  
}
\label{circle_illustration}
\end{figure}

We now outline the main points of the article. In Section \ref{operator_review}, we review the definition and basic properties of the spectral fractional Laplacian with nonzero Dirichlet boundary conditions. In Section \ref{feller}, we review the theory of Feller semigroups, culminating in the abstract Cauchy problem, that will provide the setup for proving stochastic solution formulas. In Section \ref{rn}, we introduce a technique based on subordination and a result of Balakrishnan for obtaining fractional formulas from classical ones, and use this technique to prove the Feynman-Kac formula for the fractional Laplacian $(-\Delta)^{\alpha/2}$ on $\mathbb{R}^d$. This formula itself is well-known and can be shown using the L\'evy-Khinchtine formula or other techniques, but the proof presented here serves an illustration of the method. 
In Section \ref{rn_duhamel}, we extend this result to equations with nonzero forcing term using Duhamel's principle. 
Following this blueprint, in Sections \ref{cauchy} and \ref{cauchy_duhamel} we prove the first main result of the article. 
For the fractional parabolic mixed initial-boundary value problem in $\Omega$, 
\begin{align}
\begin{cases}
&\partial_t u(t,x) + (-\Delta_{\Omega,g})^{\alpha/2}u (t,x) 
= r(x,t)\text{ for $t > 0, x \in \Omega$} \\
&u(0,x) = f(x) \text{ for $x \in \Omega$, }\quad u(t,x) = g(x) 
 \text{ for $x \in \partial\Omega$}, 
\end{cases}
\end{align} 
we obtain the stochastic solution formula
\begin{multline}
u(t,x) = \mathbb{E}_{X^{\Omega,\alpha}_0 = x} 
\left[ f(X^{\Omega,\alpha}_t) \chi_{\tau_\Omega > T_{\alpha/2}(t)}
+ g(X^{\Omega,\alpha}_t) \chi_{\tau_\Omega \le T_{\alpha/2}(t)}\right] \\
+ \mathbb{E}_{X^{\Omega,\alpha_0} = x}  \left[ \int_0^t  r(t-s, X^{\Omega,\alpha}_{s})  \chi_{\tau_\Omega > T_{\alpha/2}(t-s)}  ds \right].
\end{multline}
This formula involves the process $X^{\Omega,\alpha}_t$, which is constructed by first stopping standard isotropic 2-stable motion $X^2_t$ at the boundary of $\Omega$ at sample exit time $\tau_\Omega$, then subordinating by the standard $\alpha/2$-stable subordinator:
\begin{equation}
X^{\Omega,\alpha}_t = X^{\Omega, 2}_{T_{\alpha/2}(t)}, \quad
X^{\Omega, 2}_t = X^{2}_{t \wedge \tau_\Omega}, \quad \tau_\Omega = 
\inf \{ s : X^2_s \not\in \Omega\}. 
\end{equation}  
Standard isotropic 2-stable motion $X^2_t$ is equivalent to $\sqrt{2}B_t$, i.e., Brownian motion scaled by $\sqrt{2}$, so we shall refer to the process $X^{\Omega,\alpha}_t$ as ``subordinate stopped Brownian motion'' throughout this article. The process is illustrated in Figure \ref{circle_illustration}.

For the same problem above, in Section \ref{regularity_section}, we then study a classical solution formula, smoothness for $t > 0$, convergence as $t \rightarrow 0+$ to the initial condition and convergence as $t \rightarrow \infty$ to a steady state which satisfies a fractional elliptic equation. This is the second main result of this article. From this analytical result and the stochastic solution formula for the parabolic problem proved in the previous section, we obtain exponential decay in time of the ``survival probability'' of subordinate stopped Brownian paths $X^{\Omega,\alpha}_t$ -- that is, the probability that such a path has not reached the boundary $\partial \Omega$ (being stopped there) within a given time. These results will are then used to prove the stochastic solution formula for the time-independent problem in the Section \ref{dirichlet}, where, for the fractional elliptic boundary value problem in $\Omega$ 
\begin{align}
\begin{cases}
&(-\Delta_{\Omega,g})^{\alpha/2}u (x) = r(x) \text{ for $t > 0, x \in \Omega$} \\
&u(x) = g(x) \quad \text{ for $x \in \partial\Omega$}, 
\end{cases}
\end{align}
we establish the stochastic solution formula 
\begin{equation}
u(x) = \mathbb{E}_{X^{\Omega,\alpha}_0 = x} \left[ g\left(X^{\Omega,\alpha}_{T^{-1}_{\alpha/2}(\tau_\Omega)}\right) \right]
+ \mathbb{E}_{X^{\Omega,\alpha}_0 = x}  \left[ \int_0^{T^{-1}_{\alpha/2}(\tau_\Omega)}  r(X^{\Omega,\alpha}_{s}) ds \right].
\end{equation}
This is the third main result of the article. Finally, in Section \ref{numerics}, we discuss the implementation of these formulas and provide numerical simulations that verify them for benchmark problems on the $2$-dimensional unit square and $3$-dimensional unit cube. We discuss how to discretize the process $X^{\Omega,\alpha}_t$ and study the convergence of the stochastic solution formula both with the number of paths and the time step size $dt$ used in the discretization.   

There is a rich literature on the potential theory of subordinate killed Brownian motion from the probabilistic perspective \cite{song2003potential, kim2018potential, kim2017boundary, song2008relationship}. In \cite{song2003potential}, it is mentioned that the generator of subordinate killed Brownian motion is the spectral fractional Laplacian with zero Dirichlet boundary conditions. To the best of our knowledge, there have been no studies or reported stochastic solution formulas in this literature for the recently introduced spectral fractional Laplacian $(-\Delta_{\Omega,g})^{\alpha/2}$ with nonzero boundary conditions. This represents a novel contribution of our article. 

\section{The Spectral Fractional Laplacian with Nonzero Dirichlet Boundary Conditions.}
\label{operator_review}
Denoting the Dirichlet eigenvalues and eigenfunctions of $(-\Delta)$ by $\lambda_i$ and $e_i$ , respectively, the spectral fractional Laplacian for zero Dirichlet boundary values is defined as
\begin{align}\label{antil_def}
	(-\Delta_{\Omega,0})^{\alpha/2}u = \sum_{i=1}^\infty \lambda_i^{\alpha/2} (u,e_i)_{L^2(\Omega)} e_i,
\end{align}
for $0 < \alpha < 2$.
See Appendix \ref{appendix} for the basic properties of $(\lambda_i, e_i)$. 
This operator has attracted significant attention, both theoretical and numerical, in the past decade; for discussions of the many works relating to this article, see \cite{big_laplacian, StingaTorrea2010_ExtensionProblemHarnacksInequalityFractionalOperators, nochetto_three, bonito2017approximation} and references therein. However, the generalization of this operator to the case of nonzero boundary conditions was until recently an open problem. 

In 2017, the spectral fractional Laplacian with nonzero boundary conditions (Dirichlet and Neumann) was introduced and studied \cite{AntilPfeffererRogovs, Cusimano2017}. For functions $u$ on $\overline{\Omega}$ satisfying $u|_{\partial\Omega} = g$, Antil, Pfefferer, and Rogovs \cite{AntilPfeffererRogovs} proposed the operator defined by the spectral expansion 
\begin{align}\label{antil_def}
	(-\Delta_{\Omega,g})^{\alpha/2}u = \sum_{i=1}^\infty \left(\lambda_i^{\alpha/2} (u,e_i)_{L^2(\Omega)} - \lambda_i^{\alpha/2-1} \left(u,\frac{\partial e_i}{\partial n}
\right)_{L^2(\partial \Omega)} \right)e_i,
\end{align}
On the other hand, Cusimano, Del Teso, Gerardo-Giorda, and Pagnini \cite{Cusimano2017} defined an operator in the case of nonzero Dirichlet boundary conditions as 
\begin{equation}\label{cusimano_def}
(-\Delta_{\Omega,g})^{\alpha/2} u =
-\frac{1}{\Gamma(-\alpha/2)} \int_0^\infty t^{-\alpha/2 - 1}
(e^{t\Delta_{\Omega,g}}- I) u (x) dt,
\end{equation}
where $I$ denotes the identity operator and $e^{t\Delta_{\Omega,g}}$ is the heat semigroup, i.e., $w(x,t) = e^{t\Delta_{\Omega,g}}u(x)$ is defined as the solution to the problem
\begin{align}
\begin{cases}
&\partial_t w  - \Delta w = 0 \text{ for } x \in \Omega, t > 0  \\
&w(t,x) = g(x) \text{ for } x \in \partial\Omega, t > 0   \\
&w(0,x) = u(x) \text{ for } x \in \Omega. 
\end{cases}
\end{align}

In \cite{big_laplacian}, it was pointed out that the proposed operators \eqref{antil_def} and \eqref{cusimano_def} are essentially the same, as the same characterization in terms of (standard) harmonic lifting was proven by their respective authors of \cite{AntilPfeffererRogovs, Cusimano2017}. More specifically, 
\begin{align}\label{lifting_characterization_for_harbir}
	(-\Delta_{\Omega,g})^{\alpha/2}u &= (-\Delta_{\Omega,0})^{\alpha/2}[u - v],
\end{align}
where $v$ solves 
\begin{equation}
\begin{cases}\label{standard_harmonic_eq}
&-\Delta v = 0 \\ 
&v \big|_{\partial \Omega} = g.
\end{cases}
\end{equation}
Moreover, it was shown in \cite{big_laplacian} that the same operator can be obtained by first taking the spectral power of the inverse fractional Laplacian 
\begin{equation}
(-\Delta)^{-\beta} u = \sum_{i=1}^\infty \lambda^{-\beta} (u,e_i)_{ L^2(\Omega) } e_i, \quad
\beta > 0.
\end{equation}
then defining $(-\Delta_{\Omega,g})^{\alpha/2} =  (-\Delta)^{\alpha/2-1} (-\Delta) u$, noting that $\alpha/2 - 1 < 0$,  
and using the classical formula
\begin{equation}\label{E:series_inhomo_standard_laplacian}
(-\Delta) u = \sum_i^\infty  \left( \lambda_i (u,e_i)_{L^2(\Omega)} - \left( u ,  \frac{\partial e_i}{\partial n}  \right)_{L^2(\partial\Omega)} \right)e_i, 
\end{equation} 
which is valid for any $C^2(\Omega)$ function $u$, regardless of boundary values.
Thus, there is no ambiguity in the use of the symbol $(-\Delta_{\Omega,g})^{\alpha/2}$.

The solutions of boundary value problems involving the operator $(-\Delta_{\Omega,g})^{\alpha/2}$ have corresponding harmonic lifting characterizations as well. 
From \cite{AntilPfeffererRogovs}, the solution of
\begin{align}
\begin{split}
\begin{cases}
&(-\Delta_{\Omega,g})^{\alpha/2} u(x) = f(x), \hspace{15pt} x \in \Omega, \ \ \ \alpha \in (0,2), \\
&u(x) = g(x), \hspace{15pt} x \in \partial \Omega.
\end{cases}
\end{split}
\end{align}
can be written as
\begin{align}
\label{harbir_splitting}
u(x) = w(x) + v(x),
\end{align}
where $v$ again solves \eqref{standard_harmonic_eq} and $w$ solves the problem
\begin{align}
\label{harbir_splitting_zero_bc}
\begin{split}
\begin{cases}
&(-\Delta_{\Omega,0})^{\alpha/2} w = f \hspace{20pt} \text{in} \ \Omega, \\
&w \big|_{\partial \Omega} = 0 \hspace{20pt} \text{on} \ \partial \Omega,
\end{cases}
\end{split}
\end{align}
provided the problem \eqref{standard_harmonic_eq} for $v$ has a classical solution.
In fact, \cite{AntilPfeffererRogovs} showed that this decomposition holds even if the problem \eqref{standard_harmonic_eq} admits a solution in the \emph{very weak variational sense}. However, 
in the present article, we will only deal with classical solutions. 
In the following sections, we will prove Feynman-Kac formulas for parabolic and elliptic problems posed with $(-\Delta_{\Omega,g})^{\alpha/2}$ and Dirichlet boundary conditions.

\section{Markov Semigroups, Generators, Feller Processes, and the Abstract Cauchy Problem.}
\label{feller}
The following is setup of Feller semigroup theory needed for the main result, distilled from the first two chapters of Mark Freidlin's book \cite{freidlin}.  Let $\Omega$ be a locally compact subset of $\mathbb{R}^d$; this includes open and closed subsets of $\mathbb{R}^d$. Let $X_t : \mathbb{R}^+ \rightarrow \Omega$ be a Markov process in $\Omega$. Define the one-parameter family of operators $\mathcal{T}_t$ for $t \ge 0$, acting on functions $f : \Omega \rightarrow \mathbb{R}$, as follows:
\begin{equation}
\mathcal{T}_t f : \Omega \rightarrow \mathbb{R}, \quad \mathcal{T}_t f (x) = 
\mathbb{E}_{X_0 = x} \left[ f(X_t) \right]
\end{equation}
Then
$\mathcal{T}_t$ satisfies the semigroup properties in $t$, and is linear in $f$:
\begin{align}
\label{identity}
&\mathcal{T}_0 f (x) = f(x) \\
\label{semigroup}
&\mathcal{T}_t \mathcal{T}_s f (x) = \mathcal{T}_{t+s} f(x) \\
\label{linear}
&\mathcal{T}_t [f + \tilde{f}] (x) = \mathcal{T}_t f (x) + \mathcal{T}_t \tilde{f} (x)
\end{align}
Here $t,s > 0$.
Property \eqref{identity} is obvious. Property \eqref{semigroup} follows from the Markov property of $X_t$. Property \eqref{linear} is a property of mathematical expectation. If, in addition, $X_t$ is a \emph{Feller Process}, then the semigroup is \emph{contractive} and \emph{right continuous}:
\begin{align}
\tag{F1}
\label{contraction}
&\| \mathcal{T}_t f \|_{C_0(\Omega)} \le \| f \|_{C_0(\Omega)}  \\
\tag{F2}
\label{rightcontinuous}
&\lim_{t \rightarrow 0}\| \mathcal{T}_t f - f \|_{C_0(\Omega)} = 0  
\end{align}
Here, $C_0(\Omega)$ denotes the space of continuous functions vanishing at infinity, i.e., real-valued continuous functions $f$ on $\Omega$ with the property that for every $\epsilon>0$, there exists a compact subset $K \subset \Omega$ such that $|f(x)| < \epsilon$ for all $x$ outside of $K$. The norm of $C_0(\Omega)$ is the sup-norm on $\Omega$:
\begin{equation}
\| f \|_{C_0(\Omega)} = \sup_{x \in \Omega} |f(x)|.
\end{equation}
Thus, $\{\mathcal{T}_t\}$ is a family of (uniformly) bounded linear operators from
$C_0(\Omega)$ to itself. We refer to $\{\mathcal{T}_t f\}$ as the \emph{Feller semigroup} 
of the process $X_t$. We note that L\'evy Processes are Feller (Theorem 3.1.9, \cite{applebaum}).

Associated to every Markov semigroup $\mathcal{T}_t$ is an \emph{infinitesimal generator} or \emph{infinitesimal operator} $\mathcal{A}$, again acting on functions $f :\Omega \rightarrow \mathbb{R}$, defined by 
\begin{equation}
\label{generator1}
\mathcal{A} f :\Omega \rightarrow \mathbb{R}, \quad
\mathcal{A} f (x) = \lim_{t \rightarrow 0+} \frac{\mathcal{T}_t f (x) - f(x)}{h},
\end{equation}
where the limit is taken in $C_0(\Omega)$, i.e., 
\begin{equation}
\label{generator2}
\left\|  \frac{\mathcal{T}_t f (x) - f(x)}{h} - \mathcal{A} f (x) \right\|_{C_0(\Omega)}
\rightarrow 0 \text{ as } t \rightarrow 0+.
\end{equation}
The domain $\text{Dom}(\mathcal{A})$ is the space of all functions $f$ where the above limit exists. 

In some works, the semigroup is jumped over, and one speaks of $\mathcal{A}$ as the \emph{infinitesimal generator} of the process $X_t$ without introducing the notation $\mathcal{T}_t$. However, the semigroup is the key to proving stochastic solution formulas. It features in the following important result:
\begin{lemma}
\label{lemma_acp}
Let $X_t$ be a Feller Process.
The function $u(t,x) = \mathcal{T}_t f(x)$ solves the \textbf{abstract Cauchy problem},
\begin{align}
\label{acp}
\begin{cases}
&\frac{d}{dt} u(t,x) = \mathcal{A} u(t,x) \text{ for $t > 0$}, \\
&\lim_{t \rightarrow 0} u(t,x) = f(x), f \in \text{\textup{Dom}}(\mathcal{A}).
\end{cases}
\end{align}
It is the unique solution in the class of functions that grow at most exponentially, i.e., those solutions $u(t,x)$ for which there exist constants $P, k$ such that 
$ \| u(t, \cdot) \|_{C_0(\Omega)} \le P e^{k t}   . $
\end{lemma}

\begin{remark}
We understand the statement of problem \eqref{acp} to mean: the limit in $C_0(\Omega)$ of $\frac{1}{h} \left[ u(t+h,\cdot) - u(t,\cdot) \right]$ as $h \rightarrow 0$ exists and is equal to $\mathcal{A} u(\cdot,x)$. This convergence in $C_0(\Omega)$ (i.e., uniform convergence) of the limit quotient defining $\frac{d}{dt} u(t,x)$ is a stronger requirement than is typically understood when seeking a classical solution of a differential equation. It must be met to apply the uniqueness result in the lemma.
\end{remark}

This result is transcribed from the book of Mark Freidlin \cite{freidlin}.  For full details of this theory, the reader should consult that text and references therein.

{
}

\section{Subordination, Balakrishnan's Theorem, and the fractional Feynman-Kac formula on $\mathbb{R}^d$.}\label{rn}

In this section, we introduce background material on subordination of processes together with results of Balakrishnan that connect this to concept to fractional calculus. Rather than listing these results in isolation, we introduce them as steps used in a proof of the Feynman-Kac formula for the fractional Laplacian $(-\Delta)^{\alpha/2}$ in $\mathbb{R}^d$ by leveraging  the Feynman-Kac formula for the classical Laplacian $(-\Delta)$. This will mirror our approach for obtaining our main results with the operator $(-\Delta_{\Omega,g})^{\alpha/2}$ in later sections. The formula itself \eqref{levy_semigroup} in $\mathbb{R}^d$ is well-known; for example it appears in \cite{chen2012space} where it is obtained by appealing to Dirichlet form theory \cite{fukushima2010dirichlet}. 

We consider the following problem
\begin{align}
\label{problem_rn}
\begin{cases}
&\partial_t u(t,x) = -(-\Delta)^{\alpha/2}u (t,x) \text{ for $t > 0$} \\
&u(0,x) = f(x).
\end{cases}
\end{align}
We prove the following 
\begin{theorem}
\label{theorem_rn}
Let $u(t,x)$ be a solution to \eqref{problem_rn} such that 
\begin{align}
\label{rn_regularity_assumption_1}
&\text{for every $t \ge 0$, $u(t,\cdot) \in C^2_x(\mathbb{R}^d)$, and}\\
\label{rn_regularity_assumption_2}
&\text{$\partial_t u(\cdot,\cdot)$ is uniformly continuous in spacetime $[0,\infty) \times \mathbb{R}^d$.}
\end{align}
Then 
\begin{equation}
\label{levy_semigroup}
u(t,x) = \mathcal{T}^\alpha_t f (x) := \mathbb{E}_{X^\alpha_0 = x} \left[ f(X^\alpha_t) \right]
\end{equation}
where $X^\alpha_t$ is standard isotropic $\alpha$-stable L\'evy motion in $\mathbb{R}^d$. 
\end{theorem}

To prove the theorem, we will show that the assumption \eqref{rn_regularity_assumption_1} implies that $u(t,\cdot) \in \text{Dom}(\mathcal{A^\alpha})$ and 
\begin{equation}
\label{equivalence_alpha}
-(-\Delta)^{\alpha/2} u = \mathcal{A}^\alpha u, 
\end{equation}
where $\mathcal{A}^\alpha$ is the generator of $\alpha$-stable L\'evy motion $X^\alpha_t$:
\begin{equation}
\label{generator_levy}
\mathcal{A}^\alpha f (x) = 
\lim_{t \rightarrow 0+} \frac{\mathcal{T}^\alpha_t f (x) - f(x)}{h}.
\end{equation} 
Since $X_t^\alpha$ is a Feller process, this result \eqref{equivalence_alpha} together with equation \eqref{problem_rn} and the assumption \eqref{rn_regularity_assumption_2} will imply that $u$ satisfies the abstract Cauchy problem \eqref{acp} of the operator $\mathcal{A}^\alpha$ . By the uniqueness asserted in Lemma \ref{lemma_acp}, $u(t,x)$ must then be given by the semigroup \eqref{levy_semigroup}, which is Theorem \ref{theorem_rn}.

Now, the result \eqref{equivalence_alpha} is classical in the case $\alpha = 2$, where it can be shown using Ito's rule as follows \cite{freidlin}. Noting that $X^2_t$ is scaled Brownian motion $\sqrt{2} B_t$, for $u \in C^2(\mathbb{R}^d)$ we have
\begin{equation}
u(X^2_t) - u(x) = \int_0^t \sum_{i=1}^n \frac{\partial u}{\partial x_i} (X_s) (dX^2_s)_i + \int_0^t \Delta u(X^2_s) ds.
\end{equation}
Taking expectations, 
\begin{equation}
\mathbb{E} [ u(X^2_t) ] - u(x) =  \int_0^t \Delta u(X^2_s) ds.
\end{equation}
Dividing by $t$ and taking $t \rightarrow 0$ gives
\begin{equation}
\mathcal{A}^2 u(x) = \Delta u(x), \quad (u \in C^2)
\label{ito_rule}
\end{equation}
This suggests one way to prove Theorem \ref{theorem_rn} is to prove \eqref{equivalence_alpha} using some sort of Ito rule for L\'evy processes. However, we will instead use \emph{subordination} to upgrade the classical $\alpha = 2$ result \eqref{ito_rule} to general $\alpha < 2$. 

We shall use the fact that $\alpha$-stable L\'evy motion $X^\alpha_t$ is equivalent to \emph{subordinate} Brownian motion, i.e., 
\begin{equation}
\label{subordination_process}
X^\alpha_t = X^{2}_{T_{\alpha/2}(t)}, 
\end{equation}
where $T_\beta(t)$ is the standard $\beta$-stable subordinator starting at zero, an increasing L\'evy process which can most easily be described as having a 
probability density function $h_\beta(t)$ with Laplace transform 
\begin{equation}
\label{subordinator_laplace}
\mathcal{L} \left\{ h_\beta \right\} (s) = e^{- s^\beta}
\end{equation}
See (see Meerschaert and Sikorskii \cite{meerschaert_sikorskii} p. 156 or Sato \cite{sato} p. 198).
Then from \eqref{levy_semigroup} and \eqref{subordination_process}, a conditioning argument yields
\begin{align}
\label{subordinated_semigroups}
\begin{split}
\mathcal{T}^\alpha_t f(x) &= \mathbb{E}_{X^{2}_0 = x} \left[ f(X^2_{T_{\alpha/2}(t)}) \right] \\
 &= \mathbb{E}_{X^{2}_0 = x} \left[  \int \mathbb{P} \{T_{\alpha/2} = s \} f(X^2_{s})  ds \right] \\
 &= \int \mathbb{P}\{T_{\alpha/2} = s \} \mathbb{E}_{X^{2}_0 = x} \left[   f(X^2_{s})  \right]  ds\\
 &= \int \mathbb{P}\{T_{\alpha/2}  = s \} \mathcal{T}^2_s f(x)  ds \\
 & = \int h_{\alpha/2}(s) \mathcal{T}^2_s f(x)  ds.
\end{split}
\end{align}

We require the following theorem of Balakrishnan \cite{balakrishnan}. The notion of \emph{equicontinuous semigroup of class $C_0$} is a semigroup $\mathcal{T}_t$ that satisfies \eqref{identity} and \eqref{semigroup}, as well as the condition 
\begin{equation}
\lim_{t \rightarrow t_0} \| \mathcal{T}_t f - \mathcal{T}_{t_0} f \|_{C_0} = 0
\end{equation}
We note that a Feller semigroup is equicontinuous of class $C_0$, as for $h > 0$, the contraction property \eqref{contraction} gives
\begin{align}
\| \mathcal{T}_{t+h} f  - \mathcal{T}_{t} f \|_{C_0} &\le 
\| \mathcal{T}_t [\mathcal{T}_h f - f ]\|_{C_0}
\le \| \mathcal{T}_h f - f \|_{C_0}, \\
\| \mathcal{T}_{t-h} f  - \mathcal{T}_{t} f \|_{C_0} &\le 
\| \mathcal{T}_{t-h} [ f - \mathcal{T}_h f ]\|_{C_0}
\le \|f - \mathcal{T}_h f \|_{C_0}. 
\end{align}
As $h \rightarrow 0$, the right continuity property \eqref{rightcontinuous} implies that both of the above bounds tend to zero. 
The following lemma is taken from Yosida \cite{yosida}, pages 259 \& 264, with a slight change in notation:
\begin{lemma}
\label{balakrishnan_lemma}
(Balakrishnan \cite{balakrishnan})
Let $\mathcal{T}_t$ be an equicontinuous semigroup of class $C_0$. Let
\begin{equation}
\label{balakrishnan_subordinator}
g_{t,\beta}(\lambda) = 
\begin{cases}
\frac{1}{2 \pi i} \int_{\sigma - i\infty}^{\sigma + i\infty}
e^{t\lambda - t z^\beta} dz ; \quad \sigma > 0, t > 0, \lambda \ge 0, 0<\beta<1. 
\\
0 \quad \text{ (when $\lambda < 0$) }
\end{cases}
\end{equation}
Then the operators $\widehat{\mathcal{T}}^{\alpha}_t$ defined by
\begin{equation}
\label{balakrishnan_semigroup}
\widehat{\mathcal{T}}^{{\alpha}}_t g = 
\begin{cases}
\int_0^\infty g_{t,{\alpha/2}} (s) \mathcal{T}_s f ds \quad (t > 0) \\
f \qquad (t = 0)
\end{cases}
\end{equation}
constitute an equicontinuous semigroup of class $C_0$. Moreover, the infinitesimal generator of 
$\widehat{\mathcal{T}}^\alpha_t$, denoted $\widehat{\mathcal{A}}^\alpha$, is given for 
$f \in \text{\textup{Dom}}(\mathcal{A})$ by the two equivalent
formulas
\begin{align}
\label{balakrishnan_formulas}
\widehat{\mathcal{A}}^\alpha f &=
\frac{\sin({\alpha\pi/2} )}{\pi} \int_0^\infty \lambda^{{\alpha/2} - 1} 
(\lambda I - \mathcal{A})^{-1} (-\mathcal{A} f) d\lambda
\\
&=
\frac{1}{\Gamma(-{\alpha/2})} \int_0^\infty \lambda^{-{\alpha/2} - 1}
(\mathcal{T}_\lambda- I) f d\lambda.
\end{align}
\end{lemma}

Since Feller semigroups are equicontinuous, we may apply this result to the semigroup $\mathcal{T}^2_t$ generated by Brownian motion 
$X^2_t$, in which case we know that
\begin{equation}
\mathcal{A}^2 u = \Delta u, \quad \text{ if } u(t,\cdot) \in C^2,
\end{equation}
which is the classical result \eqref{ito_rule}. 
Next, the function $g_{t,\beta}$ in equation \eqref{balakrishnan_subordinator} with parameters $t = 1$
and $\beta = \alpha/2$ is simply the inverse Laplace transform of \eqref{subordinator_laplace}; thus, it coincides with $h_\beta$. By the derivation \eqref{subordinated_semigroups}, we see that the semigroup 
$\widehat{\mathcal{T}}^\alpha_t$ in equation \eqref{balakrishnan_semigroup} coincides with the semigroup
$\mathcal{T}^{\alpha}_t$. 

Now the lemma says that the generator $\mathcal{A}^\alpha$
of $\mathcal{T}^{\alpha}_t$ is given by the formulas \eqref{balakrishnan_formulas},
\begin{equation}
\mathcal{A}^\alpha u(t,\cdot) = \widehat{\mathcal{A}}^{\alpha} u(t,\cdot) =
\frac{\sin(\alpha \pi/2)}{\pi} \int_0^\infty \lambda^{\alpha/2 - 1} 
(\lambda I - \Delta)^{-1} (\Delta u(t,\cdot)) d\lambda
\end{equation}
whenever $u(t,\cdot) \in C^2$.
However, we recognize the expression on the right-hand side as the usual 
Balakrishnan formula, which is equivalent to the standard fractional Laplacian
on $\mathbb{R}^d$ (Kwasnicki \cite{kwasnicki}). Therefore,
\begin{equation}
\mathcal{A}^\alpha u(t,x) = -(-\Delta)^\alpha u(t,x) \text{ when } u(t,\cdot) \in C^2.
\end{equation}
This proves \eqref{equivalence_alpha}, and completes the proof of Theorem \ref{theorem_rn}.

\begin{remark}
\normalfont
The same technique may apply for the fractional power $\mathcal{L}^\alpha$ of a constant coefficient elliptic operator 
\begin{equation}
\mathcal{L} = \frac{1}{2} \sum_{i,j = 1}^{n}
a^{ij} \frac{\partial^2}{\partial x_i}{\partial x_j} + 
\sum_{i=1}^{n} b^i(x) \frac{\partial}{\partial x_i},
\end{equation}
In the book of Freidlin \cite{freidlin} a result corresponding to \eqref{ito_rule} for such operators is proven, leading to a Feynman-Kac formula
\begin{equation}
u(t,x) = \mathcal{T}^\alpha_t f (x) := \mathbb{E}_{X_0 = x} \left[ f(X_t) \right],
\end{equation}
where $X_t$ is now the Ito process
\begin{equation}
dX_t = \sigma(X_t) dB_t + b(X_t) dt, \quad X_0 = x, \quad a^{ij}(x) = \sigma(x)\sigma^*(x). 
\end{equation}
Since spectral powers can be constructed for such operators $\mathcal{L}$ and they can also be expressed by Balakrishnan formulas, we would obtain a subordinated formula
\begin{equation}
u(t,x) = \mathbb{E}_{X_0 = x} \left[ f(X_{T_{\alpha/2}(t)}) \right] = 
 \mathbb{E}_{Z_0 = x} \left[ f(X_{Z_t}) \right]
\end{equation}
for solutions to the Cauchy problem of $L^{\alpha}$.
Here $Z_t = X_{T_{\alpha/2}(t)}$ would be a L\'evy process with measure defined by the coefficients of $L$ and $\alpha$. Thus the solution would be given by nonisotropic L\'evy process. For example, if the $a^{ij}$ are constant, we expect $Z_t = X_{T_{\alpha/2}(t)}$ to be an elliptically-contoured stable L\'evy process. 
Furthermore, it may be possible to prove Feynman-Kac formulas for fractional equations with an additional zero-order term $c(x) u(x,t)$, which Freidlin \cite{freidlin} treats by comparison to an auxillary semigroup, to obtain a formula
\begin{equation}
u(x,t) = \mathbb{E}_{X^\alpha_0 = x} \left[ f(X_t) \right]
\exp{\int_0^t c(X_s) ds}.
\end{equation}
Neither of these results are not pursued here, but would be a worthwhile extensions of the present article. 

\end{remark}


\section{The Inhomogeneous Fractional Cauchy Problem on $\mathbb{R}^d$.}
\label{rn_duhamel}
In this section, we show how a stochastic solution formula for the Cauchy problem with nonzero forcing term can be obtained from the results in the previous sections using Duhamel's principle. 
We consider the following problem
\begin{align}
\label{problem_rn_inhomo}
\begin{cases}
&\partial_t u(t,x) + (-\Delta)^{\alpha/2}u (t,x) = r(x,t) \text{ for $t > 0$} \\
&u(0,x) = f(x).
\end{cases}
\end{align}
which differs from Problem \eqref{problem_rn} by the presence of the inhomogeneity $r(x,t)$. 

Let us write the solution to \eqref{problem_rn_inhomo} as $u(t,x) = u_1(t,x) + u_2(t,x)$, 
where $u_1(t,x)$ solves
\begin{align}
\label{u1_problem_rn}
\begin{cases}
&\partial_tu_1(t,x) + (-\Delta)^{\alpha/2}u_1 (t,x) = 0 \text{ for $t > 0$} \\
&u_1(0,x) = f(x)
\end{cases}
\end{align}
and $u_2(t,x)$ solves
\begin{align}
\label{u_2_rn}
\begin{cases}
&\partial_t u_2(t,x) + (-\Delta)^{\alpha/2} u_2 (t,x) = r(t,x) \text{ for $t > 0$} \\
&u_2(0,x) = 0.
\end{cases}
\end{align}
Then by Theorem \ref{theorem_rn}, provided the regularity assumptions hold,
\begin{equation}
\label{sol_u1_rn}
u_1 =  \mathbb{E}_{X^\alpha_0 = x} \left[ f(X^\alpha_t) \right].
\end{equation}

As for $u_2$, by Duhamel's principle, 
\begin{equation}
u_2(x,t)  = \int_0^t [P^s r] (x,t) ds,
\end{equation}
where $[P^s r]$ solves, for fixed $s$,
\begin{align}
\label{u2_problem_rn}
\begin{cases}
&\partial_t u (t,x) + (-\Delta)^{\alpha/2} u (t,x) = 0 \text{ for $t > s$} \\
& u (s,x) = r(s,x) \text{ for $x \in \mathbb{R}^d$.}
\end{cases}
\end{align}
Letting $\tilde{u}(t-s,x)= u(t,x)$, this problem becomes 
\begin{align}
\label{u2_problem_rn_shifted}
\begin{cases}
&\partial_t \tilde{u} (t,x) + (-\Delta)^{\alpha/2} \tilde{u} (t,x) = 0 \text{ for $t > 0$} \\
& \tilde{u} (0,x) = r(s,x) \text{ for $x \in \mathbb{R}^d$.}
\end{cases}
\end{align}
Again, by Theorem \ref{theorem_rn}, provided the regularity assumptions hold,
\begin{equation}
\tilde{u}(x,t) = \mathbb{E}_{X^\alpha_0 = x} \left[ r(s,X^\alpha_t) \right]. 
\end{equation}
so that
\begin{equation}
[P^s r](x,t) = \tilde{u}(x,t-s) = \mathbb{E}_{X^\alpha_0 = x} \left[ r(s,X^\alpha_{t-s}) \right],
\end{equation}
and therefore
\begin{equation}
\label{sol_u2_rn}
u_2(x,t) = \int_0^t \mathbb{E}_{X^\alpha_0 = x} \left[ r(s,X^\alpha_{t-s}) \right] ds
= \mathbb{E}_{X^\alpha_0 = x}  \left[  \int_0^tr(s,X^\alpha_{t-s}) ds \right].
\end{equation}
Adding \eqref{sol_u2_rn} and \eqref{sol_u1_rn} gives the following result.

\begin{theorem}
\label{theorem_rn_inhomo}
Let $u(t,x)$ solve \eqref{problem_rn_inhomo}. 
Suppose that $f$ 
and $r(s,\cdot)$, for every $s \ge 0$, are such that the problems
\eqref{u1_problem_rn} and \eqref{u2_problem_rn_shifted} each admit a solution satisfies
the regularity conditions
\eqref{rn_regularity_assumption_1} and
\eqref{rn_regularity_assumption_2}.
Then 
\begin{align}
\begin{split}
\label{sol_rn_inhomo}
u(t,x) &= \mathbb{E}_{X^\alpha_0 = x} \left[ f(X^\alpha_t) \right]
+ \mathbb{E}_{X^\alpha_0 = x}  \left[ \int_0^t  r(s, X^\alpha_{t-s}) ds \right] \\
&=
\mathbb{E}_{X^\alpha_0 = x} \left[ f(X^\alpha_t) \right]
+ \mathbb{E}_{X^\alpha_0 = x}  \left[ \int_0^t  r(t-s, X^\alpha_{s}) ds \right]
\end{split}
\end{align}
where $X^\alpha_t$ is standard isotropic $\alpha$-stable L\'evy motion in $\mathbb{R}^d$. 
\end{theorem}
We persue analogues of this result in bounded domains with the operator 
$-(\Delta_{\Omega,g})^{\alpha/2}$.


\section{The Homogeneous Spectral Fractional Cauchy Problem on a Bounded Domain, with Dirichlet Boundary Conditions.}
\label{cauchy}

Let $\Omega \subset \mathbb{R}^d$ be an open domain with Lipschitz boundary $\partial\Omega$. We consider the Cauchy problem on $\Omega$, with both initial and boundary conditions, 
\begin{align}
\label{problem_omega}
\begin{cases}
&\partial_t u(t,x) = -(-\Delta_{\Omega,g})^{\alpha/2}u (t,x) \text{ for $t > 0, x \in \Omega$} \\
&u(0,x) = f(x) \quad \text{ for $x \in \Omega$}.  \\
&u(t,x) = g(x) \quad \text{ for $x \in \partial\Omega$}. 
\end{cases}
\end{align}
We require that $f$ and $g$ are continuous, and for any $x_0 \in \partial \Omega$
\begin{equation}
\label{continuity}
\lim_{\Omega \ni x \rightarrow x_0} f(x) = g(x_0).
\end{equation}
We prove the following 
\begin{theorem}
\label{theorem_omega}
Let $u(t,x)$ be a solution to \eqref{problem_omega}, \eqref{continuity} such that
\begin{align}
\label{omega_regularity_assumption_1}
&\text{for every $t \ge 0$, $u(t,\cdot) \in C^2_x(\Omega)$, and}\\
\label{omega_regularity_assumption_2}
&\text{$\partial_t u(\cdot,\cdot)$ is uniformly continuous in spacetime $[0,\infty) \times \Omega$.}
\end{align}
Then 
\begin{equation}
\label{levy_semigroup_omega}
u(t,x) = \mathbb{E}_{X^{\Omega,\alpha}_0 = x} 
\left[ f(X^{\Omega,\alpha}_t) \chi_{\tau_\Omega > T_{\alpha/2}(t)}
+ g(X^{\Omega,\alpha}_t) \chi_{\tau_\Omega \le T_{\alpha/2}(t)}\right]
\end{equation}
where $X^{\Omega,\alpha}_t$ is standard $\alpha/2$-subordinate stopped isotropic Brownian (or rather, standard $2$-stable) motion in 
$\Omega$:
\begin{equation}
\label{stopped}
X^{\Omega,\alpha}_t = X^{\Omega, 2}_{T_{\alpha/2}(t)}, \quad
X^{\Omega, 2}_t = X^{2}_{t \wedge \tau_\Omega}, \quad \tau_\Omega = 
\inf \{ s : X^2_s \not\in \Omega\}. 
\end{equation}  
\end{theorem}

Thus, paths of $X^{\Omega,\alpha}_t$ are simply paths of Brownian motion (scaled by $\sqrt{2}$) that are stopped upon reaching the boundary $\partial\Omega$, then time-changed by the standard $\alpha/2$ stable subordinator.
To prove the theorem, the same strategy is used as in Section \ref{rn}. We start with proven results for the standard Laplacian/Brownian motion. Then, we use the subordination lemma of Balakrishnan to obtain the result for the fractional equation. 

Define the semigroup $\mathcal{T}^{\Omega, \alpha }_t$ on $C_0(\overline{\Omega})$ by
\begin{equation}
\label{F_semigroup}
\mathcal{T}^{\Omega, \alpha }_t F (x) 
=
 \mathbb{E}_{X^{\Omega,\alpha}_0 = x} 
\left[ F(X^{\Omega,\alpha}_t) \right], \quad F \in C_0(\overline{\Omega}). 
\end{equation}
This is well defined, since the process $X^{\Omega,\alpha}_t$ remains in $\overline{\Omega}$ for all time. Moreover, $\mathcal{T}^{\Omega, \alpha }_t$ is Feller semigroup of the same process.  The reason for introducing this semigroup is that 
the condition \eqref{continuity} imples that the initial condition $f$ can be extended to the $C_0(\overline{\Omega})$ function $\overline{f}$ on $\overline{\Omega}$:
\begin{align}
\label{F_def}
\overline{f}(x) = 
\begin{cases}
&f(x), x \in \Omega \\
&g(x) = \lim_{\Omega \ni x \rightarrow x_0} f(x), x \in \partial\Omega
\end{cases}
\implies 
\overline{f} \in C_0(\overline{\Omega}).
\end{align}
Then the proposed solution \eqref{levy_semigroup_omega} can be written as the Feller semigroup 
\begin{equation}
\label{levy_semigroup_omega_F}
u(t,x) = \mathcal{T}^{\Omega, \alpha }_t \overline{f} (x).
\end{equation}
We shall now prove that $u(t,x)$ in this form satisfies the problem. First, we prove the following
\begin{lemma}
\label{lemma_generator_omega}
The generator $\mathcal{A}^{\alpha}$ of the semigroup $\mathcal{T}^{\Omega, \alpha }_t$ on $C_0(\overline{\Omega})$ defined by \eqref{F_semigroup} agrees with the operator 
$-(-\Delta_{\Omega,g})^{\alpha/2}$ on 
\begin{equation}
C_0(\overline{\Omega}) \cap C^2(\Omega) 
\equiv 
\left\{
u \in C_0(\overline{\Omega}) 
\text{ such that } 
u|_{\Omega} \in C^2(\Omega)
\right\}. 
\end{equation}
\end{lemma}

To prove this lemma, we begin by referring to Chung and Zhao \cite{chung_zhao}, Theorem 2.13, p. 56, where it is shown under more general conditions for the case $\alpha = 2$. In that case, the operator $-(-\Delta_{\Omega,g})^{\alpha/2}$ reduces to the standard Laplacian $\Delta$ and the process 
$X^{2,\Omega}_t$ to stopped standard isotropic $2$-stable motion. 
Thus, 
\begin{equation}
C_0(\overline{\Omega}) \cap C^2(\Omega) \subset \text{Dom}({\mathcal{A}^2}), \text{ and }
\mathcal{A}^2 = \Delta \text{ on }C_0(\overline{\Omega}) \cap C^2(\Omega). 
\end{equation}
Then, the result of Balakrishnan (Lemma \ref{balakrishnan_lemma}) tells us that the semigroup
\begin{equation}
\widehat{\mathcal{T}}^{\Omega,\alpha}_t F = 
\begin{cases}
\int_0^\infty h_{\alpha/2} (s) \mathcal{T}^{\Omega,2}_s F ds \quad (t > 0) \\
F \qquad (t = 0)
\end{cases}
\end{equation}
has infinitesimal generator, for $F \in \text{Dom}(\mathcal{A})$, and in particular,  
${F \in C_0(\overline{\Omega}) \cap C^2(\Omega)}$,
\begin{align}
\mathcal{A}^{\alpha} F (x) = 
\frac{1}{\Gamma(-\alpha/2)} \int_0^\infty \lambda^{-\alpha/2 - 1}
(\mathcal{T}^{\Omega,2}_\lambda- I) F (x) d\lambda.
\end{align}
We begin with the right-hand side of this equation. 
In Freidlin \cite{freidlin}, formula \eqref{levy_semigroup_omega_F} is shown to satisfy Problem \eqref{problem_omega} when $\alpha = 2$. In other words, 
\begin{equation}
\label{heat_semigroup}
\mathcal{T}^{\Omega, 2 }_t \overline{f} (x) = e^{t\Delta_{\Omega,g}} f(x), x \in \Omega,
\end{equation}
Therefore,  for $F \in C_0(\overline{\Omega}) \cap C^2(\Omega)$  and $x \in \Omega$, by \eqref{heat_semigroup}, we have 
\begin{align}
\mathcal{A}^{\alpha} F (x) &= 
\frac{1}{\Gamma(-\alpha/2)} \int_0^\infty \lambda^{-\alpha/2 - 1}
(e^{t\Delta_{\Omega,g}}- I) F (x) d\lambda \\
&= 
-(-\Delta_{\Omega,g})^{\alpha/2} F (x).
\end{align}
It remains to show that $
\widehat{\mathcal{T}}^{\Omega,\alpha}_t = \mathcal{T}^{\Omega,\alpha}_t.
$
In order to do this, we will use the definition \eqref{stopped} of the stopped process 
$X^{\Omega, \alpha}_t$ as subordinate stopped Brownian motion, and the tower property:
\begin{align}
\begin{split}
\mathcal{T}^{\Omega,\alpha}_t F(x) &= \mathbb{E}_{X^{2}_0 = x} \left[ F(X^{\Omega,\alpha}_{t}) \right] \\
&= \mathbb{E}_{X^{2}_0 = x} \left[ F\left(X^{\Omega,2}_{T_{\alpha/2}(t)}\right) \right] \\
 &= \mathbb{E}_{X^{2}_0 = x} \left[  \int \mathbb{P} \{T_{\alpha/2} = s \} F(X^{\Omega,2}_{s})  ds \right] \\
 &= \int \mathbb{P}\{T_{\alpha/2} = s \} \mathbb{E}_{X^{2}_0 = x} \left[   F(X^{\Omega,2}_{s})  \right]  ds\\
 &= \int \mathbb{P}\{T_{\alpha/2}  = s \} \mathcal{T}^2_s F(x)  ds \\
 & = \int h_{\alpha/2}(s) \mathcal{T}^2_s F(x)  ds \\
& = \widehat{\mathcal{T}}^{\Omega,\alpha}_t F(x) .
\end{split}
\end{align}
The proof of the lemma is complete. 

Let us now prove Theorem \ref{theorem_omega}. 
We know by Lemma \ref{lemma_acp} that $\mathcal{T}^{\Omega,\alpha}_t \overline{f}(x)$ is the unique solution to the abstract Cauchy problem
\begin{align}
\begin{cases}
&\frac{d}{dt} u(t,x) = \mathcal{A}^{\alpha} u(t,x) \text{ for $x \in \overline{\Omega}, t > 0$} \\
&\lim_{t \rightarrow 0} u(t,x) = \overline{f}(x).
\end{cases}
\end{align}
For any path starting at $x \in \partial\Omega$, $\tau_\Omega = 0$, so that
$X^{\Omega, \alpha}_t = x$,  and therefore
\begin{equation}
\mathcal{T}^{\Omega, \alpha }_t \overline{f} (x)  =
\mathbb{E}_{X^{\Omega,\alpha}_0 = x} 
\left[ \overline{f}(x) \right]
=
\mathbb{E}_{X^{\Omega,\alpha}_0 = x} 
\left[ g(x) \right]
= g(x). 
\end{equation}
This shows that
\begin{equation}
\frac{d}{dt} \mathcal{T}^{\Omega, \alpha }_t \overline{f} (x) = 0 \text{ for } x \in \partial\Omega. 
\end{equation}
Therefore, an equivalent statement is that $\mathcal{T}^{\Omega, \alpha }_t \overline{f}$ is the unique solution to the problem  
\begin{align}
\label{modified_acp}
\begin{cases}
&\frac{d}{dt} u(t,x) = \mathcal{A}^{\alpha} u(t,x) \text{ for $x \in \Omega, t > 0$} \\
&\frac{d}{dt} u(t,x) = 0 \text{ for $x \in \partial \Omega, t > 0$} \\
&\lim_{t \rightarrow 0} u(t,x) = \overline{f}(x) 
\end{cases}
\end{align}
On the other hand, we note that, by Lemma \ref{lemma_generator_omega}, 
\begin{equation}
u(t,\cdot) \in  C_0(\overline{\Omega}) \cap C^2(\Omega) \subset \text{Dom}(\mathcal{A}^\alpha).
\end{equation}
and the differential equation in the problem \eqref{problem_omega} together with the regularity condition \eqref{omega_regularity_assumption_2} is equivalent to the first differential equation in the problem \eqref{modified_acp}. Also, $u(t,x)$ clearly satisfies the remaining two equations in \eqref{modified_acp}, so by uniqueness, 
\begin{equation}
u(t,x) = \mathcal{T}^{\Omega, \alpha }_t \overline{f} (x). 
\end{equation}
This completes the proof. 

\section{The Inhomogeneous Spectral Fractional Cauchy Problem on a Bounded Domain, with Dirichlet Boundary Conditions.}
\label{cauchy_duhamel}

We consider, on bounded domain $\Omega \subset \mathbb{R}^d$, the mixed problem with both initial and boundary conditions: 
\begin{align}
\label{problem_omega_inhomo}
\begin{cases}
&\partial_t u(t,x) + (-\Delta_{\Omega,g})^{\alpha/2}u (t,x) = r(x,t) \text{ for $t > 0, x \in \Omega$} \\
&u(t,x) = g(x) \quad \text{ for $t > 0, x \in \partial\Omega$} \\
&u(0,x) = f(x) \quad \text{ for $x \in \Omega$}. 
\end{cases}
\end{align}
Again, we require that $f$ and $g$ are continuous, and for any $x_0 \in \partial \Omega$
\begin{equation}
\lim_{\Omega \ni x \rightarrow x_0} f(x) = g(x_0).
\end{equation}

Let us write the solution to \eqref{problem_omega_inhomo} as $u(t,x) = u_1(t,x) + u_2(t,x)$, 
where $u_1(t,x)$ solves
\begin{align}
\label{u1_problem_omega}
\begin{cases}
&\partial_t u_1(t,x) + (-\Delta_{\Omega,g})^{\alpha/2}u_1 (t,x) = 0 \text{ for $t > 0, x \in \Omega$} \\
&u_1(t,x) = g(x) \quad \text{for $t > 0, x \in \partial\Omega$} \\
&u_1(0,x) = f(x) \quad \text{ for $x \in \Omega$}
\end{cases}
\end{align}
and $u_2(t,x)$ solves
\begin{align}
\label{u_2_omega}
\begin{cases}
&\partial_t u_2(t,x) + (-\Delta_{\Omega,0})^{\alpha/2} u_2 (t,x) = r(t,x)\text{ for $t > 0, x \in \Omega$} \\
&u_2(t,x) = 0 \quad \text{ for $t > 0, x \in \partial\Omega$} \\
&u_2(0,x) = 0  \quad \text{ for $x \in \Omega$}.
\end{cases}
\end{align}
Then by Theorem \ref{theorem_omega}, provided the regularity assumptions hold,
\begin{equation}
\label{sol_u1_omega}
u_1 =   \mathbb{E}_{X^{\Omega,\alpha}_0 = x} 
\left[ f(X^{\Omega,\alpha}_t) \chi_{\tau_\Omega > T_{\alpha/2}(t)}
+ g(X^{\Omega,\alpha}_t) \chi_{\tau_\Omega \le T_{\alpha/2}(t)}\right].
\end{equation}
By Duhamel's principle, 
\begin{equation}
u_2(x,t)  = \int_0^t [P^s r] (x,t) ds,
\end{equation}
where $[P^s r]$ solves, for fixed $s$,
\begin{align}
\begin{cases}
&\partial_t u (t,x) + (-\Delta_{\Omega,0})^{\alpha/2} u (t,x) = 0 \text{ for $t > 0, x \in \Omega$} \\
&u(t,x) = 0 \quad \text{ for $t > 0, x \in \partial\Omega$} \\
&u (s,x) = r(s,x) \text{ for $x \in \Omega$.}
\end{cases}
\end{align}
Letting $\tilde{u}(t-s,x)= u(t,x)$, this problem becomes 
\begin{align}
\label{u2_problem_omega_shifted}
\begin{cases}
&\partial_t \tilde{u} (t,x) + (-\Delta_{\Omega,0})^{\alpha/2} \tilde{u} (t,x) = 0 \text{ for $t > 0, x \in \Omega$} \\
& u(t,x) = 0 \quad \text{ for $t > 0, x \in \partial\Omega$}. \\
& \tilde{u} (0,x) = r(s,x) \text{ for $x \in \Omega$.}
\end{cases}
\end{align}
Again, by Theorem \ref{theorem_omega},  provided the regularity assumptions hold,
\begin{equation}
\tilde{u}(x,t) =  \mathbb{E}_{X^{\Omega,\alpha}_0 = x} 
\left[ r(s, X^{\Omega,\alpha}_t) \chi_{\tau_\Omega > T_{\alpha/2}(t)}\right]
\end{equation}
so that
\begin{equation}
[P^s r](x,t) = \tilde{u}(x,t-s) = \mathbb{E}_{X^\alpha_0 = x} \left[ r(s,X^{\Omega,\alpha}_{t-s}) \chi_{\tau_\Omega > T_{\alpha/2}(t-s)} \right],
\end{equation}
and therefore
\begin{align}
\label{sol_u2_omega}
u_2(x,t) &= \int_0^t \mathbb{E}_{X^\alpha_0 = x} \left[ r(s,X^{\Omega,\alpha}_{t-s})  \chi_{\tau_\Omega > T_{\alpha/2}(t-s)}  \right] ds \\
&= \mathbb{E}_{X^{\Omega,\alpha}_0 = x}  \left[  \int_0^tr(s,X^{\Omega,\alpha}_{t-s}) \chi_{\tau_\Omega > T_{\alpha/2}(t-s)} ds \right].
\end{align}
Adding \eqref{sol_u2_omega} and \eqref{sol_u1_omega} gives the following result. 

\begin{theorem}
\label{theorem_omega_inhomo}
Let $u(t,x)$ solve \eqref{problem_omega_inhomo}. 
Suppose that $f$ 
and $r(s,\cdot)$ for every $s \ge 0$, are such that the problems
\eqref{u1_problem_omega} and \eqref{u2_problem_omega_shifted} each admit a solution satisfies
the regularity conditions
\eqref{omega_regularity_assumption_1} and
\eqref{omega_regularity_assumption_2}.
Then 
\begin{align}
\begin{split}
\label{sol_omega_inhomo}
u(t,x) &= \mathbb{E}_{X^{\Omega,\alpha}_0 = x} 
\left[ f(X^{\Omega,\alpha}_t) \chi_{\tau_\Omega > T_{\alpha/2}(t)}
+ g(X^{\Omega,\alpha}_t) \chi_{\tau_\Omega \le T_{\alpha/2}(t)}\right] \\
&\qquad \qquad \qquad \qquad + \mathbb{E}_{X^\alpha_0 = x}  \left[ \int_0^t  r(s, X^{\Omega,\alpha}_{t-s})  \chi_{\tau_\Omega > T_{\alpha/2}(t-s)}  ds \right] \\
&=
\mathbb{E}_{X^{\Omega,\alpha}_0 = x} 
\left[ f(X^{\Omega,\alpha}_t) \chi_{\tau_\Omega > T_{\alpha/2}(t)} 
+ g(X^{\Omega,\alpha}_t) \chi_{\tau_\Omega \le T_{\alpha/2}(t)}\right] \\
&\qquad \qquad \qquad \qquad + \mathbb{E}_{X^\alpha_0 = x}  \left[ \int_0^t  r(t-s, X^{\Omega,\alpha}_{s})  \chi_{\tau_\Omega > T_{\alpha/2}(s)}  ds \right]
\end{split}
\end{align}
where $X^{\Omega,\alpha}_t$ is subordinate stopped standard isotropic $2$-stable L\'evy motion in $\Omega$. 
\end{theorem}

\section{Regularity and Steady-State for the Spectral Fractional Heat Equation in a Bounded Domain and Survival Probability of Subordinate Stopped Brownian Motion.}
\label{regularity_section}

The previous section contains stochastic solution formulas for the fractional Cauchy problem that are contingent upon the regularity of the solutions to problems \eqref{u1_problem_omega} and \eqref{u2_problem_omega_shifted} . The latter can be reduced to the case of zero boundary value by harmonic lifting \cite{Cusimano2017, AntilPfeffererRogovs}. We now investigate regularity for such problems, together with convergence to a steady state. Our results allow us to obtain exponential decay of the survival probability of subordinate stopped Brownian motion, and obtain a  stochastic solution formulas for the Dirichlet problem by applying the formula in the previous section to a related Cauchy problem and taking the limit as $t \rightarrow \infty$. We proceed by working with a classical solution formula (eigenfunction expansion) directly and adapting the proof in \cite{larsson2008partial}. The classical solution formula \eqref{w_series} appears, e.g., in \cite{bonito2017approximation}, but we did not find a precise statement of the results below in the literature. Appendix \ref{appendix} contains background material used in this section.


Consider the fractional heat equation, for bounded and smooth domain $\Omega$,
\begin{align}
\label{fractional_parabolic_ally_transient}
\begin{cases}
&\partial_t w+ (-\Delta_{\Omega,0})^{\alpha/2} w
 = 0 \text{ for $ x \in \Omega$} \\
&w (x,t) = 0 \quad \text{ for $x \in \partial\Omega$} \\
&w (x,0) = w_0 \in L^2 \quad \text{ for $x \in \Omega$}.
\end{cases}
\end{align}

\begin{theorem}
\label{steady_state}
Define 
\begin{equation}
\label{w_series}
w(t,x) 
=
\sum_{j = 1}^\infty e^{-t\lambda_j^{\alpha/2} }
(w_0, e_j)_{L^2(\Omega)} e_j(x).
\end{equation}
Then, 
\begin{enumerate}
\item
For any $t>0$, integers $K, \ell \ge 0$, and $\mu < 1$, there exists $C > 0$ such that
\begin{equation}
\|  \partial_t^\ell w(t,\cdot) \|_{H^{K}(\Omega)} \le
C
t^{-\ell-K/\alpha}
e^{-\mu \lambda_1^{\alpha/2} t} \| w_0 \|_{L^2(\Omega)}.
\end{equation}
The convention is that $H^0(\Omega) = L^2(\Omega)$. 
\item
For any $t > 0$, integers $r, \ell > 0$, and $\mu < 1$, there exists $C > 0$ such that
\begin{equation}
\| \partial_t^\ell w(t, \cdot) \|_{C^r(\Omega)} \le
t^{-\ell-{\left\lceil \frac{d}{2} + r\right\rceil}/\alpha}
e^{-\mu \lambda_1^{\alpha/2} t} \| w_0 \|_{L^2(\Omega)}.
\end{equation}
where ${\left\lceil x \right\rceil}$ is the smallest integer $\ge x$.
Therefore, $w$ is smooth in $(0,\infty) \times \Omega$.
\item
$\lim\limits_{t \rightarrow 0} w(t,x) = w_0$ in $L^2(\Omega)$. Moreover, if $w_0 \in H^k(\Omega)$ for $k > 0$ and satisfies $w_0|_{\partial \Omega} = 0$, then $\lim\limits_{t \rightarrow 0} w(t,x) = w_0$ in $H^k(\Omega)$.
\item
The series $w$ given by \eqref{w_series} is the unique classical solution to the problem \eqref{fractional_parabolic_ally_transient}.
\end{enumerate}
\end{theorem}

To prove this, we define
\begin{equation}
S^{\ell,J}
=
(-1)^\ell \sum_{j = J}^\infty \lambda_j^{\ell \alpha / 2} e^{-\lambda_j^{\alpha / 2}t}
(w, e_j) e_j(x). 
\end{equation}
Note that $S^{\ell,J=1}$ corresponds to term-by-term differentiation of
the function $w(t,x)$. However, since the validity of term-by-term differentiation is
not known \emph{a priori}, we must establish it by studying $S^{\ell, J}$ for $J \rightarrow 
\infty$. 
We begin by proving the following lemma:
\begin{lemma}
\label{S_lemma}
For $t > 0$ and any integers $k, \ell \ge 0$, 
\begin{equation}
\label{lemma_series}
(-\Delta)^k S^{\ell, J} =
(-1)^\ell \sum_{j = J}^\infty \lambda_j^{\ell \alpha / 2 + k} e^{-\lambda_j^{\alpha / 2}t}
(w, e_j) e_j(x). 
\end{equation}
For any $\mu < 1$, there exists $C$ such that 
\begin{equation}
\label{lemma_l2}
\|(-\Delta)^k S^{\ell, J} 
\|_{L^2(\Omega)}
\le
C
t^{-\ell-k(2/\alpha)} e^{-\mu \lambda_J^{\alpha/2} t} \| w_0 \|_{L^2(\Omega)}. 
\end{equation}
Moreover, $(-\Delta)^k S^{\ell, J} \in H^1_0(\Omega)$, and for any $\mu < 1$, 
there exists $C$ such that
\begin{equation}
\label{lemma_h1}
\left[
(-\Delta)^k S^{\ell, J} 
\right]_{H^1_0(\Omega)}
=
\|\nabla
(-\Delta)^k S^{\ell, J} 
\|_{L^2(\Omega)}
\le
C
t^{-\ell-(k+1/2)(2/\alpha)} e^{-\mu \lambda_J^{\alpha/2} t} \| w_0 \|_{L^2(\Omega)}. 
\end{equation}
\end{lemma}

We prove the lemma by induction. For $k = 0$, the representation \eqref{lemma_series} is just the definition of $S^{\ell, J}$. To prove the $L^2$ bound \eqref{lemma_l2}, we note
\begin{align}
\label{orthogonality}
\begin{split}
(S^{\ell, J}, e_j)^2_{L^2} 
&=
\left(
(-1)^\ell \sum_{\tilde{j} = J}^\infty \lambda_{\tilde{j}}^{\ell \alpha / 2 } e^{-\lambda_{\tilde{j}}^{\alpha / 2}t}
(w_0, e_{\tilde{j}}) e_{\tilde{j}}(x),e_j
\right)^2_{L^2} \\
&=
\lambda_j^{\ell \alpha }
e^{-2 \lambda_j^{\alpha/2} t} (w_0, e_j)^2,
\end{split}
\end{align}
and write
\begin{align}
\| S^{\ell, J }\|^2_{L^2(\Omega)} &=
\sum_{j = 1}^\infty (S^{\ell, J}, e_j)_{L^2}^2 \\
&=
\sum_{j = J}^\infty \lambda_j^{\ell \alpha }
e^{-2 \lambda_j^{\alpha/2} t} (w_0, e_j)^2.
\end{align}
Next, we need the following fact:
\begin{align}
\label{the_fact}
\begin{split}
&\text{For any $\theta \ge 0$ and $\mu < 1$, there exists a constant $C$}\\
&\text{such that } s^\theta e^{-2s} \le C e^{-2 \mu s} \text{ for all } s \in [0, \infty).
\end{split}
\end{align}
Using the fact \eqref{the_fact}, we can write, for any $\mu < 1$, 
\begin{equation}
\lambda_j^{\ell \alpha }
e^{-2 \lambda_j^{\alpha/2} t} 
=
t^{-2\ell}
\left( \lambda_j^{\alpha/2} t \right)^{2\ell}
e^{-2 \lambda_j^{\alpha/2} t} 
 \le
C
t^{-2\ell}
e^{-2 \mu \lambda_j^{\alpha/2} t}.
\end{equation}
Therefore, 
\begin{align}
\| S^{\ell, J }\|^2_{L^2(\Omega)}
& \le 
C t^{-2\ell}
\sum_{j = J}^\infty 
e^{-2 \mu \lambda_j^{\alpha/2} t}
(w_0, e_j)^2 \\
& \le 
C t^{-2\ell}
e^{-2 \mu \lambda_J^{\alpha/2} t}
\sum_{j = J}^\infty 
(w_0, e_j)^2 \\
&\le
C t^{-2\ell}
e^{-2 \mu \lambda_J^{\alpha/2} t}
\|w_0\|^2.
\end{align}
Taking the square-root yields \eqref{lemma_l2} with $k = 0$:
\begin{equation}
\| S^{\ell, J }\|_{L^2(\Omega)}
\le
C t^{-\ell}
e^{- \mu \lambda_J^{\alpha/2} t}
\|w_0\|_{L^2(\Omega)}.
\end{equation}
Next, we have $S^{\ell, J} \in H^1_0(\Omega)$ if and only if 
$\sum_{j = 1}^\infty \lambda_j (S^{\ell, J} , e_j)^2_{L^2(\Omega)}$
converges, and 
\begin{equation}
\left[
S^{\ell, J} 
\right]^2_{H^1_0(\Omega)}
=
\|\nabla S^{\ell, J} 
\|^2_{L^2(\Omega)}
=
\sum_{j = 1}^\infty \lambda_j (S^{\ell, J} , e_j)^2_{L^2(\Omega)}.
\end{equation}
See Theorem \ref{H10_theorem} in Appendix \ref{appendix}.
From \eqref{orthogonality}, have
\begin{align}
\left[
S^{\ell, J} 
\right]^2_{H^1_0(\Omega)}
&=
\sum_{J = 1}^\infty 
\lambda_j^{\ell \alpha + 1}
e^{-2 \lambda_j^{\alpha/2} t} (w_0, e_j)^2.
\end{align}
Using the fact \eqref{the_fact} again, we can write, for any $\mu < 1$, 
\begin{equation}
\lambda_j^{\ell \alpha + 1 }
e^{-2 \lambda_j^{\alpha/2} t} 
=
t^{-2 ( \ell + 1/\alpha ) }
\left( \lambda_j^{\alpha/2} t \right)^{2 ( \ell + 1/\alpha )}
e^{-2 \lambda_j^{\alpha/2} t} 
 \le
C
t^{-2(\ell+ 1/\alpha)}
e^{-2 \mu \lambda_j^{\alpha/2} t}.
\end{equation}
Therefore, 
\begin{align}
\left[
S^{\ell, J} 
\right]^2_{H^1_0(\Omega)}
& \le 
C t^{-2( \ell + 1/\alpha ) }
\sum_{j = J}^\infty 
e^{-2 \mu \lambda_j^{\alpha/2} t}
(w_0, e_j)^2 \\
& \le 
C t^{-2 ( \ell + 1/\alpha ) }
e^{-2 \mu \lambda_J^{\alpha/2} t}
\sum_{j = J}^\infty 
(w_0, e_j)^2 \\
&\le
C t^{-2 ( \ell + 1/\alpha ) }
e^{-2 \mu \lambda_J^{\alpha/2} t}
\|w_0\|_{L^2(\Omega)}^2.
\end{align}
Taking the square-root yields \eqref{lemma_h1} with $k = 0$:
\begin{equation}
\left[
S^{\ell, J} 
\right]_{H^1_0(\Omega)}
\le
C t^{-( \ell + 1/\alpha ) }
e^{- \mu \lambda_J^{\alpha/2} t}
\|w_0\|_{L^2(\Omega)}.
\end{equation}
This completes the proof of the lemma for $k = 0$.

Now we perform the induction step. Suppose the lemma is true for a certain 
$k$; we show that this implies it is true for $k + 1$. We know that 
$(-\Delta)^k S^{\ell, J} \in H^1_0(\Omega)$, i.e., $(-\Delta)^k S^{\ell, J}$
is zero on $\partial \Omega$. Therefore, the representation
\begin{equation}
(-\Delta) = \sum^\infty_{j = 1} \lambda_j (\cdot, e_j)_{L^2} e_j
\end{equation}
is valid on $(-\Delta)^k S^{\ell, J}$. Since we also have \eqref{lemma_series}, we 
see that 
\begin{align}
\label{orthogonality_not_squared}
\begin{split}
\left((-\Delta)^k S^{\ell, J}, e_j\right)_{L^2} 
&=
\left(
(-1)^\ell \sum_{\tilde{j} = J}^\infty \lambda_{\tilde{j}}^{\ell \alpha / 2 + k } e^{-\lambda_{\tilde{j}}^{\alpha / 2}t}
(w_0, e_{\tilde{j}}) , e_{\tilde{j}}
\right)_{L^2} \\
&=
(-1)^\ell
\lambda_j^{\ell \alpha / 2 + k}
e^{- \lambda_j^{\alpha/2} t} (w_0, e_j).
\end{split}
\end{align}
Therefore, 
\begin{equation}
(-\Delta)^{k+1} S^{\ell, J} = 
(-\Delta) (-\Delta)^k S^{\ell, J} = (-1)^\ell
 \sum_{j = J}^{\infty}
\lambda_j^{\ell \alpha / 2 + k + 1}
e^{- \lambda_j^{\alpha/2} t} (w_0, e_j)
e_j.
\end{equation}
This is \eqref{lemma_series} for the case $k+1$. We can use it to obtain the 
$L^2$ bound as before; first we note that the representation implies
\begin{align}
\label{orthogonality_general}
\begin{split}
\left((-\Delta)^{k+1} S^{\ell, J}, e_j\right)^2_{L^2} 
&=
\left(
(-1)^\ell \sum_{\tilde{j} = J}^\infty \lambda_{\tilde{j}}^{\ell \alpha / 2 + k + 1 } e^{-\lambda_{\tilde{j}}^{\alpha / 2} t}
(w_0, e_{\tilde{j}}) , e_{\tilde{j}}
\right)^2_{L^2} \\
&=
\lambda_j^{\ell \alpha  + 2 (k+1)}
e^{- 2 \lambda_j^{\alpha/2} t}
(w_0, e_j)_{L^2(\Omega)}^2.
\end{split}
\end{align}
Therefore, using the fact \eqref{the_fact}, 
\begin{align}
\| (-\Delta)^{k+1} S^{\ell, J }\|^2_{L^2(\Omega)} &=
\sum_{j = 1}^\infty \left((-\Delta)^{k+1} S^{\ell, J}, e_j\right)_{L^2}^2 \\
&=
\sum_{j = J}^\infty \lambda_j^{\ell \alpha + 2(k+1)}
e^{-2 \lambda_j^{\alpha/2} t} (w_0, e_j)_{L^2}^2 \\
&=
\sum_{j = J}^\infty
t^{-2 \ell - 2(k+1)(2/\alpha)}
\left( \lambda_j^{\alpha/2} t
\right)^{2 \ell + 2(k+1)(2/\alpha)}
e^{-2 \lambda_j^{\alpha/2} t} (w_0, e_j)_{L^2}^2 \\
& \le 
C t^{-2 \ell - 2(k+1)(2/\alpha)}
\sum_{j = J}^\infty 
e^{-2 \mu \lambda_j^{\alpha/2} t}
(w_0, e_j)_{L^2}^2 \\
& \le 
C t^{-2 \ell - 2(k+1)(2/\alpha)}
e^{-2 \mu \lambda_J^{\alpha/2} t}
\sum_{j = J}^\infty 
(w_0, e_j)_{L^2}^2 \\
&\le
C t^{-2 \ell - 2(k+1)(2/\alpha)}
e^{-2 \mu \lambda_J^{\alpha/2} t}
\|w_0\|_{L^2}^2.
\end{align}
Taking the square-root yields \eqref{lemma_l2} for $k+1$:
\begin{equation}
\| (-\Delta)^{k+1} S^{\ell, J }\|^2_{L^2(\Omega)}
\le
C t^{- \ell - (k+1)(2/\alpha)}
e^{- \mu \lambda_J^{\alpha/2} t}
\|w_0\|_{L^2}.
\end{equation}
Next, we have $ (-\Delta)^{k+1} S^{\ell, J} \in H^1_0(\Omega)$ if and only if 
$\sum_{j = 1}^\infty \lambda_j ((-\Delta)^{k+1} S^{\ell, J} , e_j)^2_{L^2(\Omega)}$
converges, and 
\begin{equation}
\left[
(-\Delta)^{k+1} S^{\ell, J} 
\right]^2_{H^1_0(\Omega)}
=
\|\nabla (-\Delta)^{k+1} S^{\ell, J} 
\|^2_{L^2(\Omega)}
=
\sum_{j = 1}^\infty \lambda_j \left((-\Delta)^{k+1} S^{\ell, J} , e_j\right)^2_{L^2(\Omega)}.
\end{equation}
From \eqref{orthogonality_general}, have
\begin{align}
\left[
(-\Delta)^{k+1} S^{\ell, J} 
\right]^2_{H^1_0(\Omega)}
&=
\sum_{j = J}^\infty 
\lambda_j^{\ell \alpha  + 2 (k+1) + 1}
e^{- 2 \lambda_j^{\alpha/2} t}
(w_0, e_j)_{L^2(\Omega)}^2
\end{align}
Therefore, using the fact \eqref{the_fact},
\begin{align}
\left[
(-\Delta)^{k+1} S^{\ell, J} 
\right]^2_{H^1_0(\Omega)}
&=
\sum_{j = J}^\infty \lambda_j^{\ell \alpha + 2(k+1) + 1}
e^{-2 \lambda_j^{\alpha/2} t} (w_0, e_j)^2 \\
&=
\sum_{j = J}^\infty
t^{-2 \ell - \left(2(k+1)+1\right)(2/\alpha)}
\left( \lambda_j^{\alpha/2} t
\right)^{2 \ell + \left(2(k+1)+1\right)(2/\alpha)}
e^{-2 \lambda_j^{\alpha/2} t} (w_0, e_j)^2 \\
& \le 
C t^{-2 \ell - \left(2(k+1)+1\right)(2/\alpha)}
\sum_{j = J}^\infty 
e^{-2 \mu \lambda_j^{\alpha/2} t}
(w_0, e_j)^2 \\
& \le 
C t^{-2 \ell - \left(2(k+1)+1\right)(2/\alpha)}
e^{-2 \mu \lambda_J^{\alpha/2} t}
\sum_{j = J}^\infty 
(w_0, e_j)^2 \\
&\le
C t^{-2 \ell - \left(2(k+1)+1\right)(2/\alpha)}
e^{-2 \mu \lambda_J^{\alpha/2} t}
\|w_0\|_{L^2}^2.
\end{align}
Taking the square-root yields \eqref{lemma_h1} for $k+1$:
\begin{equation}
\left[
(-\Delta)^{k+1} S^{\ell, J} 
\right]_{H^1_0(\Omega)}
\le
C t^{- \ell - \left((k+1)+1/2\right)(2/\alpha)}
e^{- \mu \lambda_J^{\alpha/2} t}
\|w_0\|_{L^2(\Omega)}.
\end{equation}
This completes the proof of the lemma. 

Now that Lemma \ref{S_lemma} has been proven, to prove the theorem 
we invoke elliptic regularity (Theorem \ref{elliptic_regularity} in Appendix \ref{appendix}) to obtain, for $\delta = 0 \text{ or }1$,
\begin{equation}
\label{repeated_regularity}
\| S^{\ell, J}  \|_{H^{2k+\delta}(\Omega)} \le C \| (-\Delta)^k S^{\ell, J}  \|_{H^\delta (\Omega)}.
\end{equation}
The two estimates in Lemma \ref{S_lemma} can be wrapped into one estimate, again for $\delta = 0$ or $1$, 
\begin{align}
\label{combined_estimate}
\begin{split}
\|(-\Delta)^k S^{\ell, J} 
\|_{H^\delta(\Omega)}
&\le
\|
(-\Delta)^k S^{\ell, J} 
\|
_{L^2(\Omega)}
+
\delta
\left[ 
(-\Delta)^k S^{\ell, J} 
\right]_{H^1_0(\Omega)}
\\
&\le 
C \left[
t^{-\ell-k(2/\alpha)}
+ \delta
t^{-\ell-(k+1/2)(2/\alpha)}
\right]
e^{-\mu \lambda_J^{\alpha/2} t} \| w_0 \|_{L^2(\Omega)}. 
\end{split}
\end{align}
Together, \eqref{repeated_regularity} and \eqref{combined_estimate} yield
\begin{equation}
\label{final_S_sobolev_estimate_delta}
\| S^{\ell, J}  \|_{H^{2k+\delta}(\Omega)} \le 
C \left[
t^{-\ell-k(2/\alpha)}
+ \delta
t^{-\ell-(k+1/2)(2/\alpha)}
\right]
e^{-\mu \lambda_J^{\alpha/2} t} \| w_0 \|_{L^2(\Omega)}.
\end{equation}
If an integer $K$ is even, we can take $K = 2k$ and $\delta = 0$ in the above estimate; the second term in the square brackets vanishes, and the first term involves a factor $t^{-\ell-K/\alpha}$. If $K$ is odd, we can take $K = 2k+1$, $\delta = 1$ in the above estimate. In this case, the second term in the square brackets can be written $t^{-\ell-K/\alpha}$, and dominates the first term for small $t$; for large $t$, by adjusting $C$ and $\mu$, the first term can be dropped. Thus, for any integer $K$, we obtain
\begin{equation}
\label{final_S_sobolev_estimate}
\| S^{\ell, J}  \|_{H^{K}(\Omega)} \le 
C 
t^{-\ell-K/\alpha}
e^{-\mu \lambda_J^{\alpha/2} t} \| w_0 \|_{L^2(\Omega)}
\end{equation}
This is a direct analogue of Eq. 8.17 in \cite{larsson2008partial}.
From Theorem \ref{standard_sobolev_inequality}, we obtain that for all integers $r, l \ge 0$, there exists $C$ such that
\begin{equation}\label{S_smooth_estimate}
\| S^{\ell, J}  \|_{C^r(\Omega)} \le
C \| S^{\ell, J}  \|_{ H^{\left\lceil \frac{d}{2} + r\right\rceil}(\Omega) }
\le 
C 
t^{-\ell-{\left\lceil \frac{d}{2} + r\right\rceil}/\alpha}
e^{-\mu \lambda_J^{\alpha/2} t} \| w_0 \|_{L^2(\Omega)}.
\end{equation}
From \eqref{final_S_sobolev_estimate} and \eqref{S_smooth_estimate},
considering the particular case $r = 0$, and taking the supremum over 
$[t_1,t_2]$ with $0 < t_1 < t_2 < \infty$, we have
\begin{equation}
\| S^{\ell, J} (\cdot, x) \|_{C^0\left([t_1,t_2] \right)} \le
\| S^{\ell, J} \|_{C^0\left([t_1,t_2] \times \Omega\right)} \le
 C e^{-\lambda_J^{\alpha/2} t_1}.
\end{equation}
This tends to zero as $J \rightarrow \infty$, which shows that $S^{\ell, J=1}$ for any $\ell \ge 0$ is uniformly convergent for $t > 0$. Hence, for any $\ell$,
\begin{equation}
\label{term-by-term_valid}
\partial^\ell_t w(t,x) = S^{\ell, J=1} (t,x). 
\end{equation}
This, together with \eqref{final_S_sobolev_estimate} and \eqref{S_smooth_estimate}, gives parts (1) and (2) of Theorem \ref{steady_state}. 

To prove part (3) of Theorem \ref{steady_state}, we let $\epsilon > 0$. 
By the assumption $w_0 \in L^2(\Omega)$, there exists an integer $N$ such that 
\begin{equation}
\sum_{j = {N+1}}^\infty
(w_0, e_j)^2_{L^2(\Omega)}
\le 
\frac{\epsilon}{2}.
\end{equation}
Since $0 < e^{-t\lambda_j^{\alpha/2}} < 1$, we have
\begin{equation}
\sum_{j = {N+1}}^\infty
\left( 1 - e^{-t\lambda_j^{\alpha/2}}
\right)^2
(w_0, e_j)^2_{L^2(\Omega)}
\le 
\frac{\epsilon}{2}.
\end{equation}
For the same $N$, for all $t$ sufficiently small,
\begin{equation}
\sum_{j = 1}^N
\left( 1 - e^{-t\lambda_j^{\alpha/2}}
\right)^2
(w_0, e_j)^2_{L^2(\Omega)}
\le 
\frac{\epsilon}{2}.
\end{equation}
Therefore, 
\begin{align}
\|
w_0 - w(t,x)
\|^2_{L^2(\Omega)}
&=
\sum_{j = {1}}^N
\left( 1 - e^{-t\lambda_j^{\alpha/2}}
\right)^2
(w_0, e_j)^2_{L^2(\Omega)}
+ \sum_{j = {N+1}}^\infty
\left( 1 - e^{-t\lambda_j^{\alpha/2}}
\right)^2
(w_0, e_j)^2_{L^2(\Omega)}
\\
&\le
\epsilon 
\end{align}
for all $t$ sufficiently small.
This completes the proof of part (3) of Theorem \ref{steady_state}.  
Part (4) of that theorem is a consequence of parts (1) -- (3).

Using the regularity results in Theorem \ref{steady_state}, we can show
\begin{theorem}\label{survival_theorem}
Define, for $t \ge 0$ and $x \in \Omega$, the survival probability 
\begin{equation}
\label{survival_probability}
w(t,x) = \mathbb{P} \left\{ 
X^{\Omega, \alpha}_t \in \Omega
\right\} =
1-  \mathbb{P} \left\{ 
X^{\Omega, \alpha}_t \in \partial\Omega
\right\}.
\end{equation} 
Here, $X^{\Omega, \alpha}_t $ is a sample path of subordinate stopped
Brownian motion that begins at $x \in \Omega$. Then, $w(t,x)$ satisfies
\begin{align}
\label{survival_equation}
\begin{cases}
&\partial_t w+ (-\Delta_{\Omega,0})^{\alpha/2} w
 = 0 \text{ for $ x \in \Omega, t > 0$} \\
&w (x,t) = 0 \quad \text{ for $x \in \partial\Omega$} \\
&w (x,0) = 1 \quad \text{ for $x \in \Omega$}.
\end{cases}
\end{align}
This implies that the survival probability satisfies all of the estimates in 
Theorem \ref{steady_state}. 
\end{theorem}

Note that this result is precisely what would follow from the stochastic representation 
in Theorem \ref{theorem_omega_inhomo} for the solution of \eqref{survival_equation}:
\begin{equation}
w(t,x) = 
\mathbb{E}_{X^{\Omega,\alpha}_0 = x} 
\left[ 1 \cdot \chi_{\tau_\Omega > T_{\alpha/2}(t)}
+ 0 \cdot \chi_{\tau_\Omega \le T_{\alpha/2}(t)}\right]
=
\mathbb{P} \left\{ 
X^{\Omega, \alpha}_t \in \Omega
\right\}.
\end{equation}
The matter is that the theorem does not directly apply, since the initial 
condition on $\overline{\Omega}$ is not continuous (i.e., the limit of the initial condition at the boundary is not consistent with the boundary conditions). However, we will use a the
regularity results just obtained for the problem \eqref{survival_equation}
to prove the desired representation. 

Note that \eqref{survival_probability} clearly satisfies the initial condition and the boundary condition of problem \eqref{survival_equation}; the initial condition is satisfied since 
$X^{\Omega, \alpha}_{t=0} = x$, and the boundary condition is satisfied since for 
$x \in \partial \Omega$, we have $\tau_\Omega = 0$. Thus, it remains to show that the solution $w(t,x)$ of problem \eqref{survival_equation} satisfies
\eqref{survival_probability} for $t > 0, x \in \Omega$. By uniqueness, Theorem \ref{survival_theorem} will follow. 

Let $t, \epsilon > 0$. Put
\begin{equation}
w_{\epsilon}(t,x) = w(t+\epsilon, x), 
\end{equation}
where $w(t,x)$ is the solution of \eqref{survival_equation}. 
Then $w_{\epsilon}$ solves 
\begin{align}
\label{survival_equation_epsilon}
\begin{cases}
&\partial_t w_{\epsilon} + (-\Delta_{\Omega,0})^{\alpha/2} w_{\epsilon}
 = 0 \text{ for $ x \in \Omega$} \\
&w_{\epsilon} (x,t) = 0 \quad \text{ for $x \in \partial\Omega$} \\
&w_{\epsilon} (x,0) = w(\epsilon,x) \quad \text{ for $x \in \Omega$}.
\end{cases}
\end{align}
Since $w(t,x) \in C^{\infty}((0,\infty) \times \Omega)$, $w_{\epsilon}$ satisfies
\begin{equation}
w_{\epsilon}(t,x) = 
\mathbb{E}_{X^{\Omega,\alpha}_0 = x} 
\left[ w(\epsilon, X^{\Omega,\alpha}_t) \cdot \chi_{\tau_\Omega > T_{\alpha/2}(t)}
\right]
\end{equation}
In other words, for $t > 0$ and $x \in \Omega$,
\begin{equation}
\label{lets_take_the_limit}
w(t+\epsilon,x) = 
\mathbb{E}_{X^{\Omega,\alpha}_0 = x} 
\left[ w(\epsilon, X^{\Omega,\alpha}_t) \cdot \chi_{\tau_\Omega > T_{\alpha/2}(t)}
\right]
\end{equation}
For any $x$, as $\epsilon \rightarrow 0$, the left-hand side $w(t+\epsilon, x) \rightarrow w(t,x)$. To evaluate the limit of the right-hand side, consider
\begin{equation}
\label{expectation_representation}
\mathbb{E}_{X^{\Omega,\alpha}_0 = x} 
\left[1 \cdot \chi_{\tau_\Omega > T_{\alpha/2}(t)}
-
w(\epsilon, X^{\Omega,\alpha}_t) \cdot \chi_{\tau_\Omega > T_{\alpha/2}(t)}
\right].
\end{equation}
This may be written
\begin{equation}
\label{integral_representation}
\int_{\Omega}
\left[1 
-
w(\epsilon, y) 
\right]
P(t,x; y) dy
\end{equation}
where $P(t, x; y)$ is the transition density of $X^{\Omega, \alpha}_t$, giving the probability of starting at $x$ and hitting $y$:
\begin{equation}
P(t, x; y) dy = \mathbb{P} \left\{
X^{\Omega,\alpha}_t \in [y, y+dy]
\right\}, \quad X^{\Omega,\alpha}_0 = x.
\end{equation}
By construction, for $y \in \Omega$, $P(t,x; y)$ is less than the transition density of standard isotropic $\alpha$-stable L\'evy motion starting at $x$. The latter is bounded uniformly by some constant $C$ (depending on $\alpha$). Thus, we may bound 
\eqref{integral_representation} by 
\begin{align}
\int_{\Omega}
\big|
1 - w(\epsilon, y)
\big|
P(t,x; y) dy
\le
\|
1 - w(\epsilon, y)
\|_{L^2(\Omega)}
\sqrt{
\int_{\Omega}
P^2(t,x; y) dy
} \\
\le
\|
1 - w(\epsilon, y)
\|_{L^2(\Omega)}
C
\sqrt{
m(\Omega)
} 
\end{align}
Since $w(\epsilon, x) \rightarrow 1$ in $L^2(\Omega)$ as $\epsilon \rightarrow 0$, we obtain that
\eqref{expectation_representation} also tends to zero as $\epsilon \rightarrow 0$. In the same limit, equation \eqref{lets_take_the_limit} therefore becomes
\begin{equation}
w(t,x) =
\mathbb{E}_{X^{\Omega,\alpha}_0 = x} 
\left[1 \cdot \chi_{\tau_\Omega > T_{\alpha/2}(t)}
\right].
\end{equation}

\section{The Spectral Fractional Dirichlet Problem.}
\label{dirichlet}
\begin{theorem}
\label{elliptic_theorem}
Let $\Omega$ be a $C^\infty$ domain. 
Let $r \in C^2(\Omega)$ and $g$ (defined on $\partial\Omega$) be continuous. 
Then the solution to 
\begin{align}
\label{fractional_elliptic}
\begin{cases}
&(-\Delta_{\Omega,g})^{\alpha/2}u (x) = r(x) \text{ for $ x \in \Omega$} \\
&u(x) = g(x) \quad \text{ for $x \in \partial\Omega$} \\
\end{cases}
\end{align}
is given by
\begin{equation}
\label{fractional_elliptic_formula}
u(x) = \mathbb{E}_{X^{\Omega,\alpha}_0 = x} 
\left[  g\left(X^{\Omega,\alpha}_{T_{\alpha/2}^{-1}(\tau_\Omega)}\right)\right]
+ \mathbb{E}_{X^\alpha_0 = x}  \left[ \int_0^{T_{\alpha/2}^{-1}(\tau_\Omega)}  r(X^{\Omega,\alpha}_{s}) ds \right]
\end{equation} 
\end{theorem}

According to \eqref{harbir_splitting}, the solution $u$ to problem \eqref{fractional_elliptic} can be written as
$u = u_1 + u_2$ where $u_1$ solves
\begin{align}
\label{fractional_elliptic_1}
\begin{cases}
&(-\Delta_{\Omega,0})^{\alpha/2}u_1 (x) = r(x) \text{ for $ x \in \Omega$} \\
&u_1(x) = 0 \quad \text{ for $x \in \partial\Omega$} \\
\end{cases}
\end{align}
and $u_2$ solves
\begin{align}
\label{fractional_elliptic_2}
\begin{cases}
&\Delta u_2 (x) = 0 \text{ for $ x \in \Omega$} \\
&u_2(x) = g(x) \quad \text{ for $x \in \partial\Omega$} \\
\end{cases}
\end{align}
The classical theory of the Dirichlet problem \cite{john1982partial,gilbarg2015elliptic} implies $u_2 \in C^2(\Omega)$, so by 
the Feynman-Kac formula, it can be written
\begin{equation}
u_2(x) =
\mathbb{E}_{X^{\Omega,2}_0 = x} 
\left[  g(X^{\Omega,2}_{\tau_\Omega})\right]
=
\mathbb{E}_{X^{\Omega,\alpha}_0 = x} 
\left[  g\left(X^{\Omega,\alpha}_{T_{\alpha/2}^{-1}(\tau_\Omega)}\right)\right]
\end{equation}
on account of $X^{\Omega,2}_{\tau_\Omega} = X^{\Omega,\alpha}_{T_{\alpha/2}^{-1}(\tau_\Omega)}$.
This is the first term in \eqref{fractional_elliptic_formula}. It remains to show that
$u_1$ coincides with the second term. 

We begin by noting, from Corollary 3.6 in \cite{grubb2016regularity}, that 
the solution $u_1$ in \eqref{fractional_elliptic_1} is $C^2(\Omega)$. 
Next, we consider the equation
\begin{align}
\label{fractional_parabolic_ally}
\begin{cases}
&\partial_t w + (-\Delta_{\Omega,0})^{\alpha/2} w = r \text{ for $ x \in \Omega$} \\
&w(x,t) = 0 \quad \text{ for $x \in \partial\Omega$} \\
&w(x,0) = 0 \quad \text{ for $x \in \Omega$}
\end{cases}
\end{align}
The solution $w$ to this equation can be written $w = w_{\text{trans}} + u_1$, where $u_1$ is the same as in equation
\eqref{fractional_elliptic_1} and therefore the transient part $w_{\text{trans}}$ solves
\begin{align}
\label{fractional_parabolic_transient}
\begin{cases}
&\partial_t w_{\text{trans}} + (-\Delta_{\Omega,0})^{\alpha/2} w_{\text{trans}} = 0 \text{ for $ x \in \Omega$} \\
&w_{\text{trans}}(x,t) =  0 \quad \text{ for $x \in \partial\Omega$} \\
&w_{\text{trans}}(x,0) = -u_1(x) \quad \text{ for $x \in \Omega$}
\end{cases}
\end{align}
Since $u_1 \in C^2$, we have from Theorem \ref{steady_state} that $w_{\text{trans}}(t,\cdot) \in C^2$ for $t \ge 0$ and therefore $w(t,\cdot) \in C^2$ for $t \ge 0$. 
From Theorem \ref{theorem_omega_inhomo}, we obtain
\begin{equation}
w(t,x)
=
\mathbb{E}_{X^\alpha_0 = x}  \left[ \int_0^t  r(X^{\Omega,\alpha}_{s})  \chi_{\tau_\Omega > {T_{\alpha/2}(s)}}  ds \right]
\end{equation}  
Next, we show, for each $x \in \Omega$,
\begin{enumerate}
\item
$w(t,x) \rightarrow u_1(x)$ as $t \rightarrow \infty$; equivalently, $w_{\text{trans}}(t,x) \rightarrow 0$. This follows directly from Theorem \ref{steady_state}.

\item
As $t \rightarrow \infty$,
$$
\mathbb{E}_{X^\alpha_0 = x}  \left[ \int_0^t  r(X^{\Omega,\alpha}_{s})  \chi_{\tau_\Omega > {T_{\alpha/2}(s)}}  ds \right]
\rightarrow 
\mathbb{E}_{X^\alpha_0 = x}  \left[ \int_0^{{T_{\alpha/2}^{-1}(\tau_\Omega)}}  r(X^{\Omega,\alpha}_{s})  ds \right].
$$
\end{enumerate}
Items (1) and (2) imply 
\begin{equation}
u_1(x) = \left[ \int_0^{{T_{\alpha/2}^{-1}(\tau_\Omega)}}  r(X^{\Omega,\alpha}_{s})  ds \right],
\end{equation}
yielding the second term in \eqref{fractional_elliptic_formula} and completing the proof of the theorem. 
To prove (2), we write using the law of total expectation,
\begin{multline}
\mathbb{E}_{X^\alpha_0 = x}  \left[ \int_0^t  r(X^{\Omega,\alpha}_{s})  \chi_{\tau_\Omega > {T_{\alpha/2}(s)}}  ds \right] 
=
\mathbb{E}_{X^\alpha_0 = x}  \left[ \int_0^t
r(X^{\Omega,\alpha}_{s})  ds \ \bigg| \ {\tau_\Omega > {T_{\alpha/2}(t)}} \right] 
\mathbb{P}\left\{ 
\tau_{\Omega} > T_{\alpha/2}(t)
\right\}
\\
+
\mathbb{E}_{X^\alpha_0 = x}  \left[ \int_0^{T^{-1}_{\alpha/2}(\tau_\Omega)}  
r(X^{\Omega,\alpha}_{s})   ds \ \bigg| \ \tau_{\Omega} \le T_{\alpha/2}(t)\right]
\mathbb{P}\left\{ 
\tau_{\Omega} \le T_{\alpha/2}(t)
\right\}.
\end{multline}
From Theorem \ref{survival_theorem}, we have
$
\mathbb{P}\left\{ 
\tau_{\Omega} > T_{\alpha/2}(t)
\right\}
\sim e^{-\lambda_1^{\alpha/2}t}
$,
so the first term in the above equation tends to zero as $t \rightarrow \infty$. Consequently, in the second term,
$
\mathbb{P}\left\{ 
\tau_{\Omega} \le T_{\alpha/2}(t)
\right\}
\rightarrow 1
$
as $t \rightarrow \infty$. By martingale convergence, the second term tends to
\begin{equation}
\mathbb{E}_{X^\alpha_0 = x}  \left[ \int_0^{T^{-1}_{\alpha/2}(\tau_\Omega)}  
r(X^{\Omega,\alpha}_{s})   ds \right].
\end{equation}
This completes the proof. 
\begin{remark}
Note that to obtain $u_2 \in C^2$ in equation \eqref{fractional_elliptic_2},  it is only required that $\Omega$ satisfy, e.g., the exterior sphere condition \cite{john1982partial}.
The condition that $\Omega$ is $C^\infty$ is used in two places. First, to invoke the regularity results of article \cite{grubb2016regularity} (Corollary 3.6) that 
the solution $u_1$ in Eq. \eqref{fractional_elliptic_1} is $C^2(\Omega)$.
This assumption on $\Omega$ is required due the techniques utilized in that article; we do not believe it is essential for 
obtaining $u_1 \in C^2(\Omega)$. For example, 
regularity results in H\"older spaces for the same equation  \eqref{fractional_elliptic_1}  were obtained in \cite{caffarelli2016fractional} assuming that $\Omega$ is Lipschitz, but only for the spaces $C^{0,r}$. If regularity in higher order $C^{2,r}$ spaces were proven with weaker conditions on $\Omega$, then the requirement that $\Omega$ is $C^\infty $ in Theorem \ref{elliptic_theorem} could be relaxed accordingly.  We believe this to be possible, but are not aware of such results in the literature.  The second place where the smoothness of the domain is used is in the proof of Theorem \ref{steady_state}, to invoke elliptic regularity (Theorem \ref{elliptic_regularity}) without restriction on $k$ and obtain smoothness of the solution $u(t,x)$. Full smoothness of the solution to the Cauchy problem is not required in the above proof, and may also be relaxed in a way that depends on the dimension $d$. For simplicity, we have not done so. Numerical simulations in Section \ref{numerics} verify both the time-dependent and time-independent stochastic solution formulas in the unit square and unit cube in benchmark examples. 
\end{remark}

\newpage
\section{Numerical Examples.}
\label{numerics}

\subsection{Implementation Details.}
\begin{figure}
\includegraphics[width=0.75\linewidth]{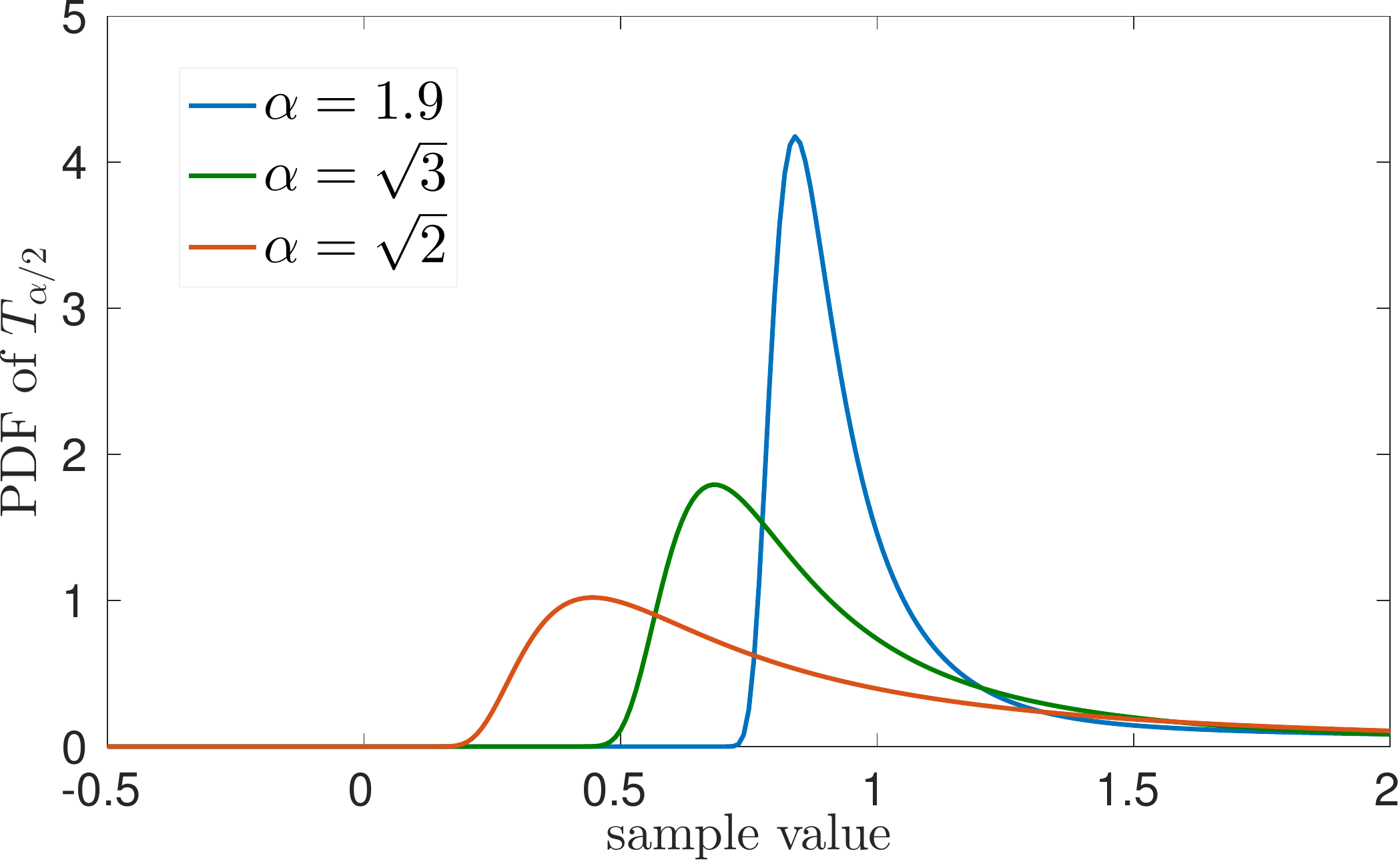}
   \caption{\small Probability density function of the standard stable subordinator $T_{\alpha/2}$ for various $\alpha$, with the other parameters fixed as discussed above.}
\label{subordinator_pdf_plots}
\end{figure}

\begin{figure}[htpb] 
  \label{fig7} 
  \begin{minipage}[b]{0.5\linewidth}
    \centering
    \includegraphics[width=1\linewidth]{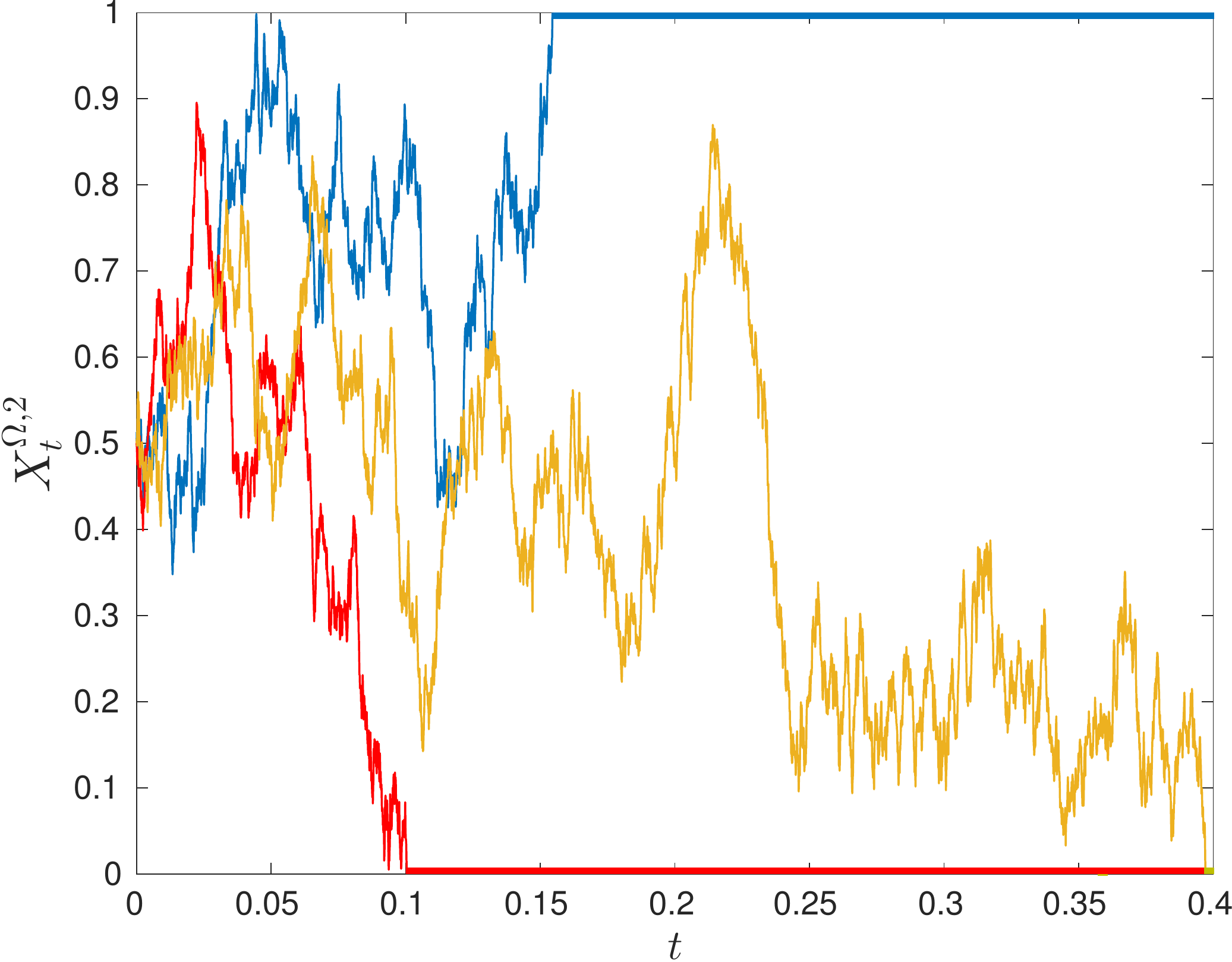} 
    \caption*{\small Original Brownian paths} 
    \vspace{4ex}
  \end{minipage}
  \begin{minipage}[b]{0.5\linewidth}
    \centering
    \includegraphics[width=1\linewidth]{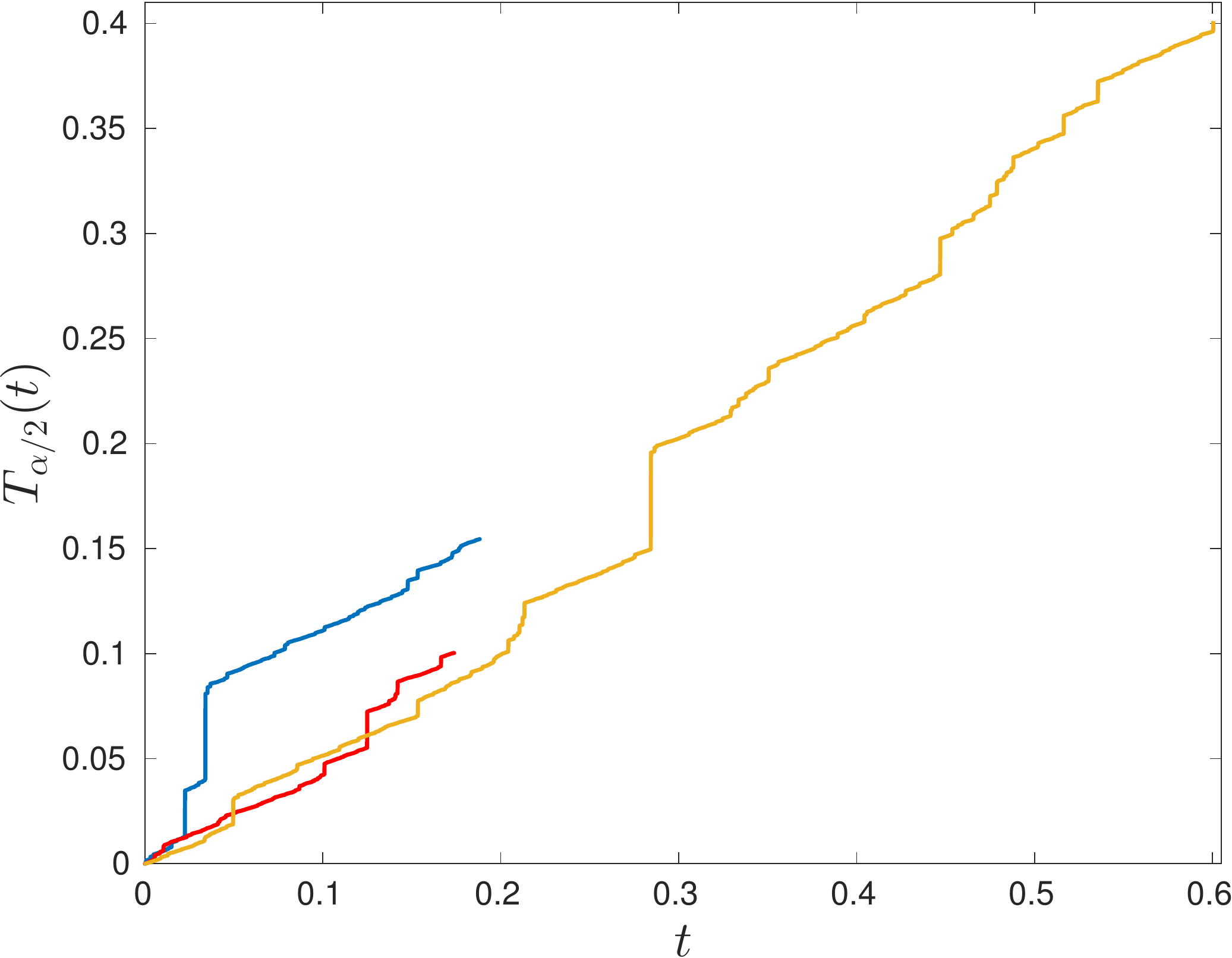} 
    \caption*{\small Subordinators, $\alpha = \sqrt{3}$} 
    \vspace{4ex}
  \end{minipage} 
  \begin{minipage}[b]{0.5\linewidth}
    \centering
    \includegraphics[width=1\linewidth]{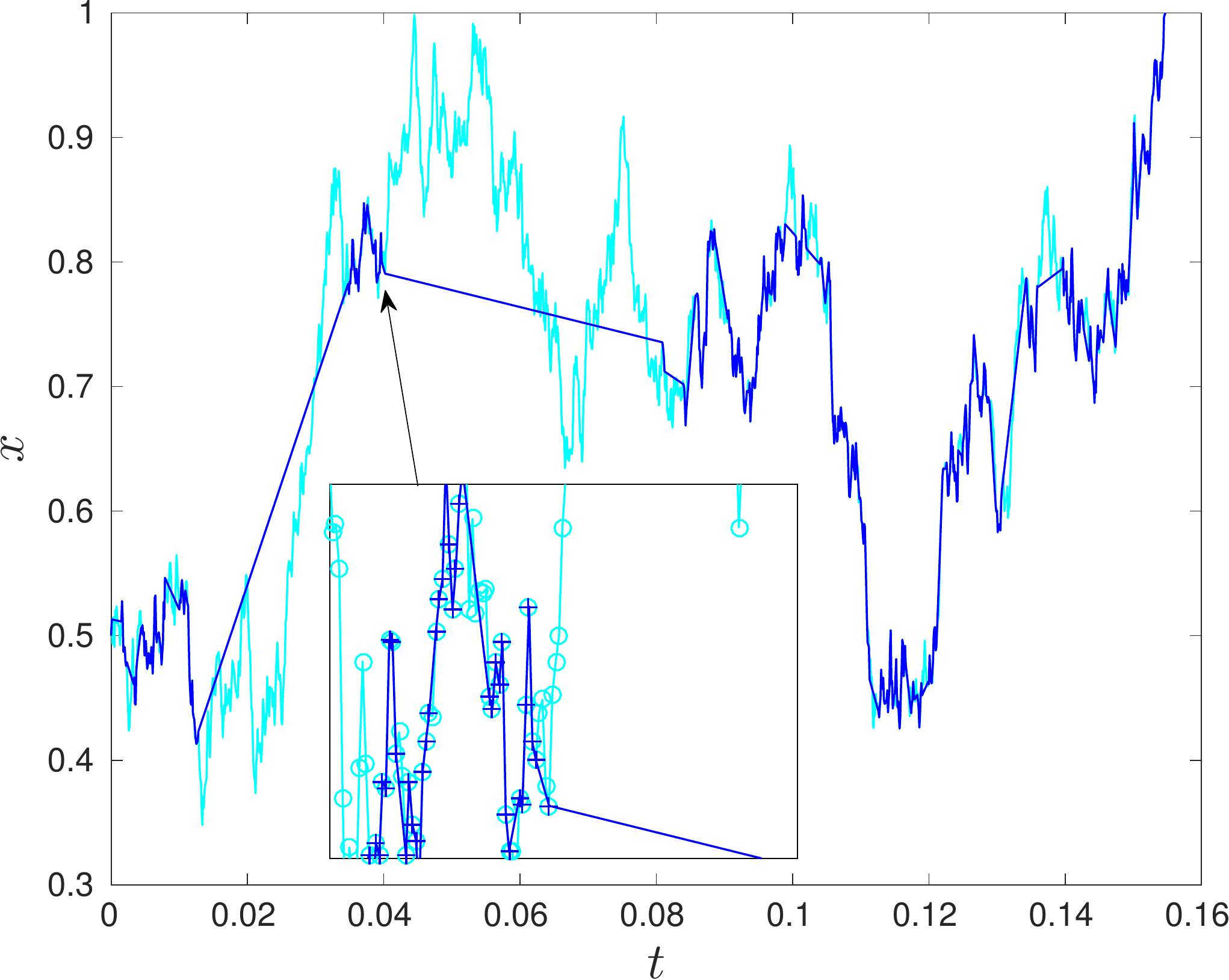} 
    \caption*{\small A subordinated path as sub-path of original path} 
    \vspace{4ex}
  \end{minipage}
  \begin{minipage}[b]{0.5\linewidth}
    \centering
    \includegraphics[width=1\linewidth]{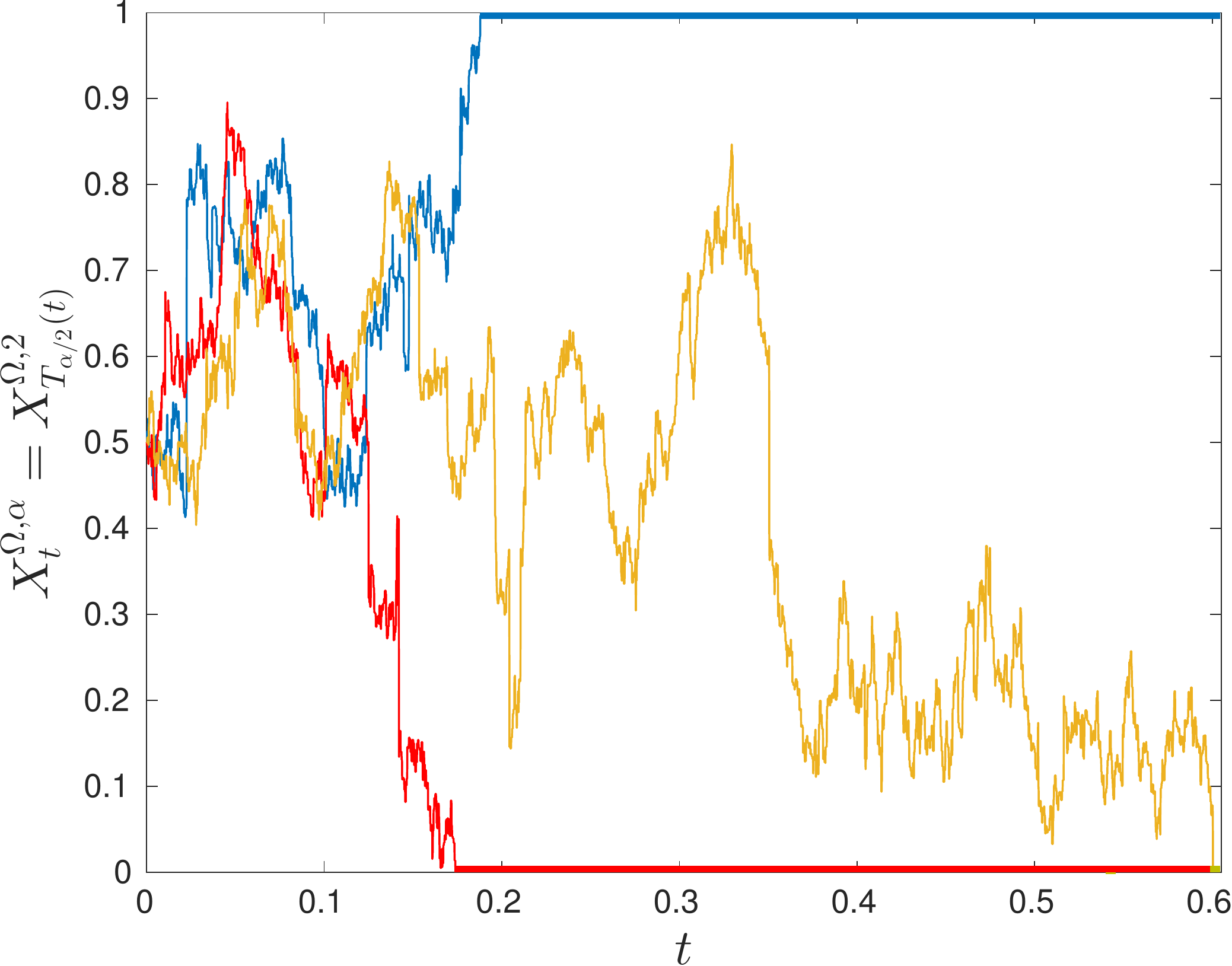} 
    \caption*{\small Final subordinated Paths from original BM paths} 
    \vspace{4ex}
  \end{minipage} 
\caption{\small \emph{Top left:} Three time series of stopped Brownian paths $X^{\Omega,2}$ in the unit interval $[0,1]$. \emph{Top right:} Three samples of the standard $\alpha$-stable subordinator $T_{\alpha/2}$ for $\alpha = \sqrt{3}$, run until the exit time of the respective stopped Brownian path shown to the left. 
\emph{Bottom Left:} The cyan process is an original BM (the blue one) from the top left plotted versus standard time $t$. The dark blue process shows the subordinated process plotted versus subordinated time $T_{\alpha/2}$. Some of the locations of the original process are skipped, and others are kept. 
\emph{Bottom right:} Final subordinated paths plotted versus standard time $t$. 
\label{1D_paths_illustration}}
\end{figure}

\begin{figure}[p] 
  \label{fig7} 
  \begin{minipage}[b]{0.5\linewidth}
    \centering
    \includegraphics[width=0.8\linewidth]{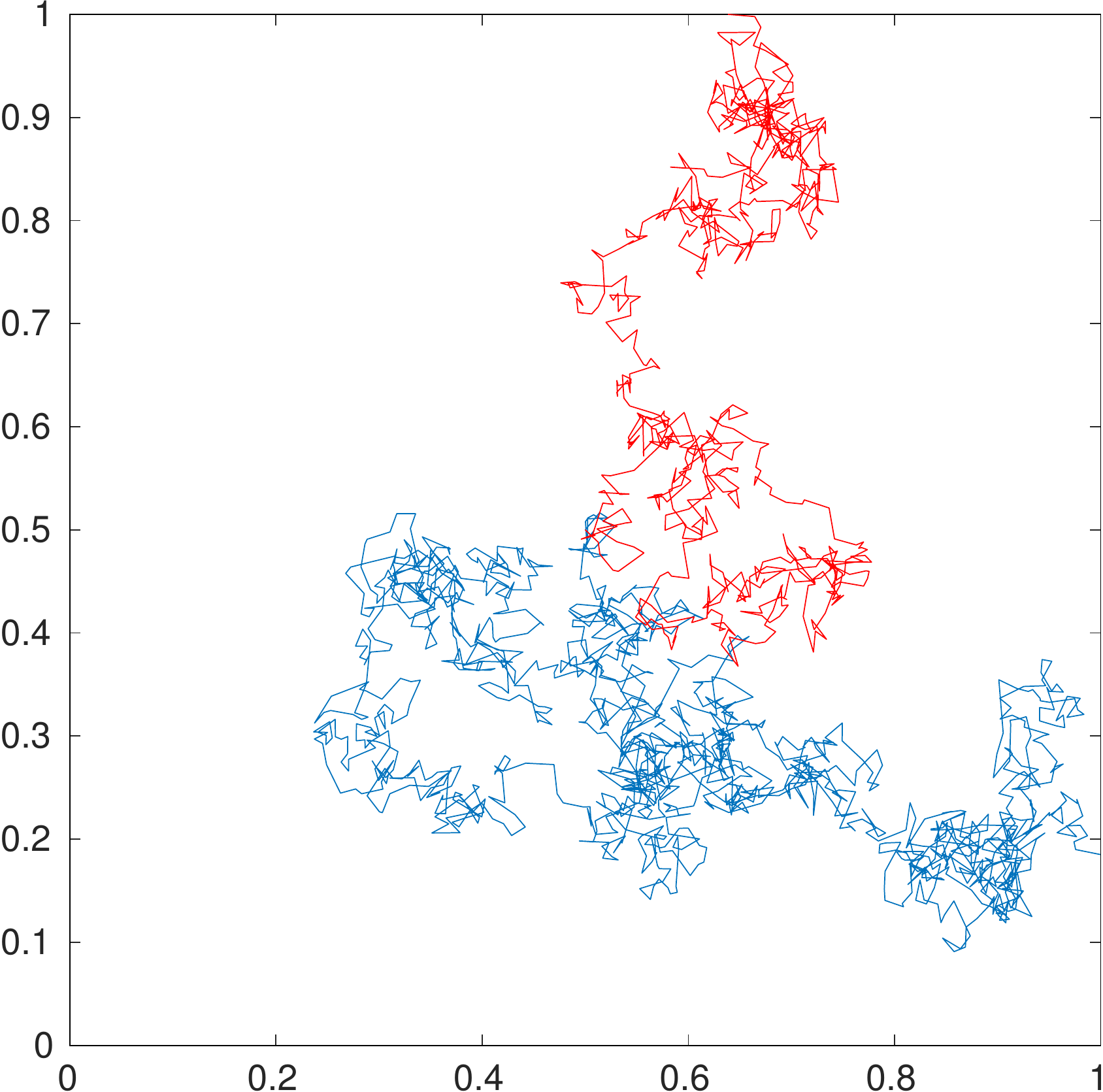} 
    \caption*{\small Stopped 2D Brownian paths in the unit square} 
    \vspace{4ex}
  \end{minipage}
  \begin{minipage}[b]{0.5\linewidth}
    \centering
    \includegraphics[width=0.8\linewidth]{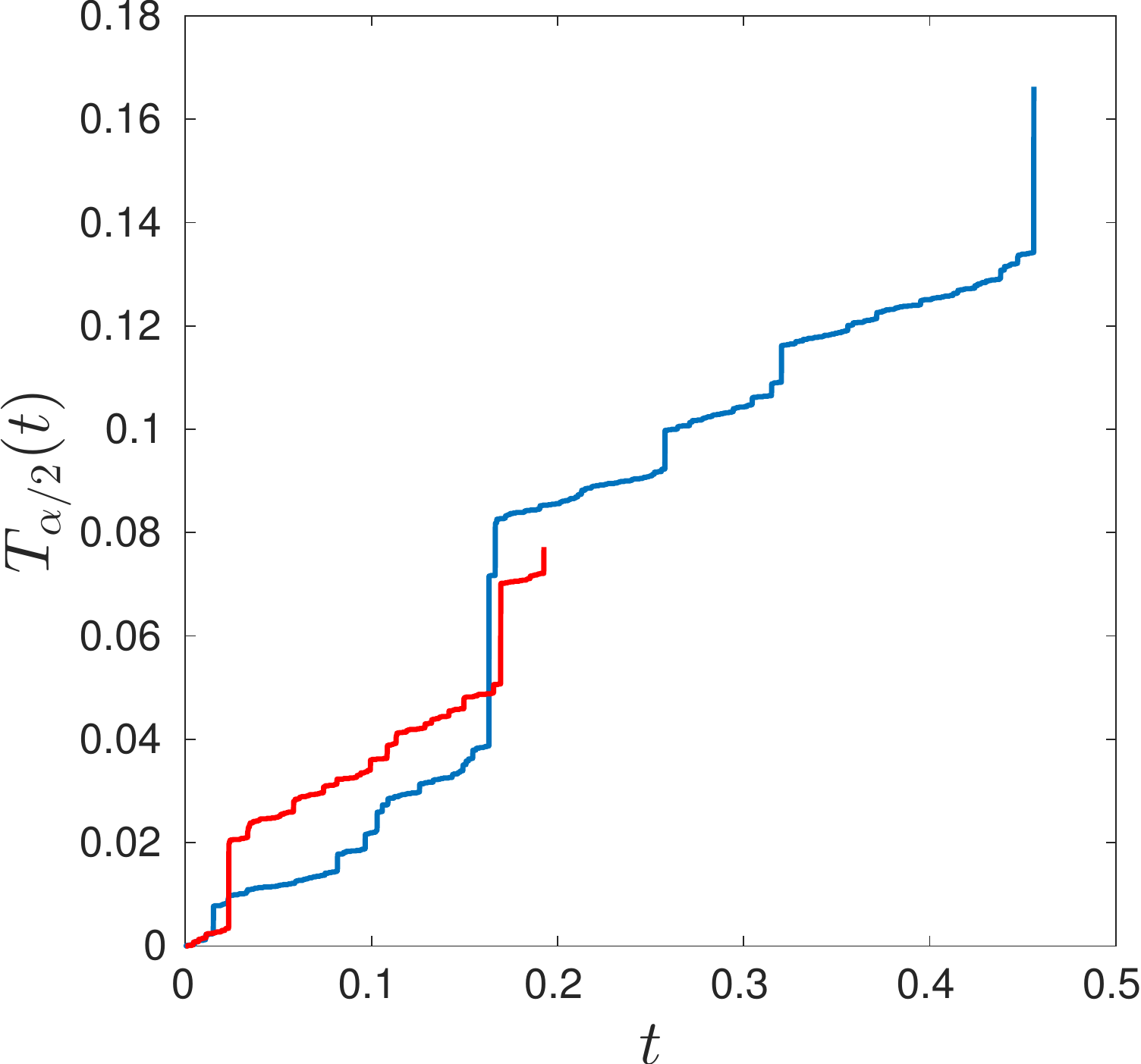} 
    \caption*{\small Subordinators, $\alpha = \sqrt{2}$} 
    \vspace{4ex}
  \end{minipage} 
  \begin{minipage}[b]{0.5\linewidth}
    \centering
    \includegraphics[width=0.8\linewidth]{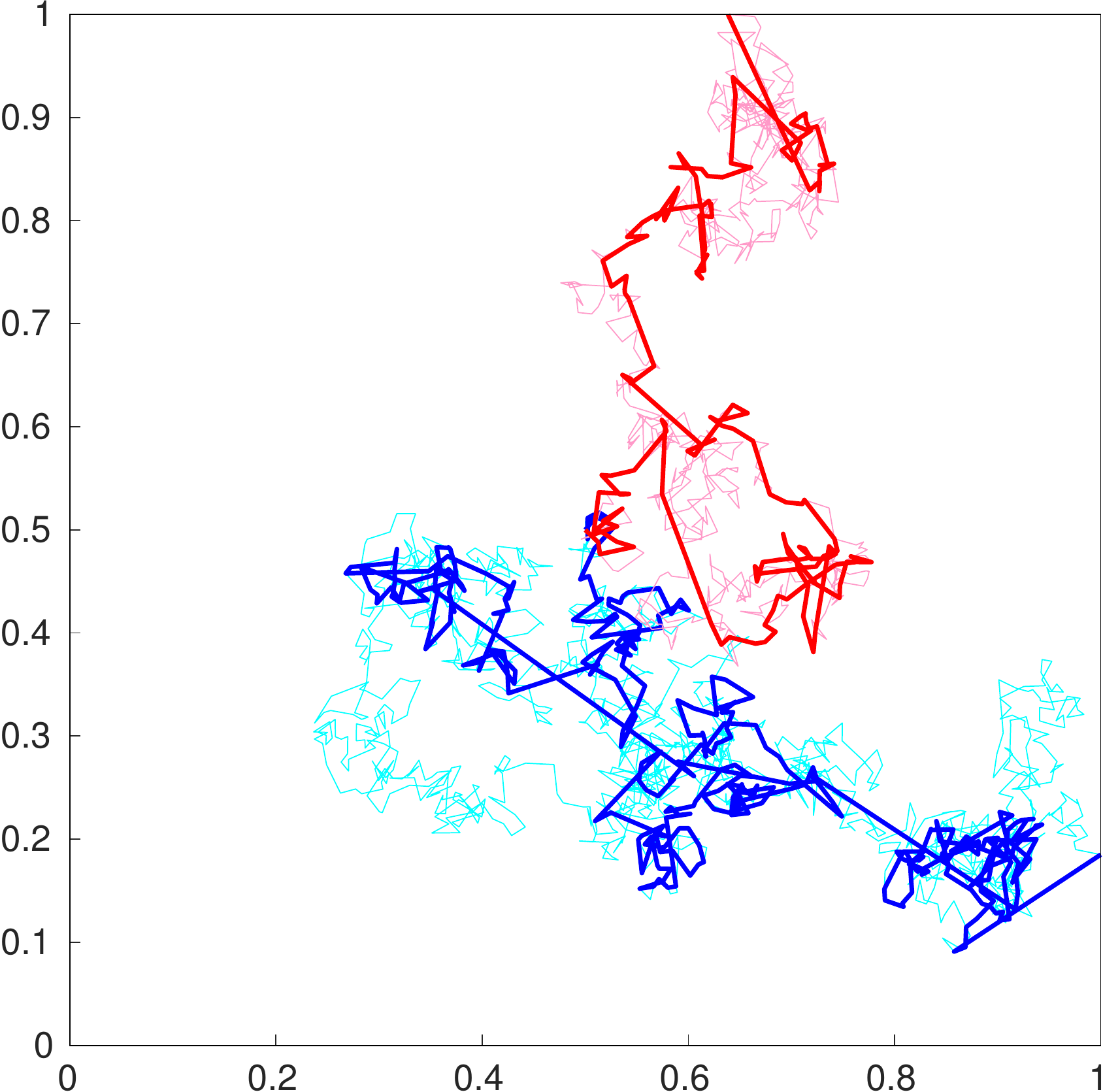} 
    \caption*{\small Subordinated paths as sub-paths of original paths} 
    \vspace{4ex}
  \end{minipage}
  \begin{minipage}[b]{0.5\linewidth}
    \centering
    \includegraphics[width=0.8\linewidth]{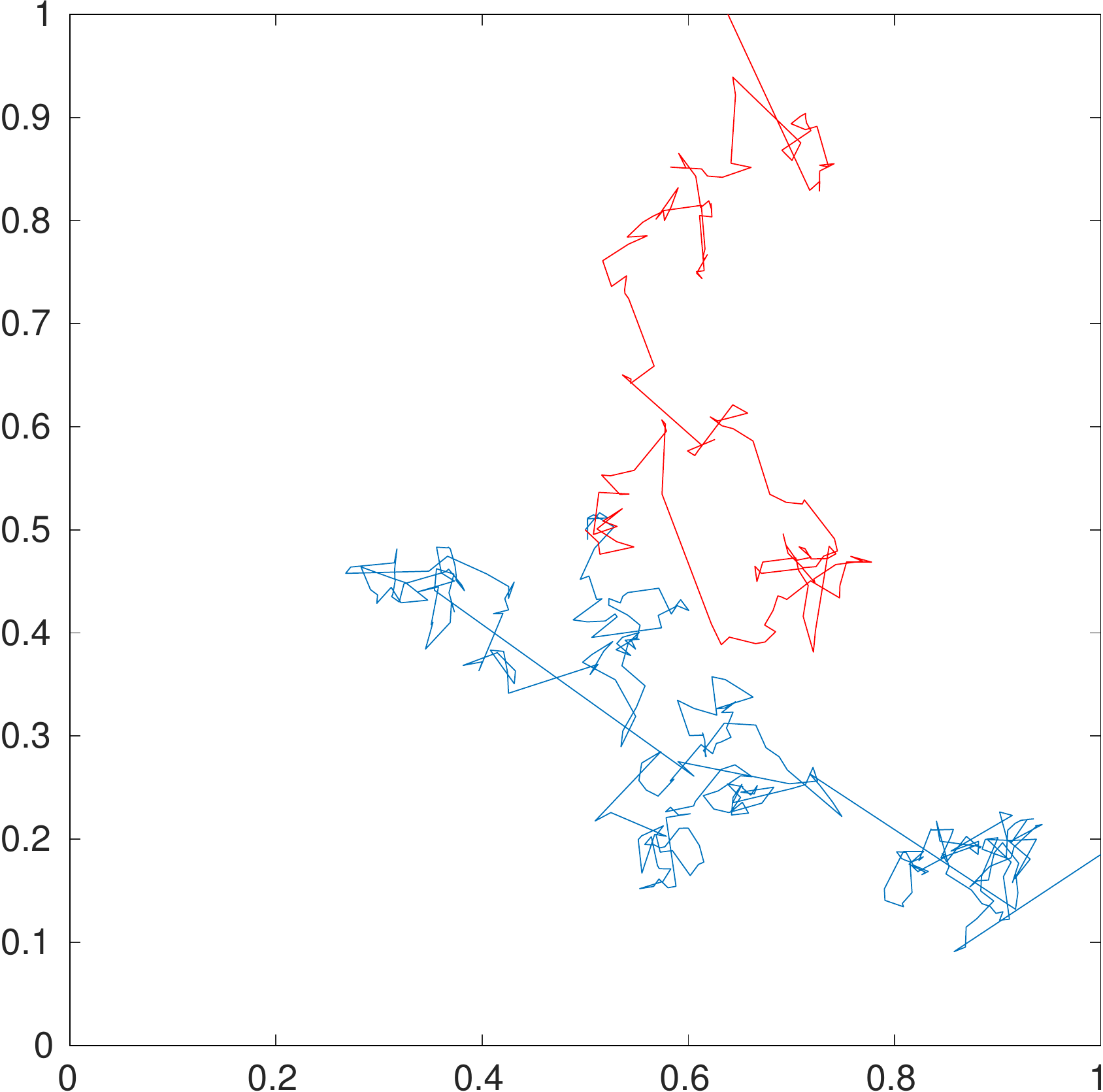} 
    \caption*{\small Final subordinated paths from original BM paths} 
    \vspace{4ex}
  \end{minipage} 
\caption{\small \emph{Top left:} Two stopped Brownian paths $X^{\Omega,2}_t$ in the unit square $[0,1] \times [0,1]$. \emph{Top right:} Two samples of the standard $\alpha$-stable subordinator $T_{\alpha/2}$ for $\alpha = \sqrt{2}$, run until the exit times of the respective stopped Brownian paths shown to the left are exceeded. 
\emph{Bottom Left:} The light-colored paths are original stopped BM from the top left; the dark-colored paths are subordinated paths $X^{\Omega,2}_{T_{\alpha/2}(t)}$. The subordinated process is seen to be a subprocess of the stopped BM. 
\emph{Bottom right:} Final subordinated paths.
\label{2D_paths_illustration}}
\end{figure}

\begin{figure}[p] 
  \label{fig7} 
  \begin{minipage}[b]{0.5\linewidth}
    \centering
    \includegraphics[width=0.8\linewidth]{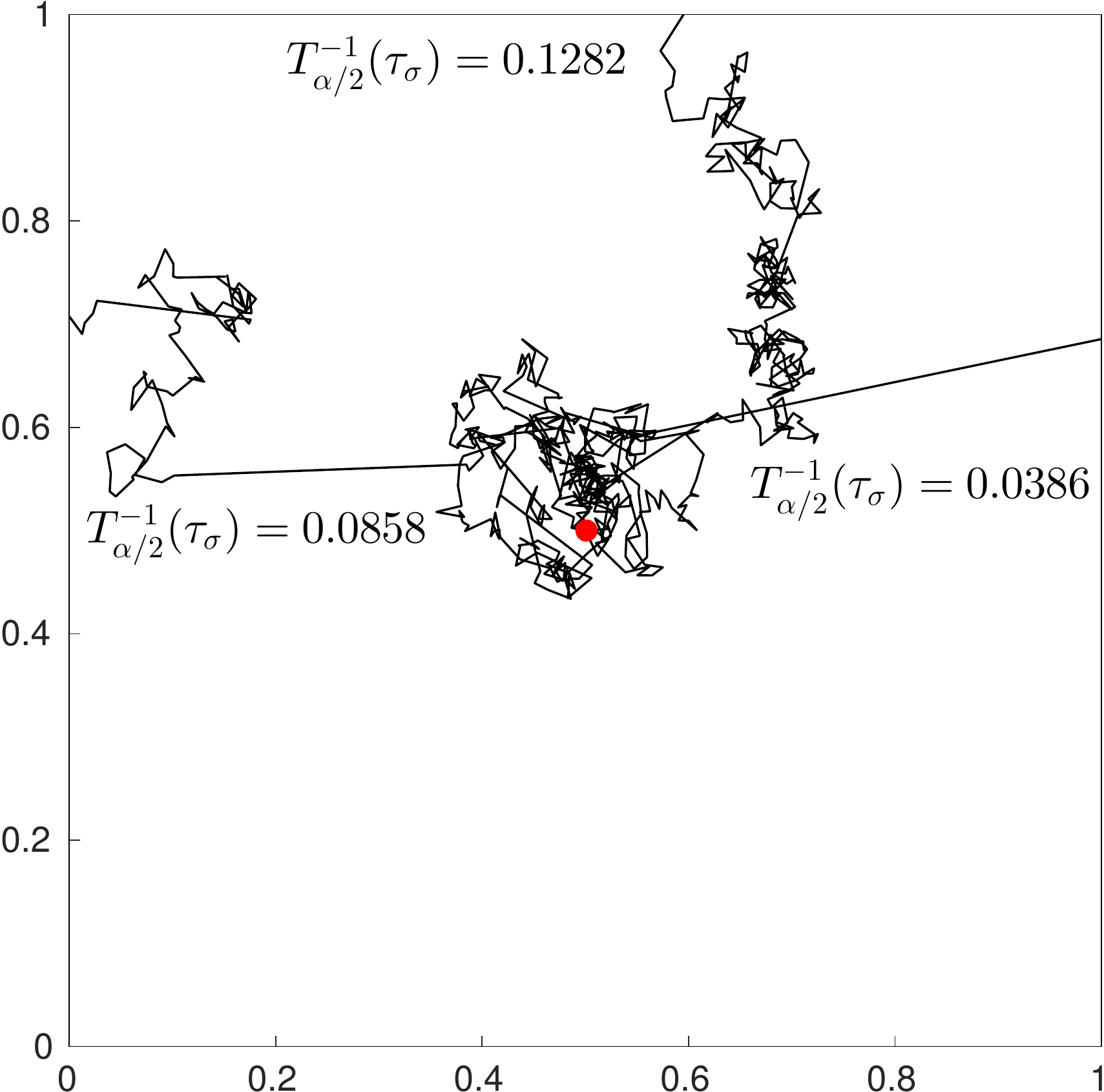} 
    \vspace{4ex}
  \end{minipage}
  \begin{minipage}[b]{0.5\linewidth}
    \centering
    \includegraphics[width=0.8\linewidth]{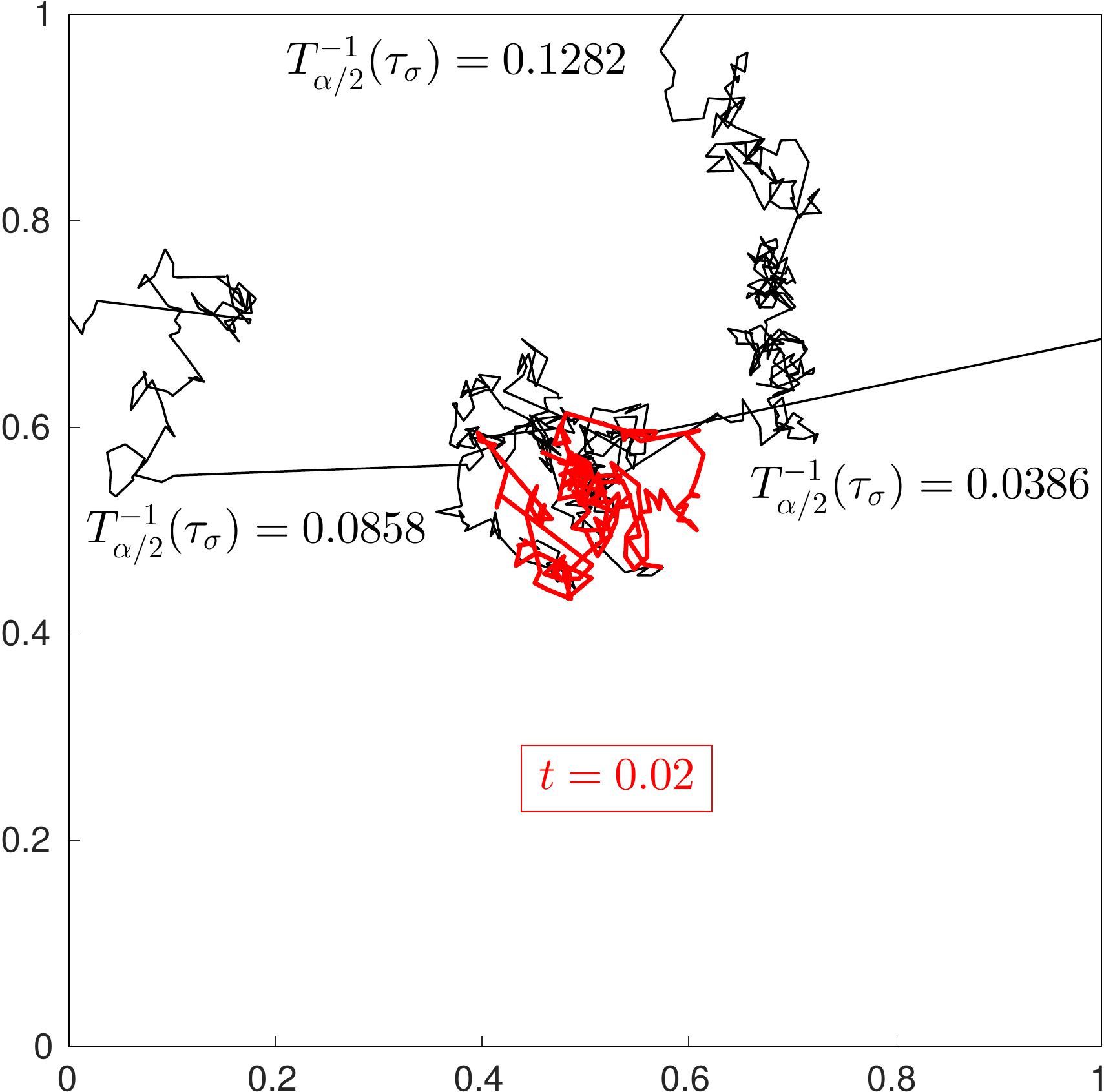} 
    \vspace{4ex}
  \end{minipage} 
  \begin{minipage}[b]{0.5\linewidth}
    \centering
    \includegraphics[width=0.8\linewidth]{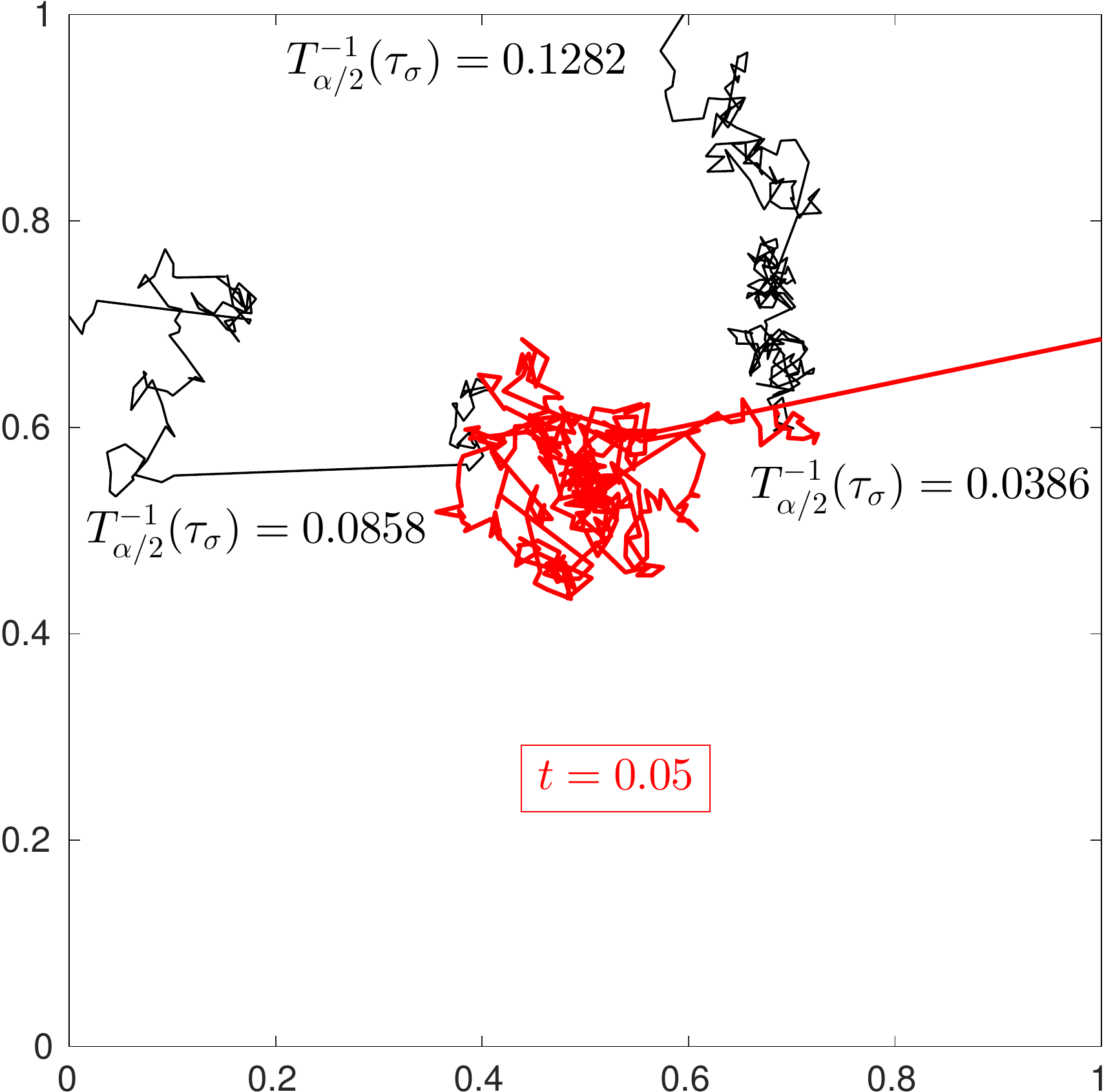} 
    \vspace{4ex}
  \end{minipage}
  \begin{minipage}[b]{0.5\linewidth}
    \centering
    \includegraphics[width=0.8\linewidth]{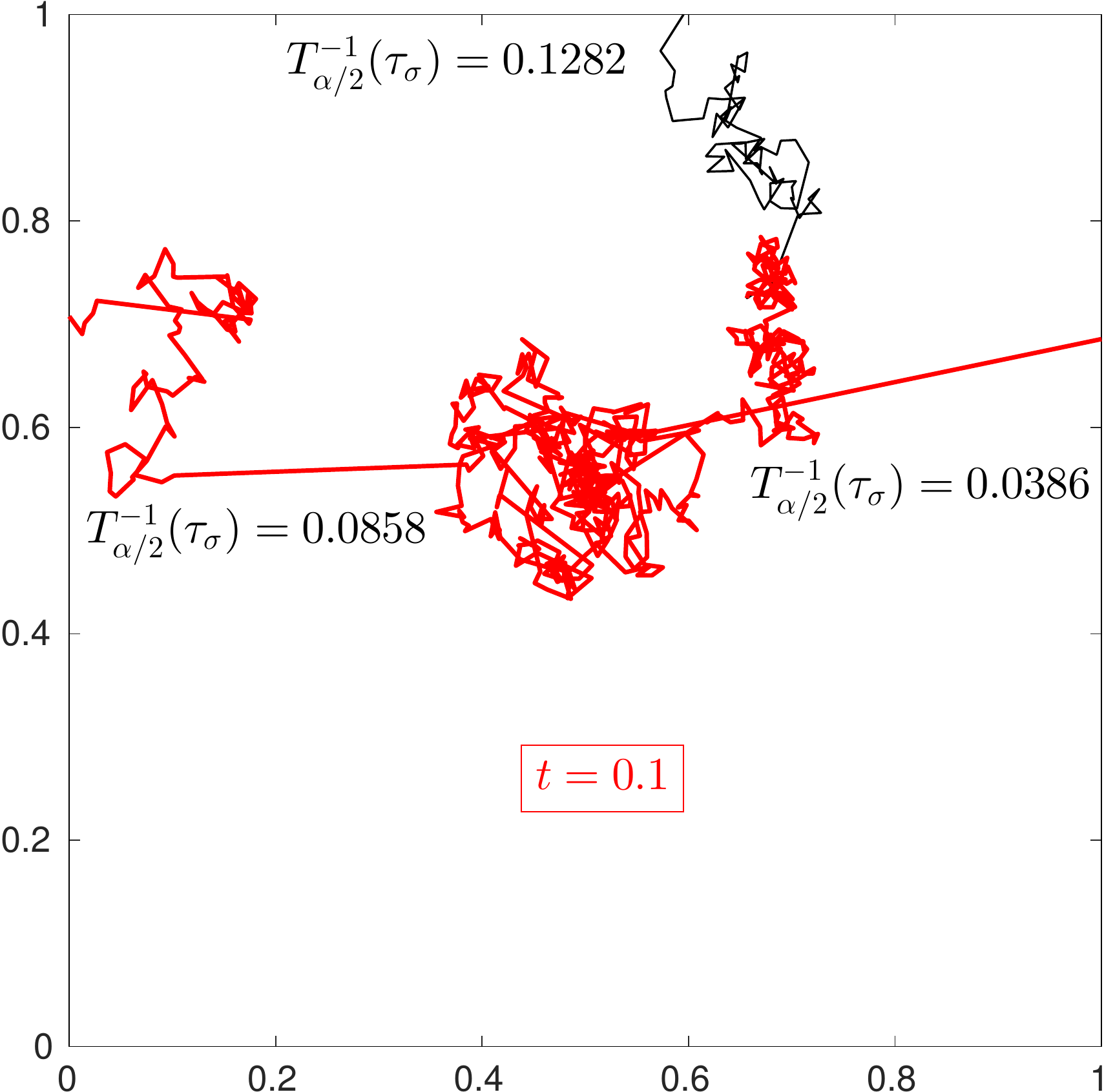} 
    \vspace{4ex}
  \end{minipage} 
\caption{\small Illustration of path integration in the unit square $[0,1] \times [0,1]$ for the stochastic solution in time. First, several paths of the subordinate stopped process are generated and stored in memory. \emph{Top left:} Three sample paths of the subordiante stopped process $X^{\Omega, \alpha}_t$ (black) for $\alpha = \sqrt{3}$ starting at $(x,y) = (0.5, 0.5)$ -- shown as a red dot -- with indicated stopping times. \emph{Top right, bottom left, bottom right:} The red tracings show the paths of integration (running along the stopped paths for specified time $t$) for computing the solution to the parabolic problem $u(t, x = 0.5, y = 0.5)$ at increasing $t$.}
\label{path_integration}
\end{figure}

We now verify the stochastic solution formulas proved in this article. Namely, for the parabolic time-dependent problem \eqref{problem_omega_inhomo} with initial condition $f$, boundary condition $g$, and time-independent right-hand side $r$, we study the formula
\begin{multline}
\label{parabolic_repeated}
u(t,x) = \mathbb{E}_{X^{\Omega,\alpha}_0 = x} 
\left[ f(X^{\Omega,\alpha}_t) \chi_{\tau_\Omega > {T_{\alpha/2}(t)}}
+ g(X^{\Omega,\alpha}_t) \chi_{\tau_\Omega \le {T_{\alpha/2}(t)}}\right] \\
+ \mathbb{E}_{X^\alpha_0 = x}  \left[ \int_0^t  r(X^{\Omega,\alpha}_s)  \chi_{\tau_\Omega > {T_{\alpha/2}(s)}}  ds \right].
\end{multline}
For the elliptic problem \eqref{fractional_elliptic} with boundary condition $g$ and right-hand side $r$, we study the formula 
\begin{equation}
\label{elliptic_repeated}
u(x) = \mathbb{E}_{X_0^{\Omega,\alpha} = x} \left[ g\left(X^{\Omega,\alpha}_{T^{-1}_{\alpha/2}(\tau_\Omega)}\right) \right]
+ \mathbb{E}_{X_0^{\Omega,\alpha} = x}  \left[ \int_0^{T^{-1}_{\alpha/2}(\tau_\Omega)}  r\left(X^{\Omega,\alpha}_{s}\right) ds \right].
\end{equation}
We start by discussing the direct discretization of the process $X^{\alpha,\Omega}_t \equiv X^{2,\Omega}_{T_{\alpha/2}(t)} = 
 X^2_{{T_{\alpha/2}(t)} \wedge \tau_\Omega}$.
\begin{enumerate}
\item
\textbf{Generation of Discrete Stopped Brownian Motion $X_t^{2,\Omega}$.} A starting point $X^{2,\Omega}_0 \in \Omega$ is specified. A discrete time step $dt$ is chosen, and we generate the path at times $t = 0, dt, 2dt, ...$ using
\begin{equation}
\label{bm_discretized}
X^{2,\Omega}_{t + dt} = 
X^{2,\Omega}_{t} + \sqrt{2} dt^{1/2} \texttt{mvnrnd}(\bm{0}, \text{\textbf{Id}}) 
\end{equation}
where $\texttt{mvnrnd}(\bm{0}, \text{\textbf{Id}})$ denotes the multivariate normal random variable with mean vector $\bm{0} = [0 \ 0 \ ... \ 0]^T$ and covariance matrix 
\begin{equation}
\text{\textbf{Id}} = 
  \left[ {\begin{array}{cccc}
   1 & 0 & ... & 0 \\
   0 & 1 & ... & 0 \\
\vdots & \vdots & ... & \vdots \\
0 & 0 & ... & 1
  \end{array} } \right]
\end{equation}
in dimension $d$ where $\Omega \subset \mathbb{R}^d$.
Equation \eqref{bm_discretized} is used to generate new points of the path while $X^{2,\Omega}_t \in \Omega$. Once we obtain a position $X^{2,\Omega}_{t + \Delta t} \not\in \Omega$, we replace this last position by the point $X^{2,\Omega}_{\tau_\Omega}$ nearest to the boundary, and end. The full path is stored and used in the next steps. The exit time $\tau_\Omega$ rounded up to the nearest $dt$, denoted $\lceil \tau_\Omega \rceil_{dt}$,  is saved as well.
\item
\textbf{Generation of Discrete Subordinator $T_{\alpha/2}(t)$ starting at zero.}
We generate the subordinator at times $t = 0, dt, 2dt, ...$ with $T_{\alpha/2}(0) = 0$ and
\begin{equation}
\label{subordinator_discretized}
T_{\alpha/2}(t+dt) = T_{\alpha/2}(t) + dt^{1/a} \text{ \texttt{stblrnd}}(a, \texttt{skewness},\texttt{scale}, \texttt{center})
\end{equation}
with parameters 
$a = \alpha/2$, $\texttt{skewness} = 1$, $\texttt{scale} = (\cos(\pi a/2))^{1/a}$, and $\texttt{center} = 0$.
For simplicity, we use the same value of $dt$ as used to discretize the Brownian motion, although this is not required. 
The parameter values of $\texttt{skewness}$, $\texttt{scale}$, and $\texttt{center}$ are the parameters in the Samorodnitsky and Taqqu parametrization of the one-dimensional stable process that yield the standard stable subordinator \cite{samoradnitsky2017stable}.
The function \texttt{stblrnd}, written by Mark Veillette \cite{stblrnd}, uses the methods of \cite{weron1995computer, chambers1976method}
to generate samples of the one-dimensional stable random variable in the Samorodnitsky and Taqqu parametrization. The resulting samples from $\text{ \texttt{stblrnd}}$ are positive; for these parameters, the stable distribution is supported on $\mathbb{R}^+$, as shown in Figure \ref{subordinator_pdf_plots}. Using \eqref{subordinator_discretized}, the subordinator is updated while $T_{\alpha/2}(t) < \tau_\Omega$.
When $T_{\alpha/2}(t+dt)$ is greater than or equal to $\tau_\Omega$ for the first time, we replace it by $\tau_\Omega$ and end. 
\item
\textbf{Subordination of Discrete Stopped Brownian Motion $X_t^{2,\Omega}$.} We now wish to insert the values of the subordinators $T_{\alpha/2}$ into the respective Brownian paths. However, since the positions of the Brownian paths are only available at times $0, dt, 2d,...$, we round the subordinators up to the nearest multiple of $dt$. We denote this by ${\lceil T_{\alpha/2} \rceil_{dt}}$. The discretized process  $X_t^{\alpha,\Omega}$ is obtained as $X_{\lceil T_{\alpha/2} \rceil_{dt}}^{2,\Omega}$. This step is illustrated in Table \ref{subordination_table}, and the process is illustrated with some example paths in Figures \ref{1D_paths_illustration} and \ref{2D_paths_illustration}.
\item
\textbf{Integration for parabolic solution \eqref{parabolic_repeated} in time}.
To calculate the first term in the solution at time $Ndt$, if a sample path of $X^{\Omega,\alpha}_{t}$ has hit the boundary at time $Ndt$, the boundary condition $g$ is evaluated at $X^{\Omega,\alpha}_{Ndt}$; otherwise, the initial condition $f$ is evaluated at $X^{\Omega,\alpha}_{N dt}$. The second term (the path integral) is discretized as
\begin{equation}
\label{discretization_of_solution}
\left[ \int_0^{N dt}  r(X^{\Omega,\alpha}_s)  \chi_{\tau_\Omega > {T_{\alpha/2}(s)}}  ds \right]
\sim
dt \sum_{n=1}^{N \wedge \lceil T_{\alpha/2}^{-1}(\tau_\Omega) \rceil_{dt} / dt} r(X^{\Omega,\alpha}_{ndt}).
\end{equation}
Importantly, the integration does not repeat over the exit point. 
Once a path hits the boundary, the path integral over that path no longer evolves. 
This is illustrated in Figure \ref{path_integration}.
\item
\textbf{Integration for elliptic solution.}
The boundary condition $g$ is evaluated at the endpoints of each path, and for the path integral, the same discretization \eqref{discretization_of_solution} is used except that all sums run to $\lceil T_{\alpha/2}^{-1}(\tau_\Omega) \rceil_{dt} / dt$, rather than $N \wedge \lceil T_{\alpha/2}^{-1}(\tau_\Omega) \rceil_{dt} / dt$. 
\end{enumerate}

\begin{remark}
Unlike in the classical $\alpha = 2$ case, computing the solution for the time-dependent problem (using the direct approach above) even for a short time $t$ requires generating Brownian paths all the way to the boundary. This is because the subordinator $T_{\alpha/2}$ may advance time significantly beyond the specified time $t$. The vast majority of computation time is spent generating paths, so once the paths are stored, we recommend computing the entire solution curve in time as well as the solution to the time independent problem.
\end{remark}

\begin{remark}
The true subordinator $T_{\alpha/2}(t)$ is strictly increasing. However, because we have to round the subordinator values to nearest multiple of $dt$, the discrete subordinator $\lceil T_{\alpha/2} \rceil_{dt}$ is merely nondecreasing, and frequently ``waits'' at the same time value for several increments.
\end{remark}

\begin{table}[htbp]
\def\arraystretch{1.25}
\begin{tabular}{|l|l|l|l|l|l|}
\cline{1-2} \cline{4-6}
$\bm{t}$                           & $\bm{X}^{\Omega,2}_t$                  &               & $\bm{t}$                                                 & $\bm{T}_{\alpha/2}(t)$                 & $\bm{X}^{\Omega,\alpha}_t$                 \\ \cline{1-2} \cline{4-6} 
$0$                                & $X_0$                                  &               & $0$                                                      & $0$                                    & $X_0$                                      \\ \cline{1-2} \cline{4-6} 
$dt$                               & $X_{dt}$                               &               & $dt$                                                     & $\lceil T_{\alpha/2}(dt) \rceil_{dt}$  & $X_{\lceil T_{\alpha/2}(dt) \rceil_{dt}}$  \\ \cline{1-2} \cline{4-6} 
$2dt$                              & $X_{2dt}$                              & $\rightarrow$ & $2dt$                                                    & $\lceil T_{\alpha/2}(2dt) \rceil_{dt}$ & $X_{\lceil T_{\alpha/2}(2dt) \rceil_{dt}}$ \\ \cline{1-2} \cline{4-6} 
$\vdots$                           & $\vdots$                               &               & $\vdots$                                                 & $\vdots$                               & $\vdots$                                   \\ \cline{1-2} \cline{4-6} 
$ndt$                              & $X_{ndt}$                              &               & $ndt$                                                    & $\lceil T_{\alpha/2}(ndt) \rceil_{dt}$ & $X_{\lceil T_{\alpha/2}(ndt) \rceil_{dt}}$ \\ \cline{1-2} \cline{4-6} 
$\vdots$                           & $\vdots$                               &               & $\vdots$                                                 & $\vdots$                               & $\vdots$                                   \\ \cline{1-2} \cline{4-6} 
$\lceil \tau_{\sigma} \rceil_{dt}$ & $X_{\lceil \tau_{\sigma} \rceil_{dt}}$ &               & $\lceil T^{-1}_\alpha(\tau_{\sigma}) \rceil_{dt} = N dt$ & $\lceil T_{\alpha/2}(Ndt) \rceil_{dt}$ & $X_{\lceil T_{\alpha/2}(Ndt) \rceil_{dt}}$ \\ \cline{1-2} \cline{4-6} 
$\vdots$                           & $X_{\lceil \tau_{\sigma} \rceil_{dt}}$ &               & $\vdots$                                                 & $\vdots$                               & $X_{\lceil T_{\alpha/2}(Ndt) \rceil_{dt}}$ \\ \cline{1-2} \cline{4-6} 
\end{tabular}
\vspace{0.1in}
\caption{\small \emph{Left:} The output of the Step 1 discussed above, in which a discrete stopped Brownian motion path is generated and stored. The process remains at the boundary point $X_{\lceil \tau_\Omega \rceil_{dt}}$ for all time greater than ${\lceil \tau_\Omega \rceil_{dt}} $, so the path need only be stored up to that point. \emph{Right:} Illustration of Steps 2 and 3, in which the discrete subordinator ${\lceil T_{\alpha/2}(Ndt) \rceil_{dt}}$ is generated (center column) and used to subordinate a discrete stopped Brownian motion paths stored in Step 1 to yield discrete subordinate stopped Brownian motion path. 
\label{subordination_table}}
\end{table}

\FloatBarrier

\FloatBarrier
\subsection{Two Dimensional Benchmark (Unit Square).}
First, we test the parabolic solution formula. 
Consider the two-dimensional unit square $R = [0,1] \times [0,1]$. Define 
$g: \partial R \rightarrow \mathbb{R}$ by
\begin{align}
\label{e:bdy_g_2d_square}
\begin{cases}
g(x, y = 0) &= 0 \\
g(x, y = 1) &= 1 \\
g(x = 0, y) &= y \\
g(x = 1, y) &= y
\end{cases} 
\end{align}
and consider the problem, for $\alpha = \sqrt{3}$,
\begin{align}
\label{2d_problem_1}
\begin{cases}
\partial_t u + (-\Delta_{R, g})^{\alpha/2} u &= \sin(\pi x) \sin( \pi y) \\
u \big|_{\partial R} &= g \\
u(t = 0, x) &= y + \frac{1}{(\pi^2 + \pi^2)^{\alpha/2}} \sin(\pi x) \sin(\pi x).
\end{cases}
\end{align}
In this problem, 
\begin{itemize}
\item
the boundary condition (BC) function $g$ is defined by \eqref{e:bdy_g_2d_square}. 
\item
the right-hand side (RHS) is
$r(x,y) = \sin(\pi x) \sin(\pi y)$, and 
\item
the initial condition (IC) is $f(x,y) = y + \frac{1}{(\pi^2 + \pi^2)^{\alpha/2}}  \sin(\pi x) \sin(\pi y) 
+ \sin(2 \pi x) \sin(2 \pi y)$. 
\end{itemize}
The exact solution is 
\begin{equation}
u(t,x) = y + \frac{1}{(\pi^2 + \pi^2)^{\alpha/2}} \sin(\pi x) \sin(\pi y)
+ e^{-\left(4\pi^2+4\pi^2\right)^{\alpha/2} t} \sin(2 \pi x) \sin(2 \pi y).
\end{equation}

\begin{figure}
\centering
\begin{subfigure}[b]{1\textwidth}
   \includegraphics[width=1\linewidth]{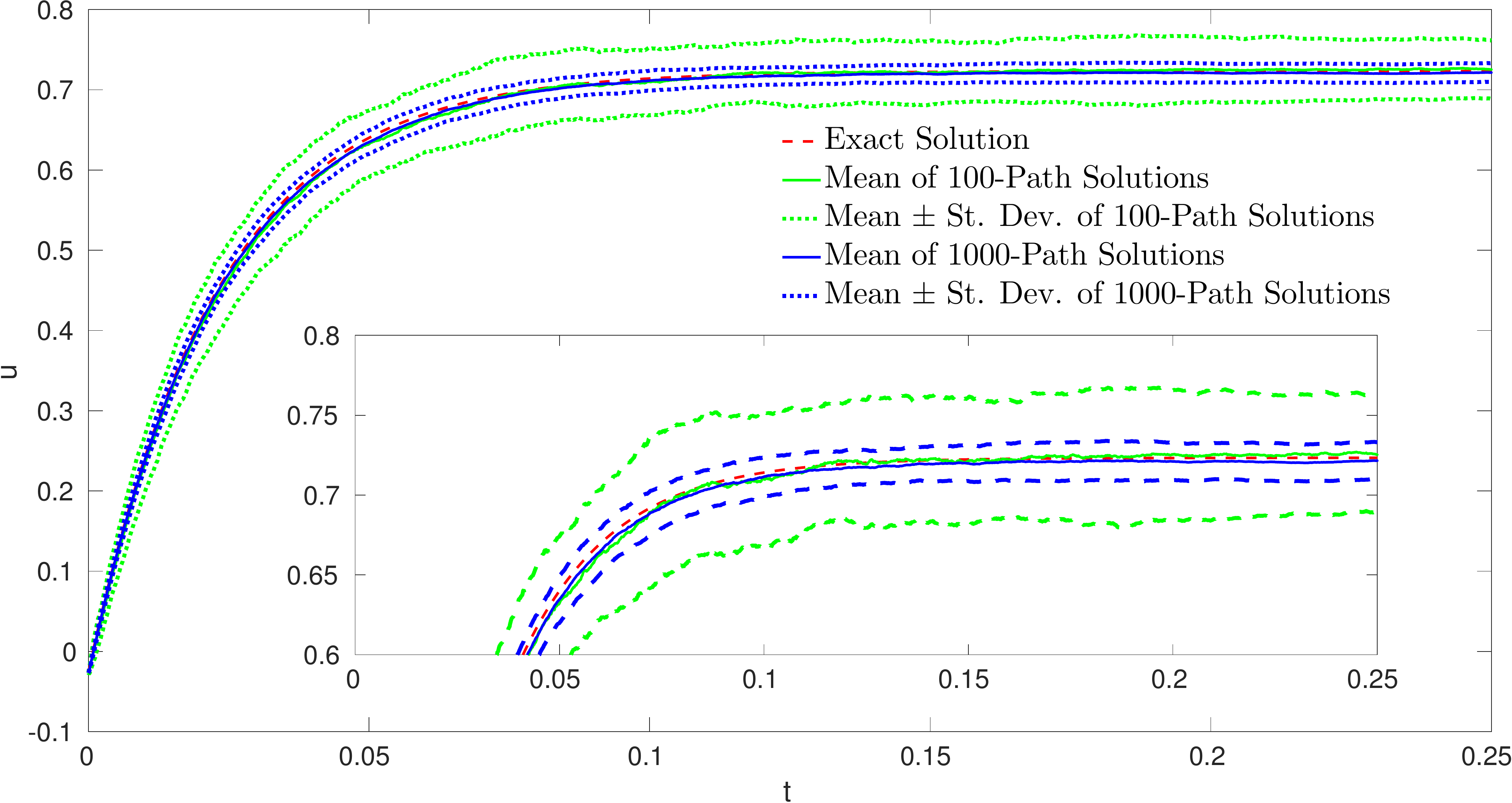}
   \label{fig:Ng1} 
\end{subfigure}

\begin{subfigure}[b]{1\textwidth}
   \includegraphics[width=1\linewidth]{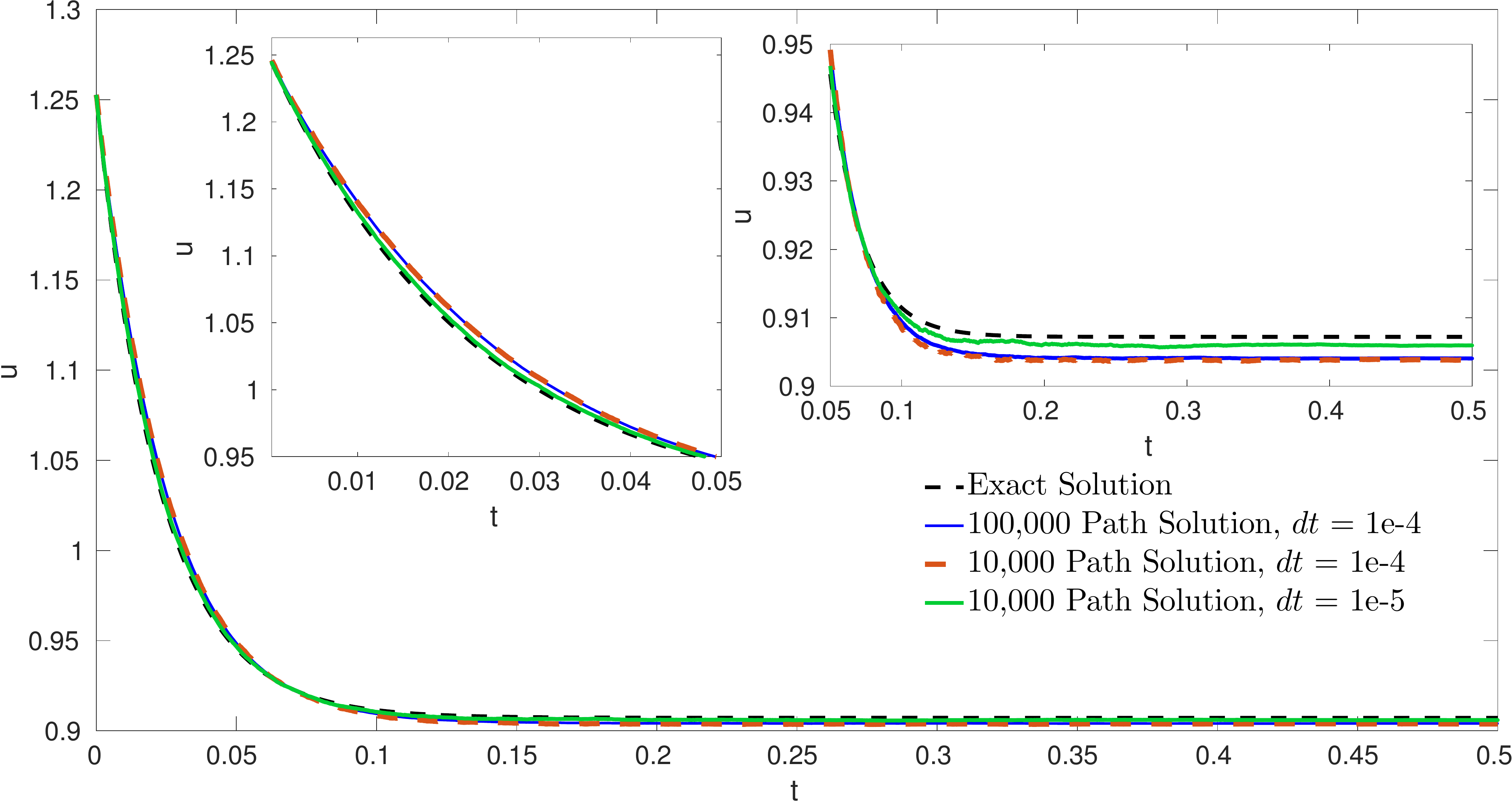}
   \label{fig:Ng2}
\end{subfigure}
\caption[3d_1o3]{
Convergence of stochastic solutions for the problem \eqref{2d_problem_1}.
\emph{Top, main:}
Convergence with respect to number of paths of the stochastic solution at (1/3, 2/3), with $dt = 10^{-4}$.  
\emph{Top, inset:}
Zoomed-in view of the main plot. 
\emph{Bottom, main:}
Bottlenecking of the convergence of the stochastic solution at (9/10,9/10), and improvement with refinement of $dt$ from $10^{-4}$ to $10^{-5}$.  
\emph{Bottom, insets:}
Zoomed-in views of the main plot.
}
\label{fig:2d_time_curves}
\end{figure}

In the top of Figure \ref{fig:2d_time_curves}, we consider the time trajectory of the stochastic solution $u(t, 1/3, 2/3)$ at the fixed point $(1/3, 2/3) \in R$. The discretized Brownian motion/subordinator time step is fixed at $dt = 10^{-4}$. The mean of 100 stochastic solutions, computed using 100 or 1,000 paths (as indicated) starting from $(1/3, 2/3)$, is plotted and compared to the exact solution. The``Mean $\pm$ Standard Deviation'' illustrates the expected variation of a stochastic solution computed using 100 or 1,000 paths, respectively, as it oscillates about the respective mean. Further, the means of the 100-path and 1000-path solutions \emph{also} represent stochastic solutions for $u(t,1/3, 2/3)$ computed using $100 \times 100 = 10,000$ paths and $100 \times 1000 = 100,000$ paths, respectively. We see that both are well-converged to the exact solution curve, although the 100,000 path solution is smoother when viewed more closely in the inset.

The next example illustrates the role of the Brownian motion/subordinator time step $dt$ when the stochastic solution formula is used near the boundary.
At the bottom of Figure \ref{fig:2d_time_curves}, we consider the time trajectory of the stochastic solution $u(t, 9/10, 9/10)$ at the fixed point $(9/10, 9/10)$ near the top-right corner of $R$. In both the main figure and the insets, the 10,000 path solution using $dt = 0.0001$ has significant discrepancy from the true solution. Increasing the number of paths (as in the top figure) to $100,000$ does not improve the accuracy the solution; in fact, it remains mostly unchanged from the 10,000 path solution. However, in this case, refining $dt$ by decreasing it one order of magnitude to $10^{-5}$ improves the solution significantly, even when using only 10,000 paths. 
This can be attributed to the closeness of the point $(9/10, 9/10)$ to the boundary. Because the paths have a high probability of exiting very close to the starting point in a few number of steps, using too high a $dt$ results in under-resolved Brownian motion that looks more like ballistic motion, bottonecking the convergence of the solution.

\begin{figure}[htpb]
   \includegraphics[width=1\linewidth]{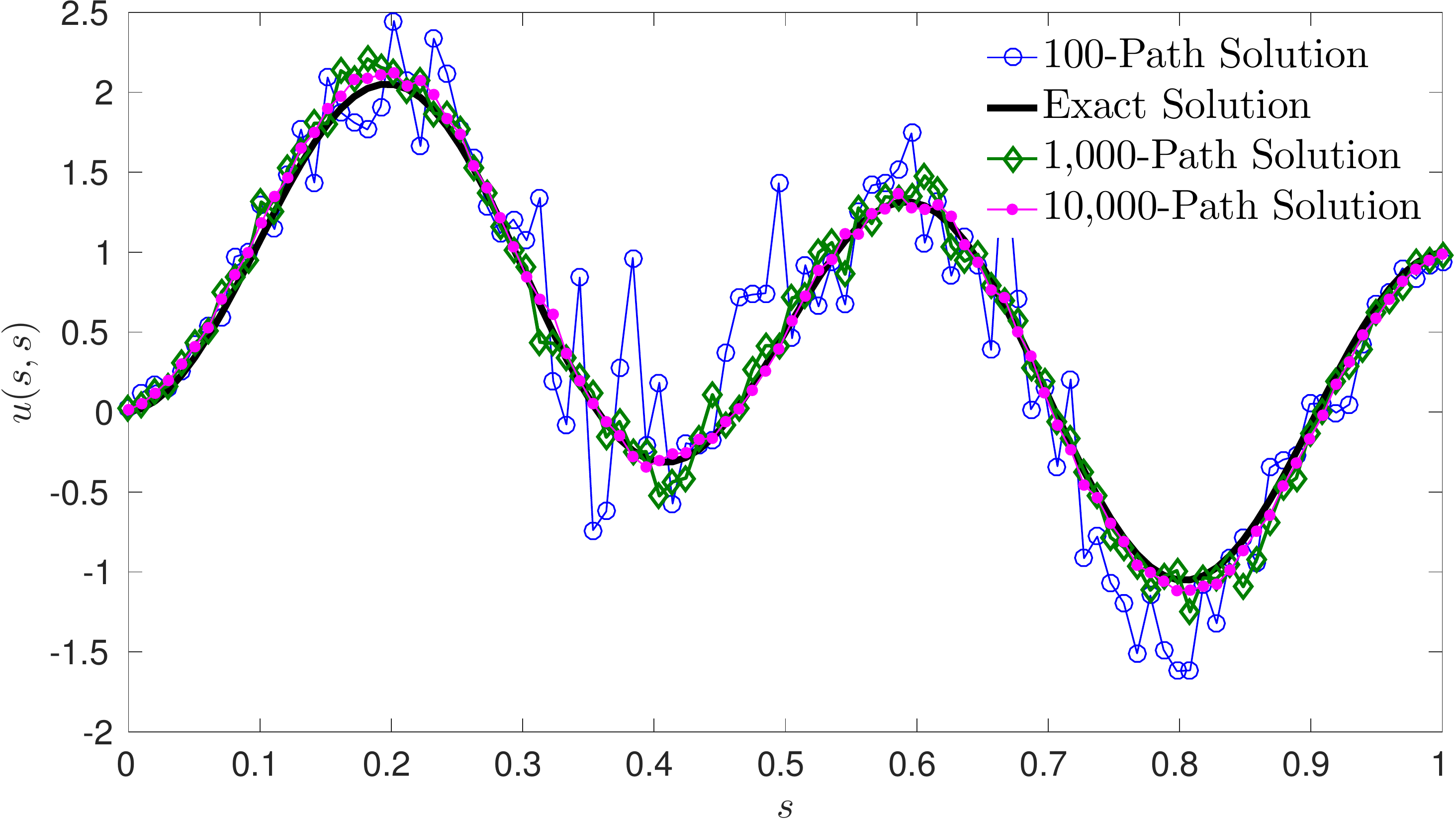}
   \caption{Stochastic solution of problem \eqref{2d_problem_2} in the 2D unit square at 100 equispaced points along the line $D = \{(s,s) \ \big| \ s \in [0,1]\}$, using 100, 1,000, and 10,000 paths starting from each point.}
   \label{fig:2d_elliptic} 
\end{figure}

\newpage
Next, we consider the elliptic problem on the unit square $R = [0,1] \times [0,1]$:
\begin{align}
\label{2d_problem_2}
\begin{cases}
(-\Delta_{R, g})^{\alpha/2} u &= 137 \sin(2\pi x)\sin(3\pi y) \\
u \big|_{\partial R}(x,y) &= y
\end{cases}
\end{align}
In this problem, the BC function $g = y$ as before and the RHS is $r(x) = 137 \sin(2\pi x)\sin(3\pi y)$. The exact solution is 
\begin{equation}
u(x,y) =  y + \frac{137}{\left[(2\pi)^2+(3\pi)^2\right]^{\alpha/2}}\sin(2\pi x)\sin(3\pi y);
\end{equation}
For $\alpha = \sqrt{3}$, we study the elliptic solution formula along the line diagonal line $D = \{(s,s) \ \big| \ s \in [0,1]\} \subset R$. Thus, $D$ runs from the bottom left corner of the square to the top right corner. We take 100 equispaced points on this line, and compute the solution at each of the points using 100, 1,000, and 10,000 paths. The discretized Browian motion/subordinator time step is fixed at $dt = 10^{-4}$. The solutions are compared to the exact solution in Figure \ref{fig:2d_elliptic}. We see that they converge to the true solution as the number of paths is increased. 

\FloatBarrier
\newpage
\subsection{Three Dimensional Benchmark (Unit Cube).}
\begin{figure}[p]
\centering
\begin{subfigure}[b]{1\textwidth}
   \includegraphics[width=1\linewidth]{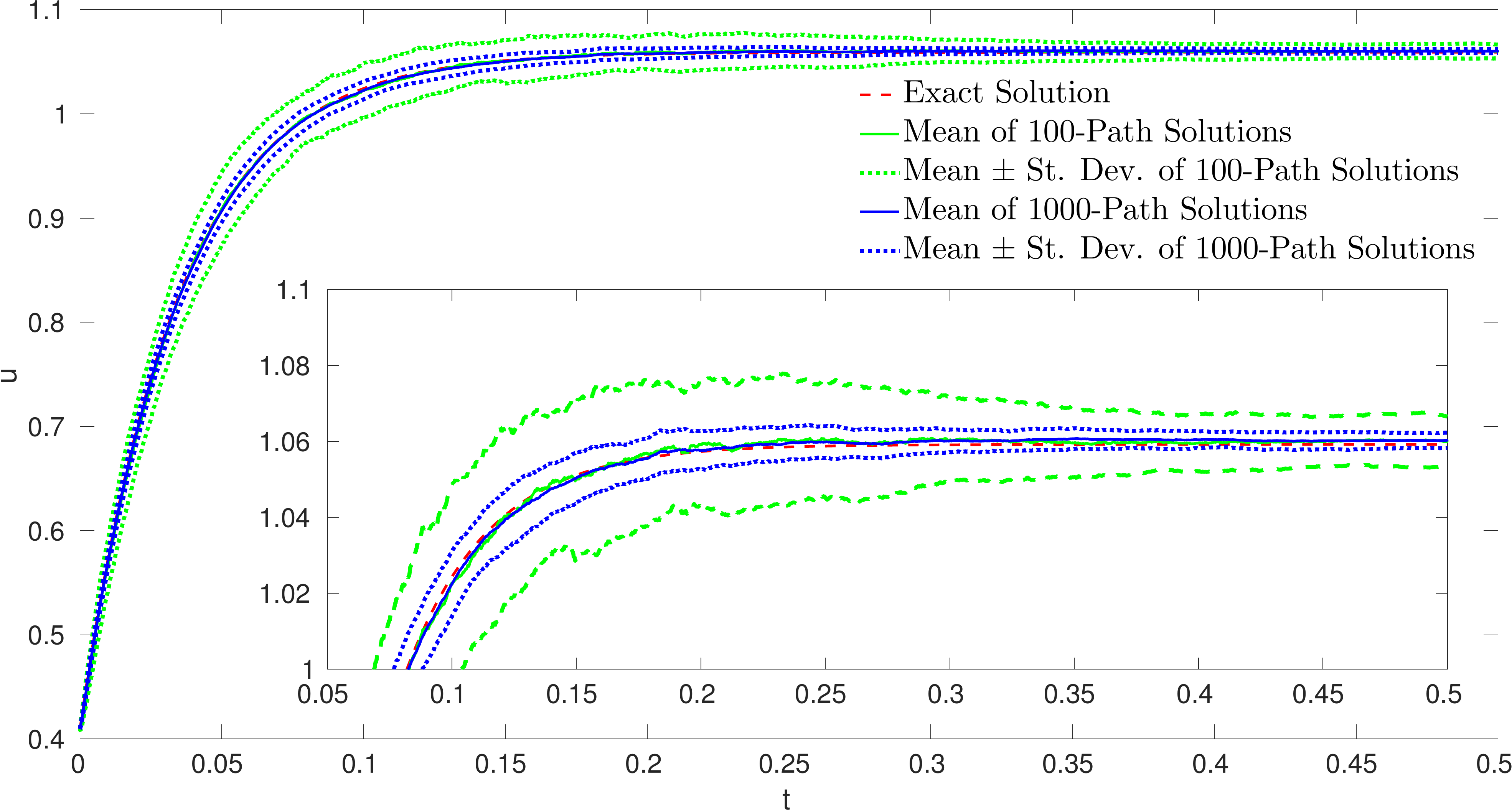}
   \label{fig:Ng1} 
\end{subfigure}

\begin{subfigure}[b]{1\textwidth}
   \includegraphics[width=1\linewidth]{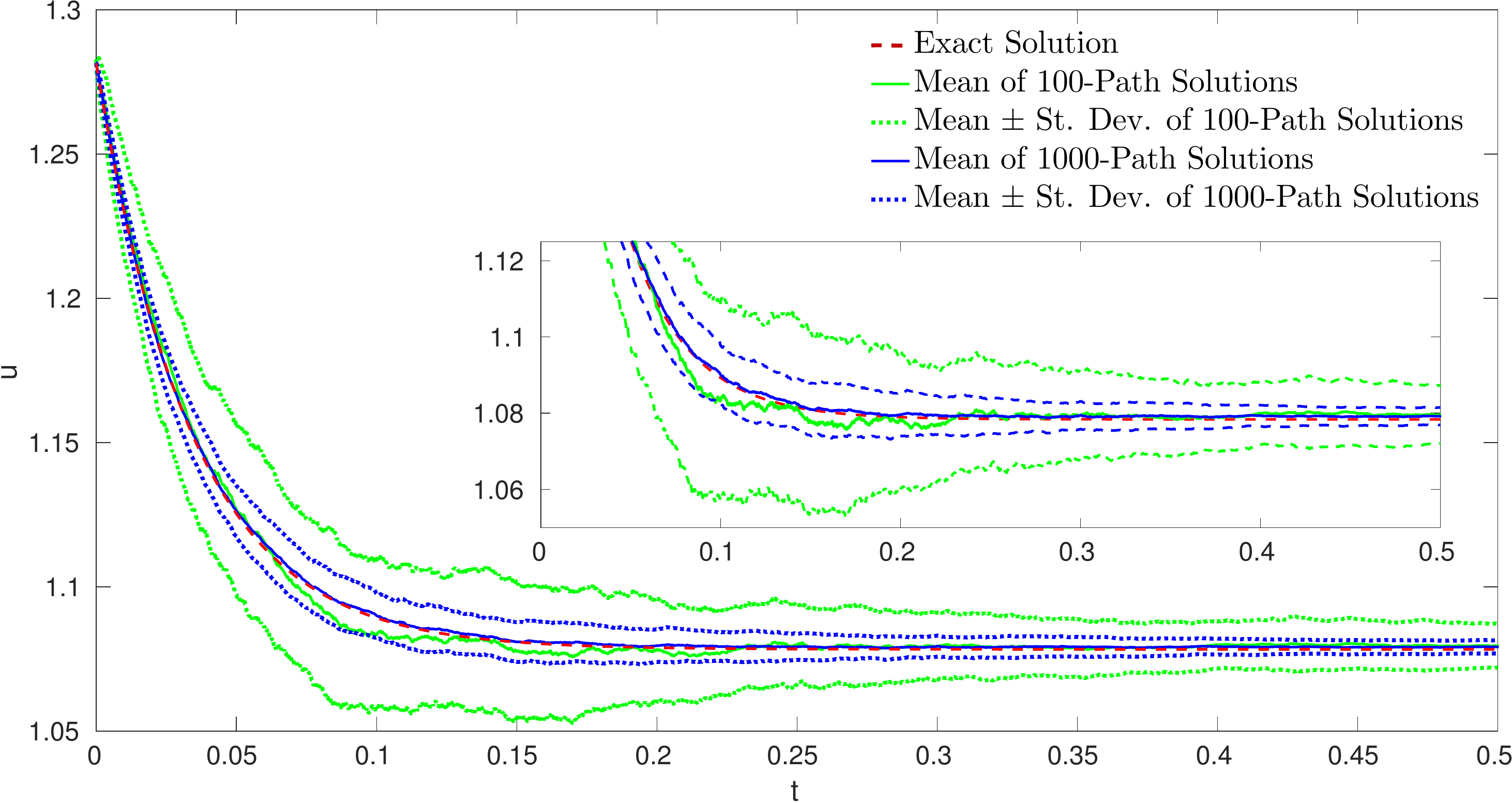}
   \label{fig:Ng2}
\end{subfigure}
\caption[3d_1o3]{
Convergence of stochastic solutions for the problem \eqref{3d_problem_1}.
\emph{Top, main:}
Convergence with respect to number of paths of the stochastic solution at (1/3, 2/3, 1/3), with $dt = 10^{-4}$.  
\emph{Top, inset:}
Zoomed-in view of the main plot. 
\emph{Bottom, main:}
Convergence with respect to number of paths of the stochastic solution at (3/5, 2/5, 3/5), with $dt = 10^{-4}$.  
\emph{Bottom, inset:}
Zoomed-in view of the main plot. 
}
\label{fig:3d_time_curves}
\end{figure}

Consider the three-dimensional unit cube $Q = [0,1]^3$. Define 
$g \equiv 1$ on $\partial Q$. We study the equation
\begin{align}
\label{3d_problem_1}
\begin{cases}
\partial_t u + (-\Delta_{Q, g})^{\alpha/2} u &= \sin(\pi x) \sin( \pi y) \sin(\pi z)\\
u \big|_{\partial Q} &= g \\
u(t = 0, x, y, z) &= 1 + \frac{1}{(\pi^2 + \pi^2  + \pi^2)^{\alpha/2}} 
\sin(\pi x) \sin(\pi y) \sin(\pi z) \\
&\qquad \qquad \qquad + \sin(2 \pi x) \sin(2 \pi y) \sin(2 \pi z).
\end{cases}
\end{align}
In this problem, 
\begin{itemize}
\item
the boundary condition (BC) function $g \equiv 1$,
\item
the right-hand side (RHS) is
$r(x,y,z) = \sin(\pi x) \sin(\pi y) \sin(\pi z)$, and 
\item
the initial condition (IC) is 
\begin{equation}
f(x,y,z) = 1 + \frac{1}{(\pi^2 + \pi^2 + \pi^2)^{\alpha/2}}  \sin(\pi x) \sin(\pi y) \sin(\pi z)
+ \sin(2 \pi x) \sin(2 \pi y) \sin(2 \pi z).
\end{equation} 
\end{itemize}
The exact solution is 
\begin{multline}
u(t,x,y,z) 
= 
1 + \frac{1}{(\pi^2 + \pi^2 + \pi^2)^{\alpha/2}} 
\sin(\pi x) \sin(\pi y) \sin(\pi z) \\
+
e^{-\left(4\pi^2+4\pi^2+4\pi^2\right)^{\alpha/2} t} 
\sin(2 \pi x) \sin(2 \pi y) \sin(2 \pi z).
\end{multline}
We fix $\alpha = \sqrt{2}$ and $dt = 10^{-4}$.  
We start by testing the parabolic solution formula in Figure \ref{fig:3d_time_curves}.
We approximate the solution $u(t,x,y,z)$ as a function of $t$ at two different points: $(x,y,z) = (1/3, 2/3, 1/3)$ and $(x,y,z) = (3/5, 2/5, 3/5)$. The mean of 100 stochastic solutions, computed using 100 or 1,000 paths (as indicated) starting from $(x,y,z)$ is plotted and compared to the exact solution. The``Mean $\pm$ Standard Deviation'' illustrates the expected variation of a stochastic solution computed using 100 or 1,000 paths, respectively, as it oscillates about the respective mean. Further, the means of the 100-path and 1000-path solutions \emph{also} represent stochastic solutions for $u(t,x,y,z)$ computed using $100 \times 100 = 10,000$ paths and $100 \times 1000 = 100,000$ paths, respectively. We see that these are well-converged to the true solution. 

\begin{figure}[p] 
  \label{fig7} 
  \begin{minipage}[b]{0.5\linewidth}
    \centering
    \includegraphics[width=\linewidth]{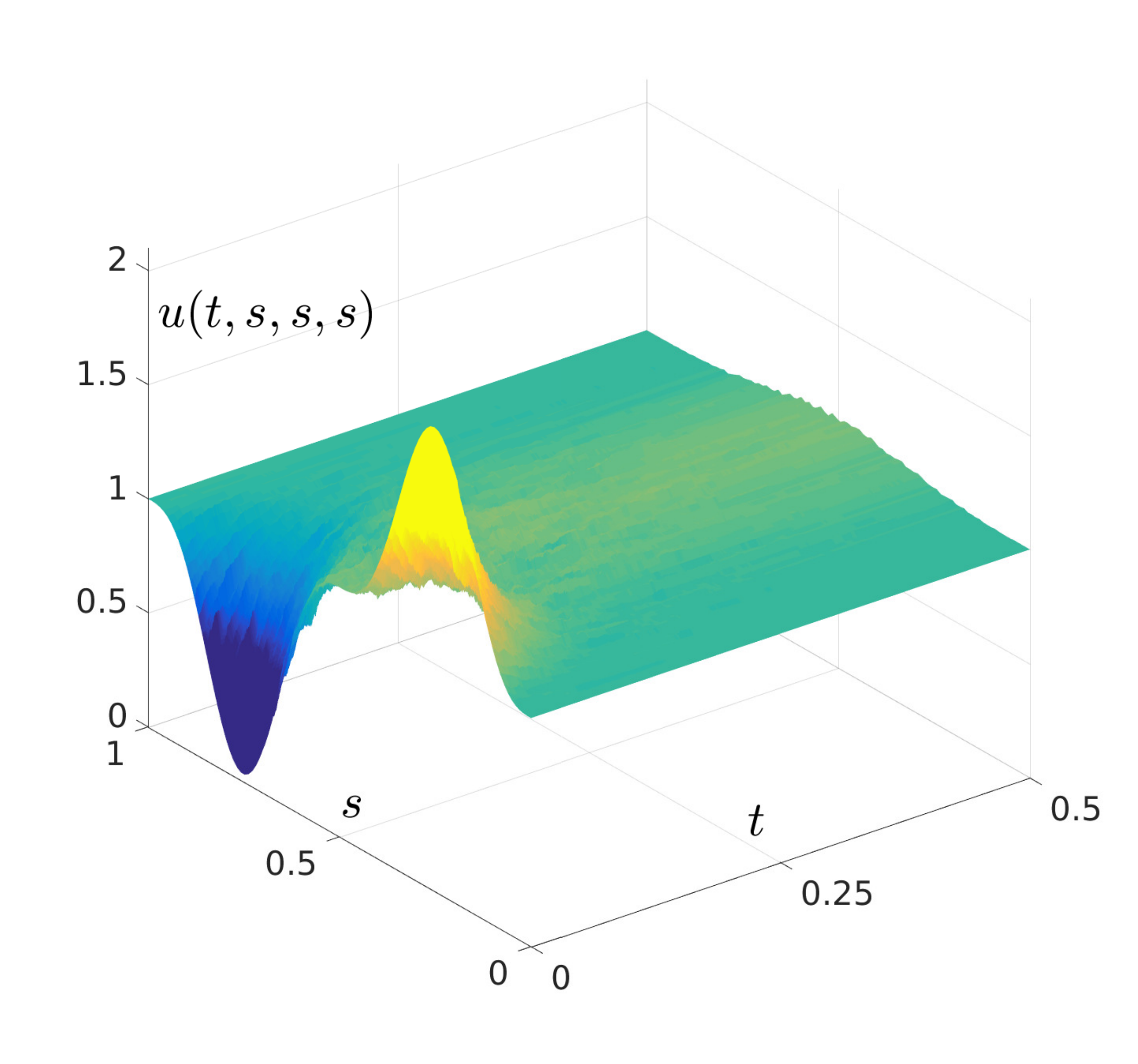} 
    \caption*{\small Stochastic solution computed using 100 paths for $t \in [0,0.5]$.} 
    \vspace{4ex}
  \end{minipage}
  \begin{minipage}[b]{0.5\linewidth}
    \centering
    \includegraphics[width=\linewidth]{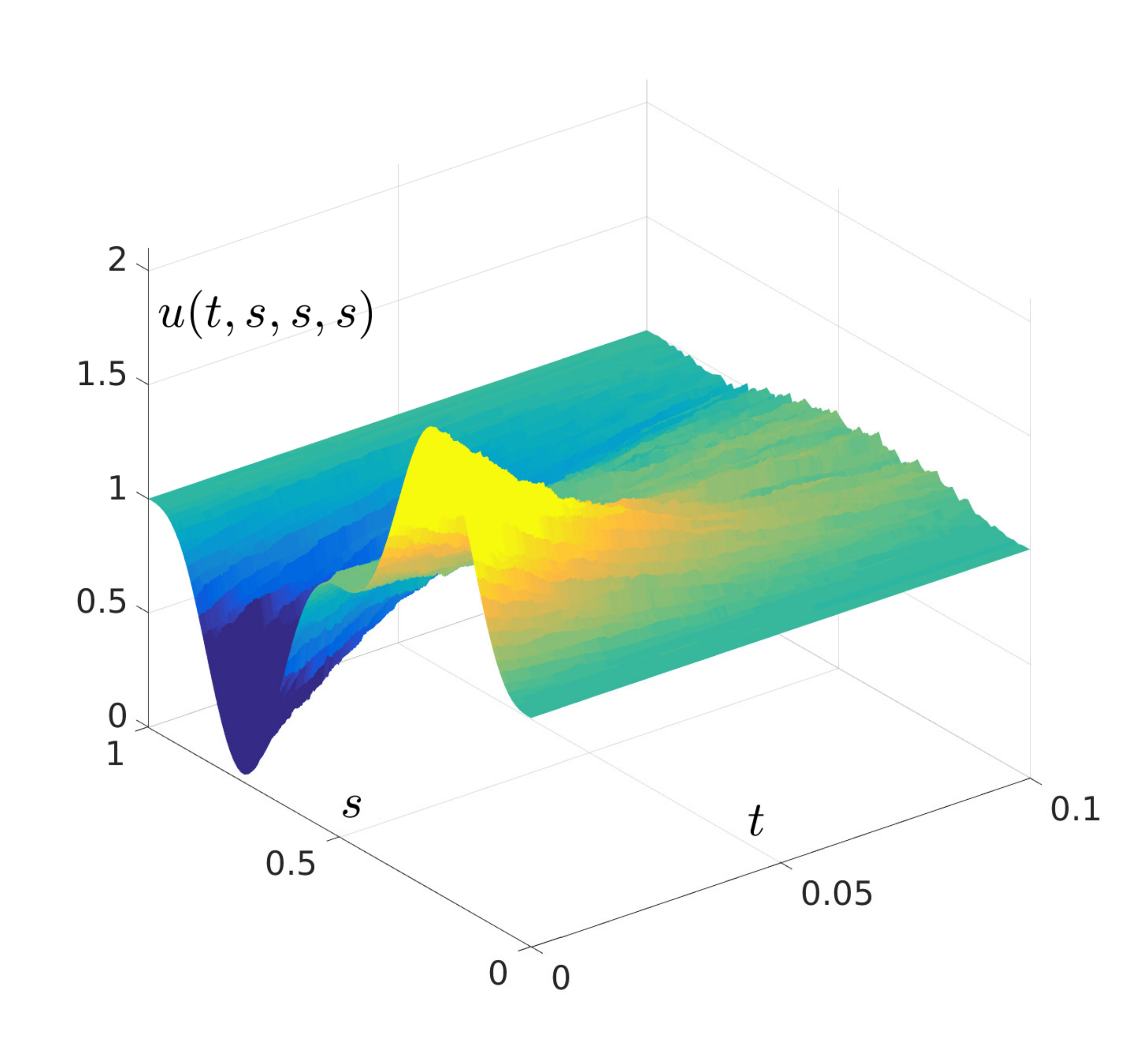} 
    \caption*{\small Stochastic solution computed using 100 paths zoomed in to a shorter time $t \in [0,0.1]$.} 
    \vspace{4ex}
  \end{minipage} 
  \begin{minipage}[b]{0.5\linewidth}
    \centering
    \includegraphics[width=\linewidth]{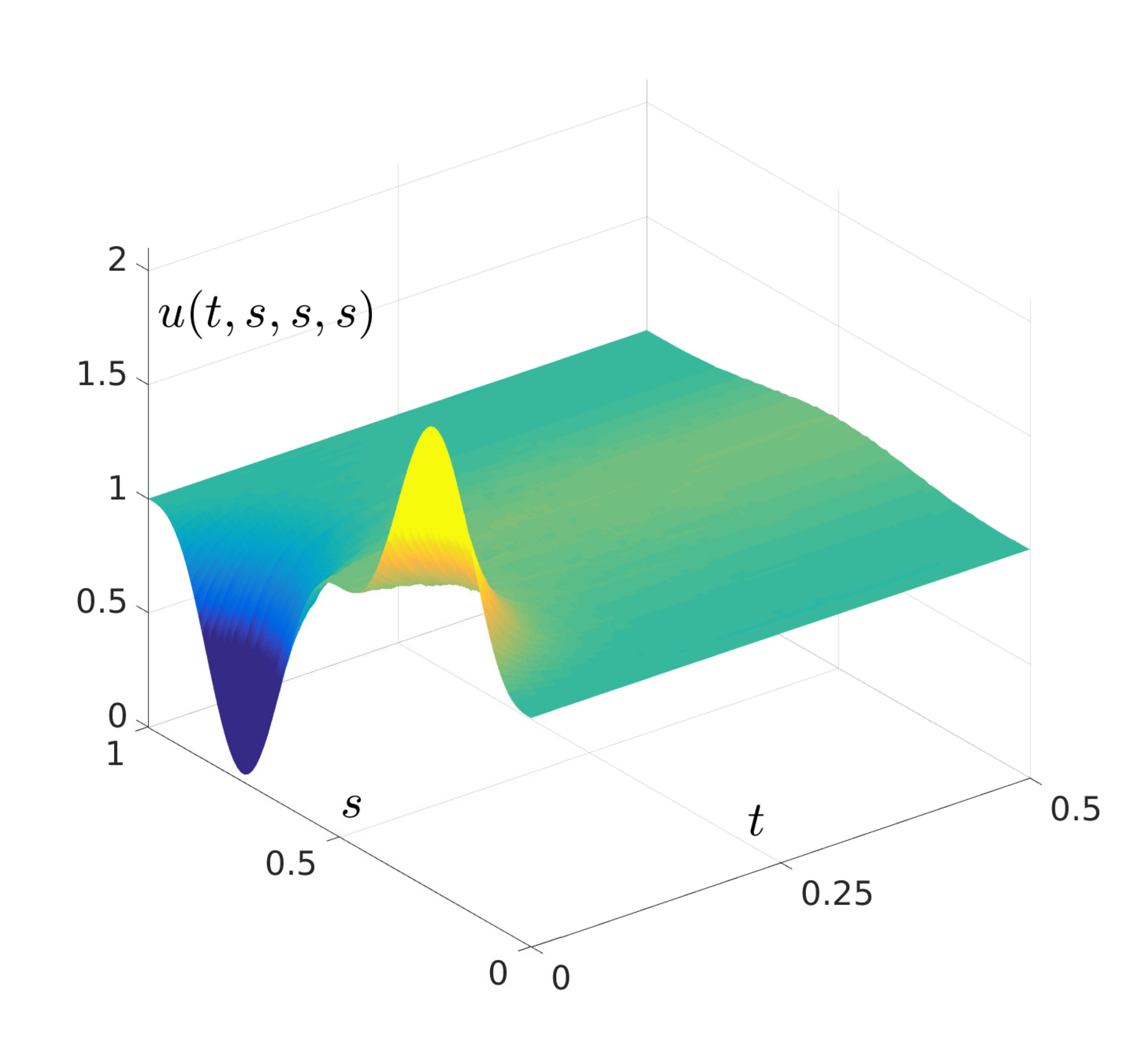} 
    \caption*{\small Stochastic solution computed using 1,000 paths for $t \in [0,0.5]$.} 
    \vspace{4ex}
  \end{minipage}
  \begin{minipage}[b]{0.5\linewidth}
    \centering
    \includegraphics[width=\linewidth]{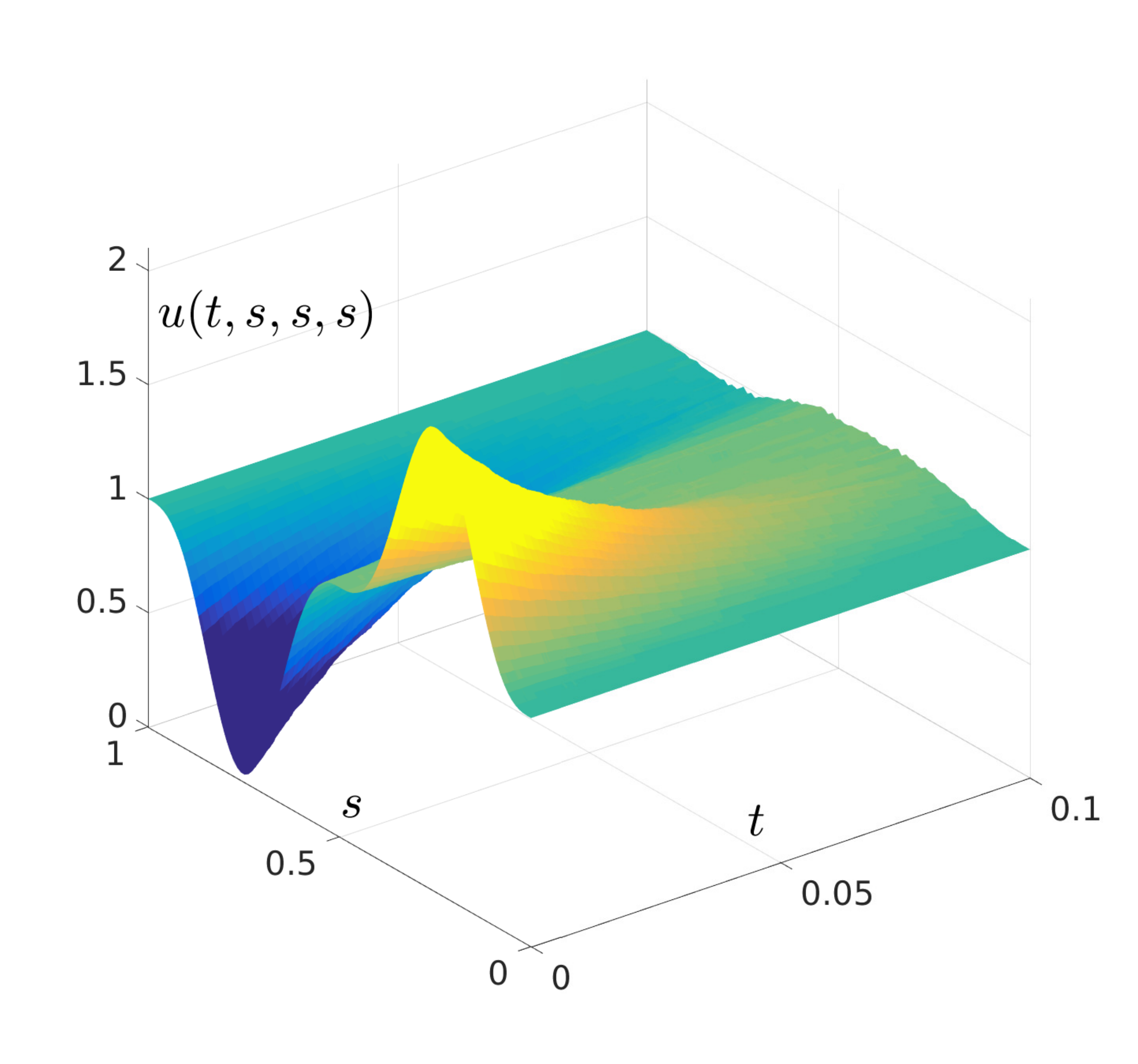} 
    \caption*{\small Stochastic solution computed using 1,000 paths zoomed in to a shorter time $t \in [0,0.1]$.}
    \vspace{4ex}
  \end{minipage} 
\caption{Stochastic solution using 100 paths (\emph{top}) and 1,000 paths  (\emph{bottom}) of problem  \eqref{3d_problem_1} in time $t$ on the diagonal $D = \{s,s,s\}$ in the unit cube.}
\label{fig:3d_spacetime_100_1000}
\end{figure}

\begin{figure}[p] 
  \label{fig7} 
  \begin{minipage}[b]{0.5\linewidth}
    \centering
    \includegraphics[width=\linewidth]{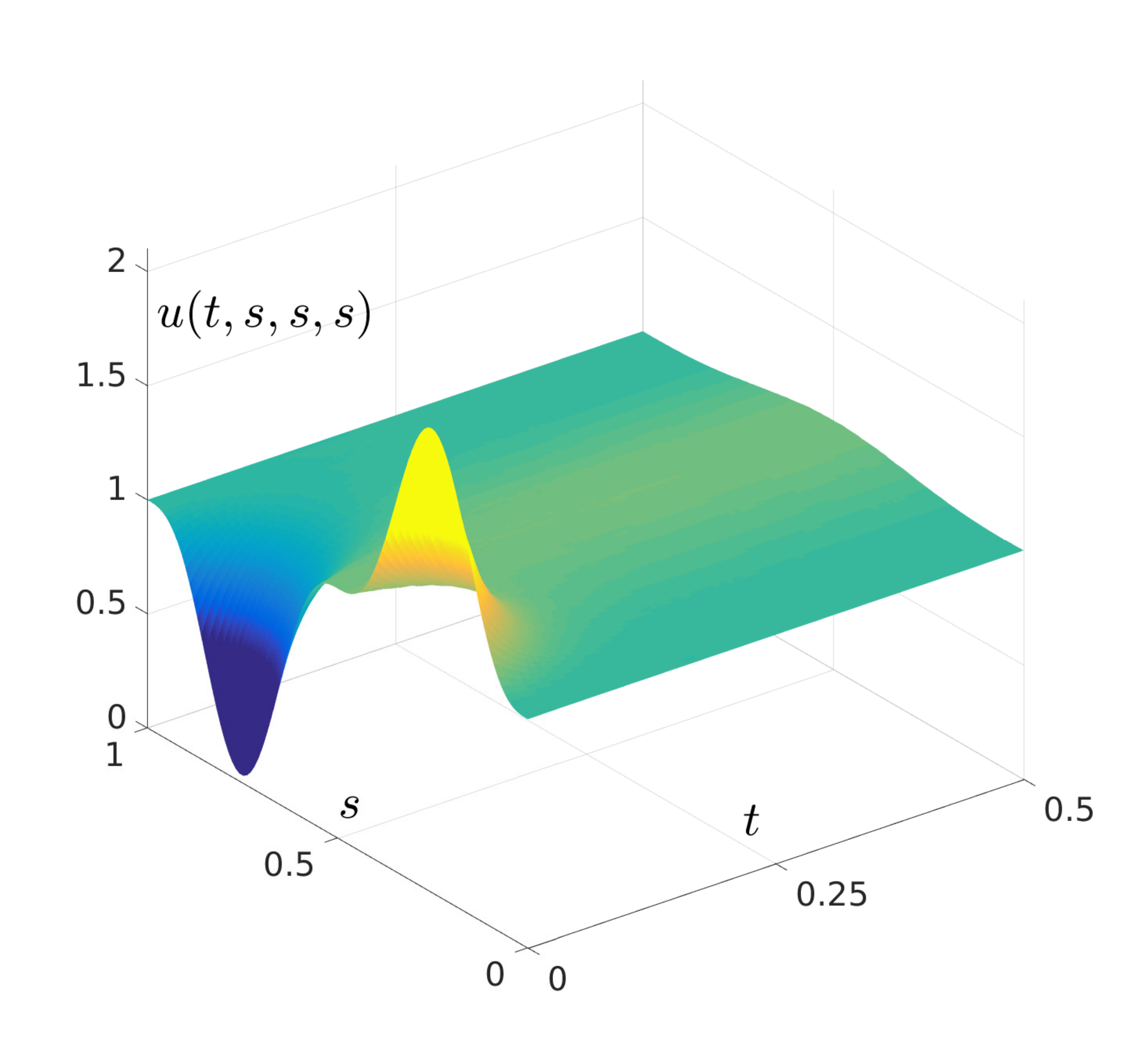} 
    \caption*{\small Stochastic solution computed using 10,000 paths for $t \in [0,0.5]$.} 
    \vspace{4ex}
  \end{minipage}
  \begin{minipage}[b]{0.5\linewidth}
    \centering
    \includegraphics[width=\linewidth]{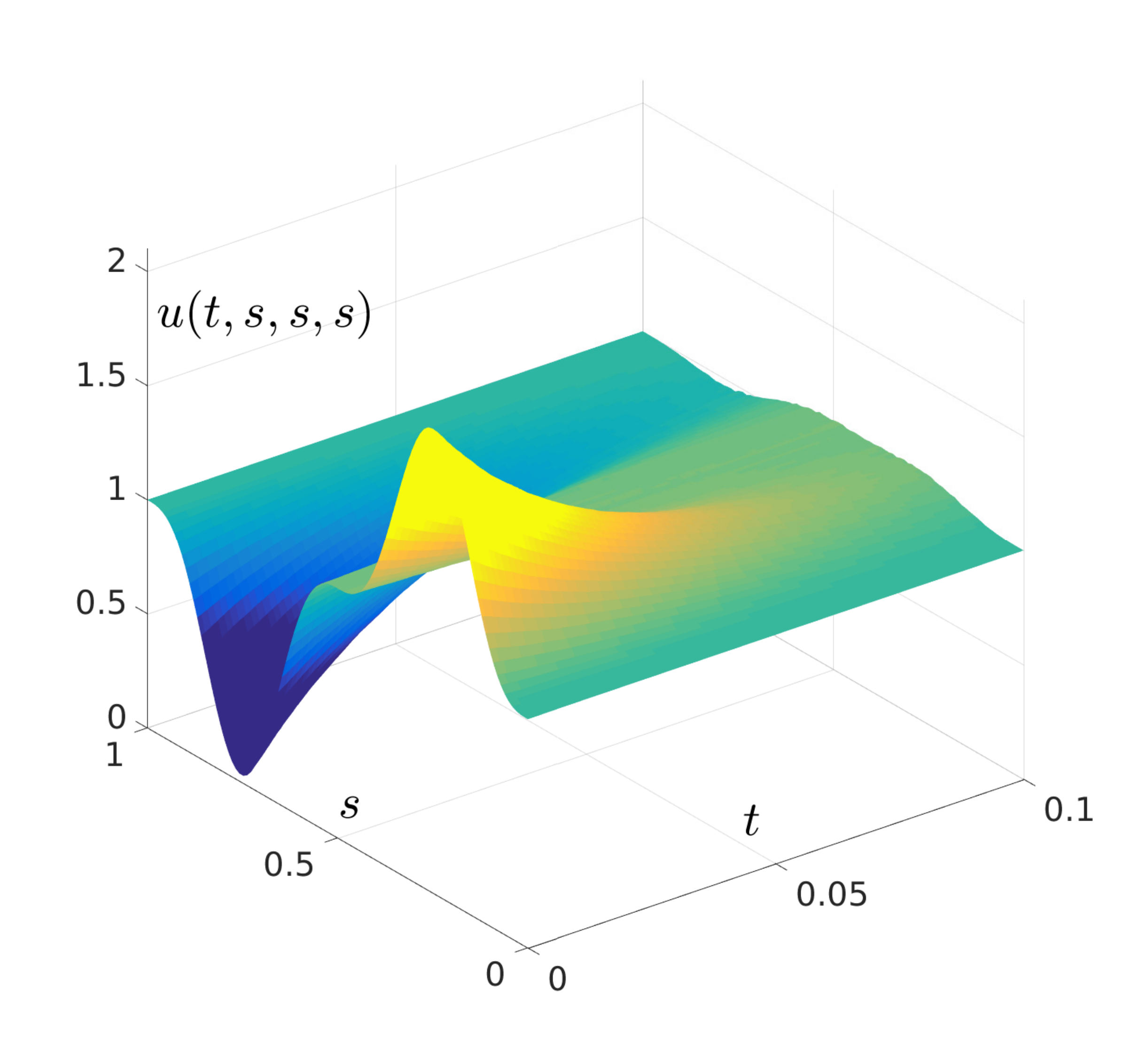} 
    \caption*{\small Stochastic solution computed using 10,000 paths zoomed in to a shorter time $t \in [0,0.1]$.} 
    \vspace{4ex}
  \end{minipage} 
  \begin{minipage}[b]{0.5\linewidth}
    \centering
    \includegraphics[width=\linewidth]{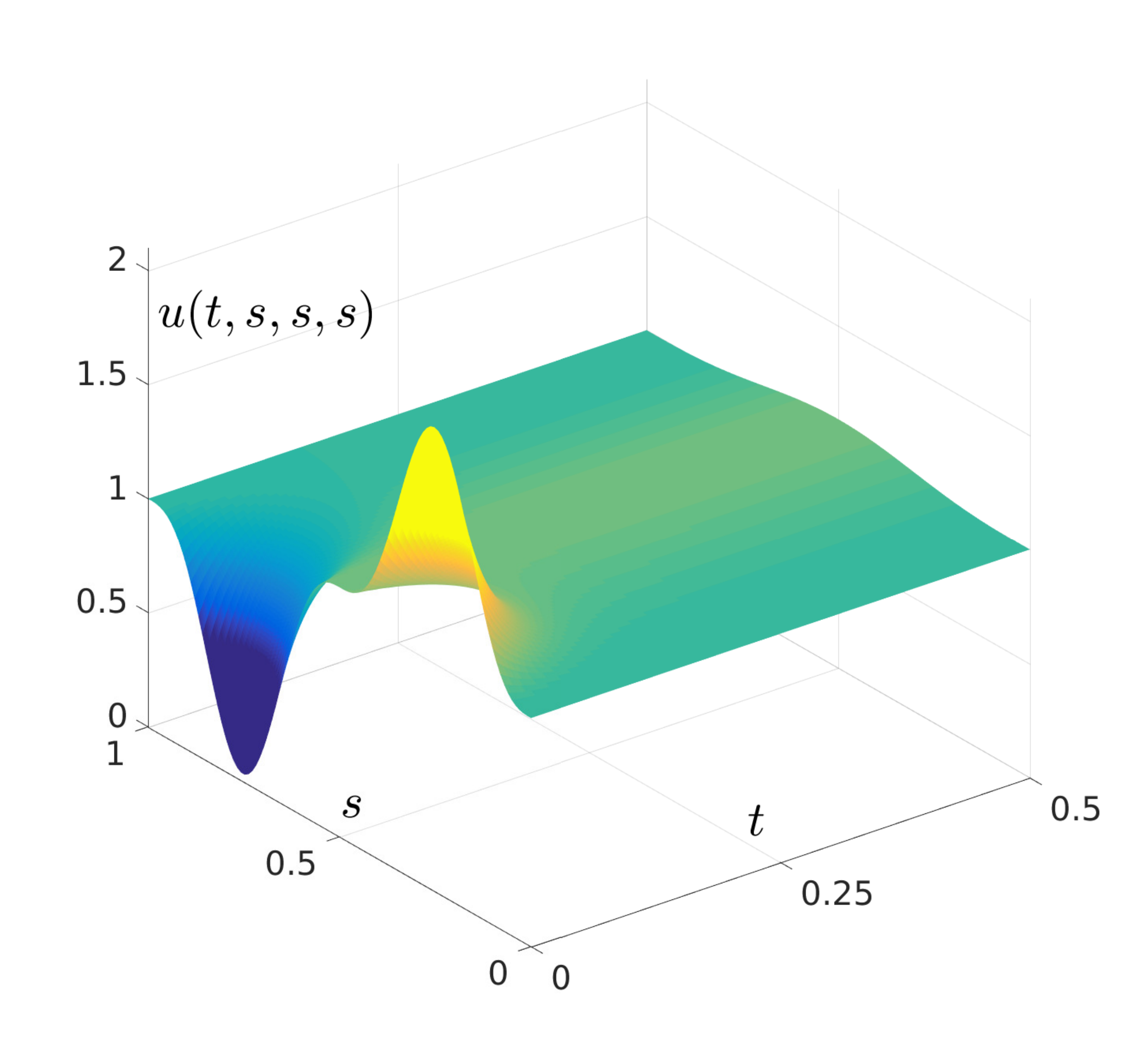} 
    \caption*{\small Exact solution for $t \in [0,0.5]$.} 
    \vspace{4ex}
  \end{minipage}
  \begin{minipage}[b]{0.5\linewidth}
    \centering
    \includegraphics[width=\linewidth]{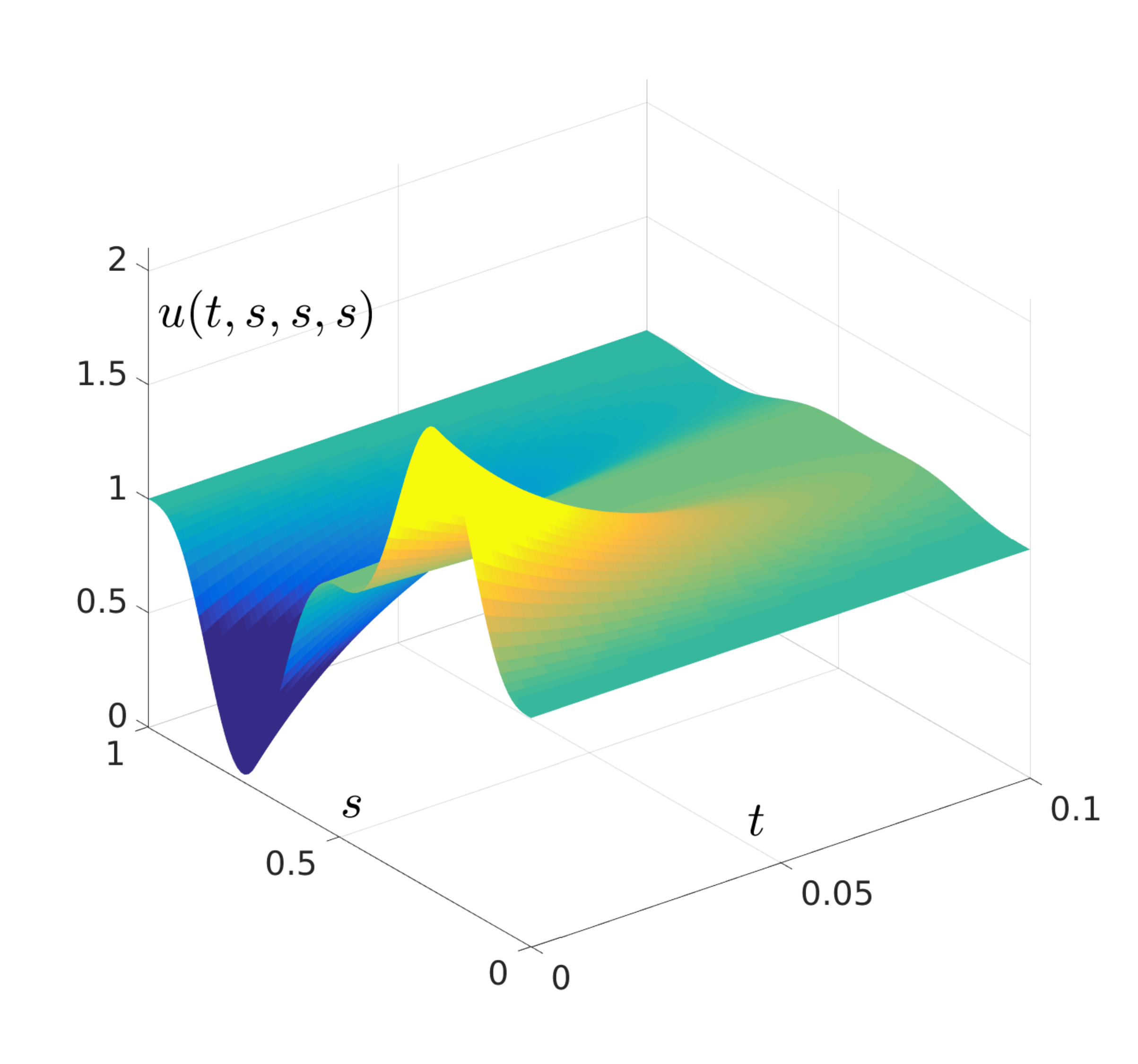} 
    \caption*{\small Exact solution zoomed in to a shorter time $t \in [0,0.1]$.} 
    \vspace{4ex}
  \end{minipage} 
\caption{Stochastic solution (using 10,000 paths, \emph{top}) and exact solution (\emph{bottom}) of problem  \eqref{3d_problem_1} in time $t$ on the diagonal $D = \{s,s,s\}$ in the unit cube.}
\label{fig:3d_spacetime_10000_exact}
\end{figure}

\begin{figure}[p] 
  \label{fig7} 
  \begin{minipage}[b]{0.5\linewidth}
    \centering
    \includegraphics[width=\linewidth]{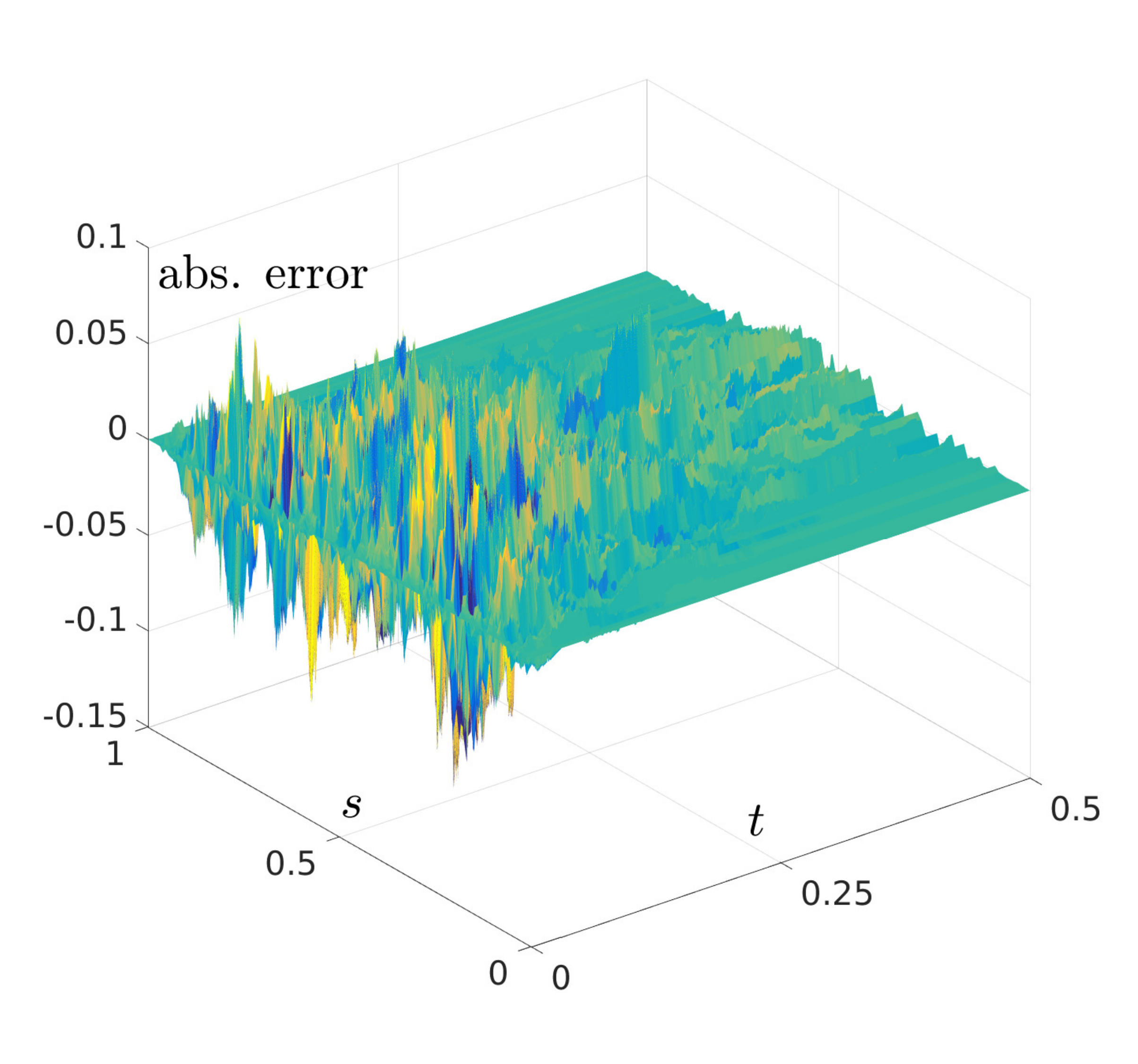} 
    \caption*{\small 100 paths for the stochastic solution.} 
    \vspace{4ex}
  \end{minipage}
  \begin{minipage}[b]{0.5\linewidth}
    \centering
    \includegraphics[width=\linewidth]{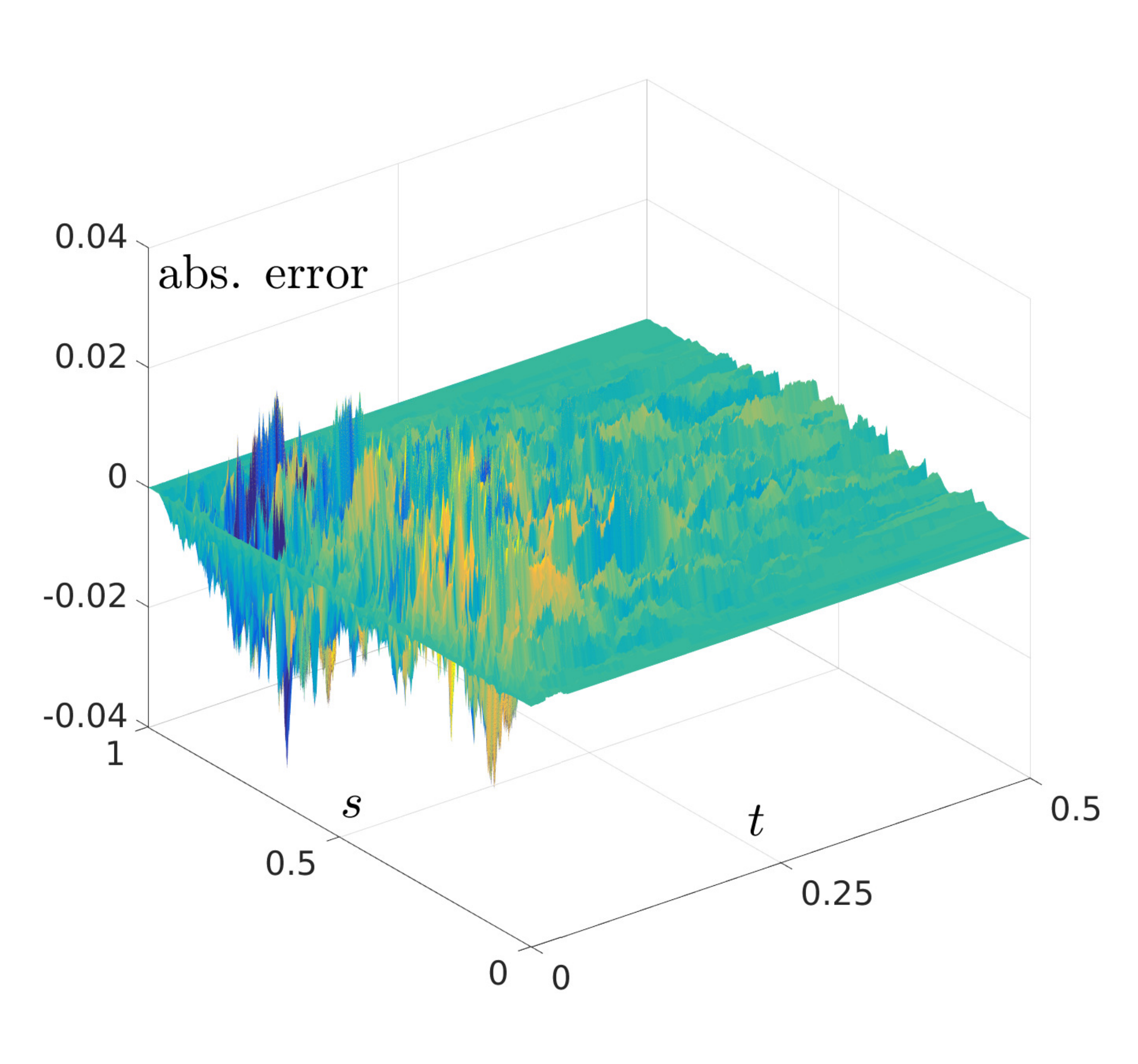} 
    \caption*{\small  1,000 paths for the stochastic solution.} 
    \vspace{4ex}
  \end{minipage} 
  \begin{minipage}[b]{0.5\linewidth}
    \centering
    \includegraphics[width=\linewidth]{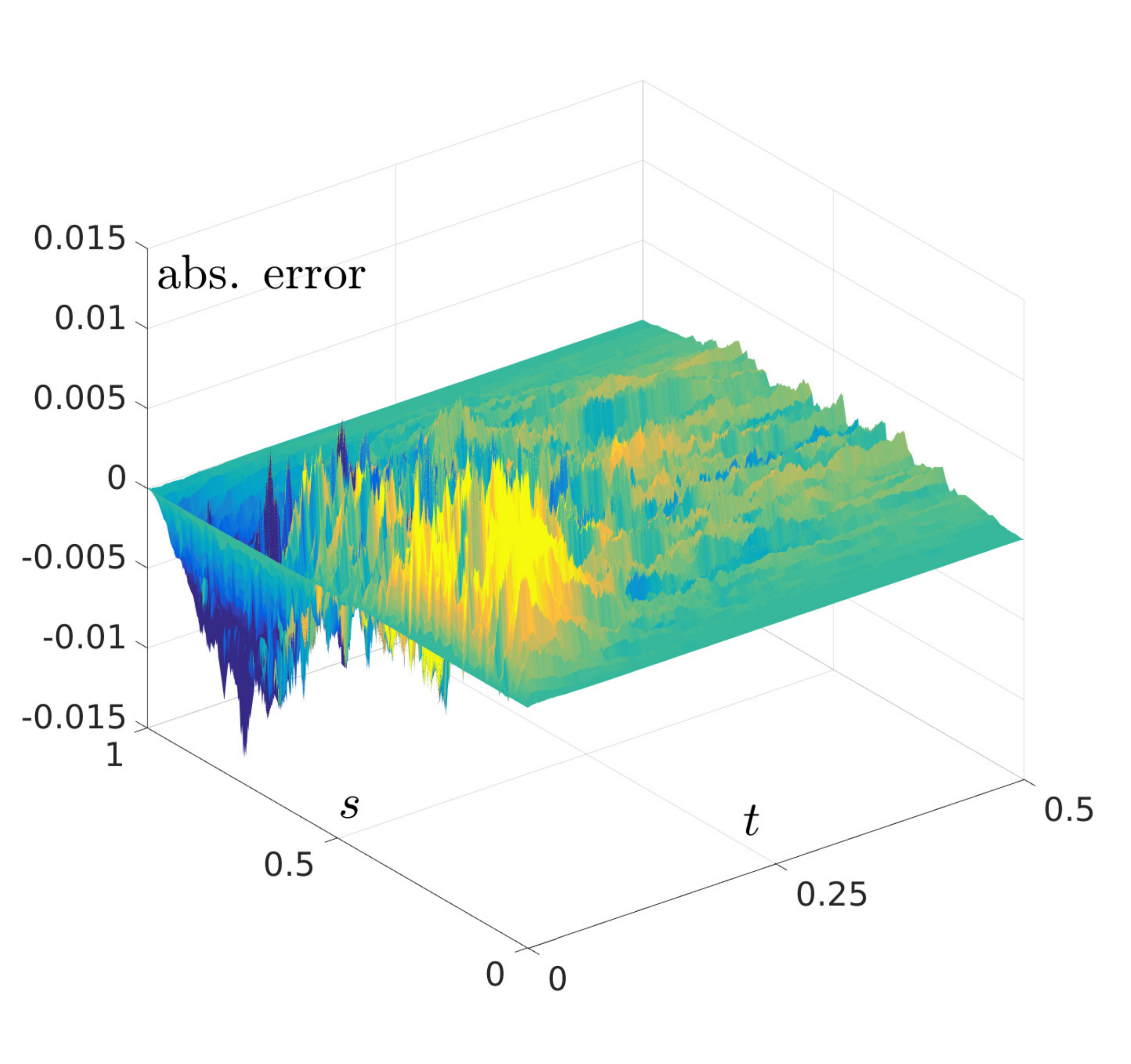} 
    \caption*{\small 10,000 paths for the stochastic solution.} 
    \vspace{4ex}
  \end{minipage}
\caption{Convergence of the absolute error $u_{\text{stochastic}}(t,s,s,s) - u_{\text{exact}}(t,s,s,s)$ for the solution of problem \eqref{3d_problem_1} in spacetime for time $t \in [0, 0.5]$ and position $s \in [0,1]$ along the diagonal $(s,s,s)$. 100, 1,000, and 10,000 paths are used in the stochastic solution.}
\label{fig:3d_spacetime_abs_errors}
\end{figure}

Next, we study the stochastic solution in spacetime along line diagonal line $D = \{(s,s,s) \ \big| \ s \in [0,1]\} \subset R$. We take 100 equispaced points on $D$, and compute the solution at each of the points using 100, 1,000, and 10,000 paths. In Figures \ref{fig:3d_spacetime_100_1000} and \ref{fig:3d_spacetime_10000_exact}, we show surface plots in $(s,t)$ of these solutions and compare to the exact solution; in Figure \ref{fig:3d_spacetime_abs_errors}, we plot the absolute error $u_{\text{stochastic}}(t,s,s,s) - u_{\text{exact}}(t,s,s,s)$ as a surface plot in $(s,t)$. We observe a decrease in the maximum value of the absolute error by a factor of $\sqrt{10}$ each time the number of paths is increased by a factor of $10$. 

\newpage
\FloatBarrier

\begin{figure}[htpb]
   \includegraphics[width=1\linewidth]{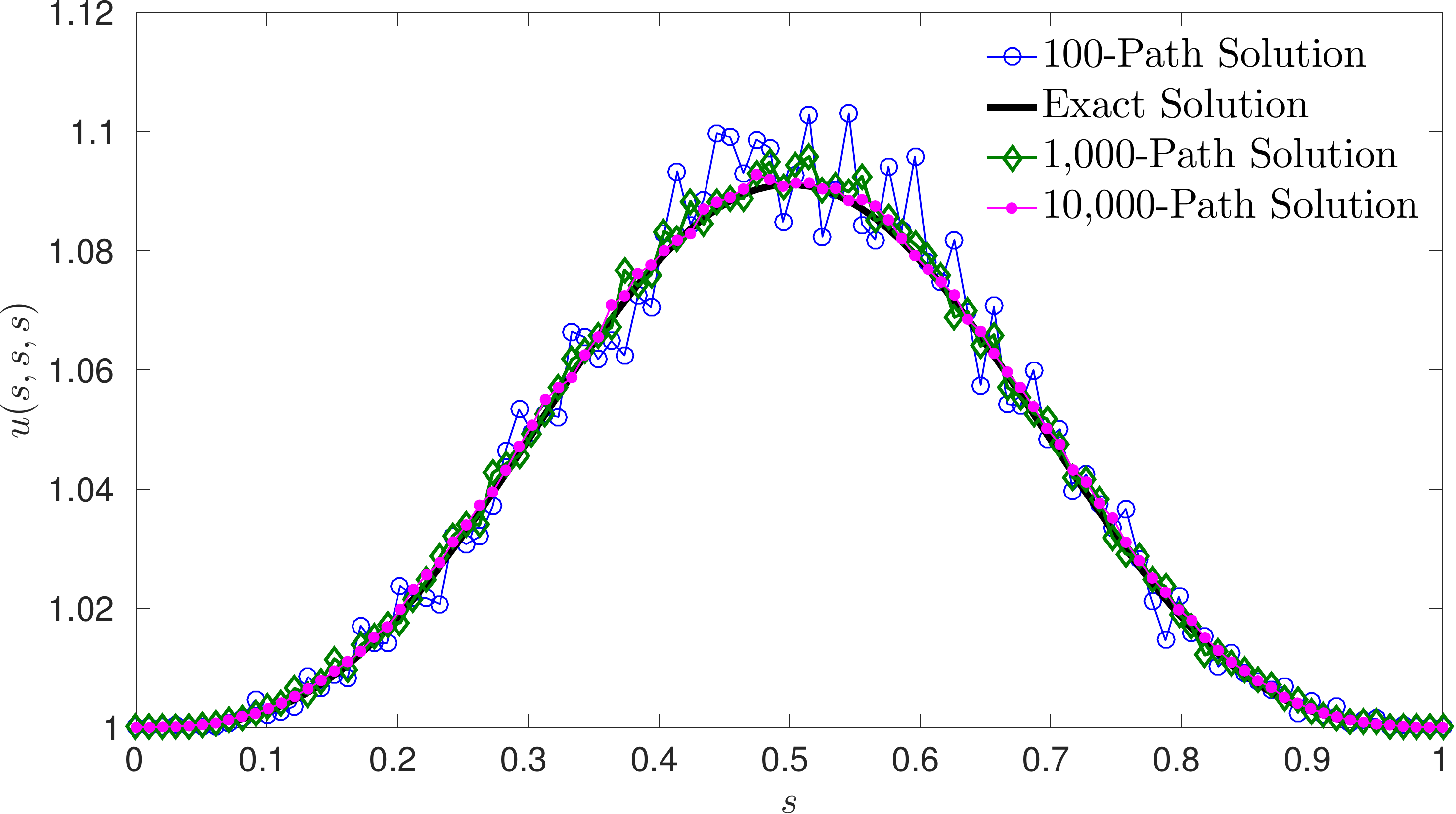}
   \caption{Stochastic solution of problem \eqref{3d_elliptic} in the 3D unit cube at 100 equispaced points along the line $D = \{(s,s,s) \ \big| \ s \in [0,1]\}$, using 100, 1,000, and 10,000 paths starting from each point.}
\label{fig:3d_elliptic}
\end{figure}

Next, we test the elliptic solution formula for the example consisting of the steady state to problem \eqref{3d_problem_1}, i.e., 
\begin{align}
\label{3d_elliptic}
\begin{cases}
(-\Delta_{Q, g})^{\alpha/2} u &= \sin(\pi x) \sin( \pi y) \sin(\pi z)\\
u \big|_{\partial Q} &= 1
\end{cases}
\end{align}
which has exact solution
\begin{equation}
u(x,y,z) 
= 
1 + \frac{1}{(\pi^2 + \pi^2 + \pi^2)^{\alpha/2}} 
\sin(\pi x) \sin(\pi y) \sin(\pi z).
\end{equation}
We keep $\alpha = \sqrt{2}$ and $dt = 10^{-4}$, and illustrate convergence by computing the solution $u(s,s,s)$ using 100, 1,000, and 10,000 paths from each of the 100 equispaced points on $D = \{(s,s,s) \ \big| \ s \in [0,1]\}$. This is shown in Figure \ref{fig:3d_elliptic}. In fact, the same exact sets of stopped paths are used to produce Figures \ref{fig:3d_spacetime_100_1000} and \ref{fig:3d_spacetime_10000_exact} are used to evaluate the stochastic solution formula for this elliptic problem.

\newpage
\FloatBarrier

\section{Conclusion and Future Work.}
\label{conclusion}
We have proven, implemented, and verified stochastic solution (Feynman-Kac) formulas for Cauchy and Dirichlet problems for the spectral fractional Laplacian $(-\Delta_{\Omega,g})^{\alpha/2}$ with nonzero Dirichlet boundary conditions. This operator was introduced recently in \cite{AntilPfeffererRogovs} and \cite{Cusimano2017}, and our article represents a novel probabilistic approach to this topic. 
The formulas, which involve subordinate stopped Brownian motion, were verified by considering a number of benchmark examples in two and three dimensions, and convergence with respect to the number of paths and time step parameter $dt$ was studied. This work validates the proposed operator $(-\Delta_{\Omega,g})^{\alpha/2}$ from a stochastic perspective. 

Stochastic solution formulas provide an attractive method to both understand fractional operators and compute solutions to associated boundary value problems. 
Such formulas for Neumann boundary value conditions, which involve boundary local time, are a worth exploring in this regard \cite{brosamler1976probabilistic, bass1991some, lions1984stochastic, bencherif2009probabilistic}. Further numerical studies may explore Monte Carlo solution of boundary value problems in complex domains. The direct SDE discretization discussed in this article is an alternative to the work of \cite{song2017computing}, as implementation requires only an efficient subroutine to determine if the process is in the domain or not. Moreover, Monte Carlo methods based on such formulas can provide efficient solutions in high dimensions, a fertile area for applications. Walk-on-spheres \cite{muller1956, cai2013, cai2017, Zhou2016} and quasi-Monte Carlo \cite{moskowitz1996smoothness} approaches may be explored to accelerate such solution methods. 

\section{Acknowledgements.}
We thank George Em Karniadakis of Brown University and Wei Cai of Southern Methodist University for their guidance since the beginning of this project.  
We also thank Mark M. Meerschaert of Michigan State University for many very helpful discussions on both the theory and the numerical implementation of this work.
Both authors acknowledge funding from MURI/ARO grant W911NF-15-1-0562.
M.G. acknowleges support from an NSF Graduate Research Fellowship and a Brown University Deans' Faculty Fellowship.

\appendix
\section{Facts about Eigenvalue Decomposition and Sobolev Spaces.}
\label{appendix}
The following results about the spectrum of $-\Delta$ and eigenfunction decomposition are transcribed from Chapter 6.1 of \cite{larsson2008partial}. 
In that text, $\Omega \subset \mathbb{R}^d$ is taken to be a domain with smooth boundary, although in general $\Omega$ can have piecewise smooth boundary \cite{courant}. Consider the Dirichlet eigenvalue problem
\begin{equation}
\label{eigenvalue_problem}
-\Delta e_\lambda = \lambda e_\lambda \text{\quad in $\Omega$}, \quad e_\lambda|_{\partial\Omega} = 0, 
\end{equation}
defining eigenvalues $\lambda$ with corresponding eigenfunctions $e_\lambda$.
\begin{theorem}
The eigenvalues $\lambda$ of \eqref{eigenvalue_problem} are real, positive, and can be ordered in a nondecreasing sequence $\lambda_k$ such that 
\begin{equation}
\lambda_k \rightarrow \infty \quad\text{as}\quad k \rightarrow \infty.
\end{equation}
Two eigenfunctions $e_{\lambda_1}$ and $e_{\lambda_2}$ corresponding to different eigenvalues $\lambda_1 \neq \lambda_2$, are orthogonal in $L^2$ and $H^1_0$. 
\end{theorem}

\begin{theorem}
\label{H10_theorem}
The sequence of eigenfunctions $\left\{e_{k}\right\}$ corresponding to the nondecreasing sequence of eigenvalues $\left\{\lambda_k\right\}$ forms an orthonormal basis for $L_2(\Omega)$. Moreover, $v \in H^1_0(\Omega)$ if and only if the series 
\begin{equation}
\sum_{k = 0}^\infty \lambda_k \left( v, e_k \right)^2 \quad \text{converges}.
\end{equation}
The following identity holds:
\begin{equation}
[v]^2_{H^1_0} = \| \nabla v \|^2_{L^2} = \sum_{k = 1}^\infty \lambda_k \left( v, e_k \right)^2.
\end{equation}
\end{theorem}

\noindent
The Sobolev inequality that we use is stated in the appendix of \cite{larsson2008partial}:
\begin{theorem}
\label{standard_sobolev_inequality}
Let $k > d/2 + r$. There exists a constant $C$ such that
\begin{equation}
\| v \|_{C^r(\Omega)} \le C \|v\|_{H^k(\Omega)}.
\end{equation}
Thus, ${H^k(\Omega)} \subset {C^r(\Omega)}$.
\end{theorem}
\noindent
This holds for domains $\Omega$ with Lipschitz boundary \cite{aubin} or which satisfy the cone condition \cite{adams}.
Finally, we require the following elliptic regularity estimate, transcribed from Chapter 3.7 of \cite{larsson2008partial}:
\begin{theorem}
\label{elliptic_regularity}
Let $\Omega$ be a smooth domain.
Then for $u \in H^{k+2}(\Omega) \cap H^1_0(\Omega) $, there exists a constant $C$ such that
\begin{equation}
\| u \|_{H^{k+2}(\Omega)} \le C \| (-\Delta) u \|_{H^k(\Omega)}. 
\end{equation}
\end{theorem}

\bibliographystyle{unsrtnat}
\bibliography{feynman_kac_ref.bib}

\end{document}